\documentclass[10pt]{amsart}
\usepackage[margin=1.2in]{geometry}
\usepackage{colortbl}
\usepackage{color}
\usepackage{cite}
\newtheorem{theorem}{Theorem}[section]

\theoremstyle{definition}

\theoremstyle{remark}
\newtheorem{remark}[theorem]{Remark}

\numberwithin{equation}{section}

\usepackage{graphicx}
\usepackage{subcaption}
\usepackage{multirow}
\usepackage{bbm}
\usepackage{lipsum}
\usepackage{amsfonts}
\usepackage{graphicx}
\usepackage{epstopdf}
\usepackage{algorithmic}
\usepackage{subcaption}
\usepackage{amssymb}
\usepackage{mathtools}
\usepackage{bbm}

\begin{document}
	
	\title{Conservative semi-Lagrangian schemes for kinetic equations \\Part II: Applications}
	
	\author[S. Y. Cho]{Seung Yeon Cho}
	\address{Seung Yeon Cho\\
		Department of Mathematics and Computer Science  \\
		University of Catania\\
		95125 Catania, Italy}
	\email{chosy89@skku.edu}
	
	\author[S. Boscarino]{Sebastiano Boscarino}
	\address{Sebastiano Boscarino\\
		Department of Mathematics and Computer Science  \\
		University of Catania\\
		95125 Catania, Italy} \email{boscarino@dmi.unict.it}

	\author[G. Russo]{Giovanni Russo}
	\address{Giovanni Russo\\
		Department of Mathematics and Computer Science  \\
		University of Catania\\
		95125 Catania, Italy} \email{russo@dmi.unict.it}
	
	\author[S.-B. Yun]{Seok-Bae Yun}
	\address{Seok-Bae Yun\\
		Department of Mathematics\\ 
		Sungkyunkwan University\\
		Suwon 440-746, Republic of Korea}
	\email{sbyun01@skku.edu}

	\maketitle
	
	\begin{abstract}
		In this paper, we present a new class of conservative semi-Lagrangian schemes for kinetic equations. They are based on the conservative reconstruction technique introduced in \cite{BCRY1}. The methods are high order accurate both in space and time. Because of the semi-Lagrangian nature, the time step is not restricted by a CFL-type condition. Applications are presented to rigid body rotation, Vlasov Poisson system and the BGK model of rarefied gas dynamics. In the first two cases operator splitting is adopted, to obtain high order accuracy in time, and a conservative reconstruction that preserves the maximum and minimum of the function is used. For initially positive solutions, in particular, this guarantees exact preservation of the $L^1$-norm. Conservative schemes for the BGK 
		model are constructed by coupling the conservative reconstruction with a conservative treatment of the collision term. High order in time is obtained by either Runge-Kutta or BDF time discretization of the equation along characteristics. Because of $L$-stability and exact conservation, the resulting scheme for the BGK model is {\em asymptotic preserving\/} for the underlying fluid dynamic limit. Several test cases in one and two space dimensions confirm the accuracy and robustness of the methods, and the AP property of the schemes when applied to the BGK model.
	\end{abstract}

\section{Introduction}
In the first part of this paper \cite{BCRY1}, we introduce and analyse a new conservative reconstruction, which is now adopted to build conservative semi-Lagrangian (SL) schemes for various equations, including rigid rotation and kinetic equations. The basic idea of the semi-Lagrangian approach is to integrate the equations along the characteristics and update solutions on a fixed grid. This allows the use of larger time steps with respect to Eulerian methods, which suffer from CFL-type time step restrictions. 

For this reason, SL methods are suitable for the numerical treatment of kinetic equations, where the collision term may be expensive, and they have been extensively adopted in the literature in the context of Vlasov-Poisson model \cite{CMS-vlasov-2010,grandgirard2006drift,qiu2017high} or BGK-type models \cite{GRS,xu2010unified,dimarco2015towards}.

Since the aforementioned kinetic equations have various conservative quantities, the conservation property is very relevant in the numerical treatment of such models. In particular, it could be crucial to maintain conservation for the long time simulation of the Vlasov-Poisson equation or capturing shocks correctly in the 
the fluid limit of the BGK model.

On the other hand, the characteristic tracing structure of SL method necessarily requires to compute the numerical solutions on off-gird points. Hence, it is important to take into account a reconstruction technique which gives high order accurate and conservative solutions. Also, it may be relevant to preserve other properties of the solution. For example, for the initial value problem for the Vlasov-Poisson system, the distribution function is carried unchanged along particle trajectories, therefore the range of the solution (the interval between minimum and maximum of the initial condition) is maintained, so it would be nice to prevent the numerical solution from exceeding such a range at a discrete level. In other cases, e.g.\ for piecewise smooth solutions with jump discontinuities for the BGK model, it is important to prevent the appearance of spurious oscillations. 
All these properties crucially depend on the details of the reconstruction, not only on its conservation property.

High order non-oscillatory WENO type interpolations, such as the one introduced in \cite{CFR}, and that we call {\em generalized WENO} (GWENO), or the one used in \cite{carrillo2007nonoscillatory}, can be adopted, for example, to prevent the appearing of spurious oscillations. They have been applied to develop a class of non-splitting SL methods \cite{GRS,RS,SP} developed for the BGK model \cite{BGK} of the celebrated Boltzmann equation \cite{Cercignani} in the rarefied gas-dynamics. However, when numerical solutions contain discontinuities, such schemes may not be conservative due to the non-linear weights used for relieving oscillations \cite{BCRY}.


As a remedy for the lack of conservation, several approaches introduce high-order conservative and non-oscillatory schemes \cite{qiu2010conservative,QS}. In \cite{FSB-vlasov-2001} the authors propose a slope correction to prevent negativity, which is further improved in \cite{umeda2012non}.

In this paper, we use a technique of conservative reconstructions introduced in \cite{BCRY1}. 
In particular, a maximum principle preserving (MPP) limiter \cite{zhang2010maximum,zhang2011maximum} is adopted in order to preserve the range of the solution, while CWENO reconstructions \cite{LPR1,LPR2,C-2008,CPSV} are taken as basic reconstructions to prevent appearing of spurious oscillations. 
A combination of CWENO and MPP limiter has been adopted in  \cite{friedrich2019maximum} in the context of Eulerian schemes. We explore here the use of this technique in the context of semi-Lagrangian schemes. 

Three problems are considered in this paper, and each of them is solved by a high-order, conservative, semi-Lagrangian, finite-difference method. The detail of the method, however, depends on the problem.

As a warm up problem we consider rigid body rotation, for which a comparison with an exact solution is possible. 
High-order in time is obtained by operator splitting. 

The second problem we deal with is the Vlasov-Poisson system in plasma dynamics. Since the system has various time-invariant quantities, conservation and positivity are very important and become more relevant for long time simulation. There is a vast literature on SL methods for the VP system. Among them, splitting methods \cite{CMS-vlasov-2010,FSB-vlasov-2001,qiu2010conservative,QS} have been often adopted due to their simplicity in treating transport-type equations (see \cite{qiu2017high} for non-splitting methods). 


To construct a high-order schemes which maintain positivity and conservation, we combine our reconstruction to high-order splitting methods proposed in \cite{CCFM,yoshida1990construction}. Note in the multidimensional case we do not need to adopt dimensional splitting for transport-type equations
(see Section \ref{SplittingMethods}).

In the context of SL methods for BGK-type models, we recently proposed a predictor-corrector scheme \cite{BCRY} adopting a flux difference form as in \cite{RQT} together with a construction technique for the Maxwellian in \cite{M}. Although the scheme allows large time steps for small collision times, a related stability analysis \cite{BCRY} shows severe stability restrictions on the time step, especially in the rarefied regime. In this paper, we will return to the standard SL approach \cite{GRS} and apply our conservative reconstruction, thus obtaining a conservative SL scheme which performs well in all  range of Knudsen number with no CFL-type stability restrictions.

%

%
%
%

The rest of this paper is organised as follows: In section 2, we recall the  two kinetic models of interest, namely Vlasov-Poisson system of plasma physics and the BGK model of rarefied gas dynamics. In Section 3, we briefly recall  the conservative reconstruction introduced in \cite{BCRY1}. Section 4 is devoted to the construction of semi-Lagrangian schemes for specific applications, such as rigid rotation, Vlasov-Poisson system and BGK model. In Section 5, we describe high order splitting SL methods for rigid rotation and Vlasov-Poisson equations, and high order non-splitting SL schemes for BGK model. Finally, in Section 6, we provide several numerical tests demonstrating the performance of our conservative reconstruction applied to the SL schemes.


\section{Applications to kinetic models}
Here we consider two kinetic equations: the Vlasov-Poisson system and the BGK model.

\subsection{Vlasov-Poisson system}
\label{sec:VP}
In the plasma physics, the movements of charged particles forced by electric field is described by the Vlasov-Poisson system: 
\begin{align}\label{VP}
\begin{split}
&\frac{\partial{f}}{\partial{t}} + v \cdot \nabla_x{f} + \mathbb{E}(x,t) \nabla_v{f}=0\cr
&\mathbb{E}(x,t)= -\nabla_x\phi(x,t), \qquad  -\Delta_x \phi(x,t)= \rho(x,t) - m,\cr
&\rho(x,t)=\int_{\mathbb{R}^{d_v}} f(x,v,t)\,dv, \quad m=\frac{1}{L^{d_x}}\int_{\mathbb{T}^{d_x}}\int_{\mathbb{R}^{d_v}} f(x,v,t) \,dv\,dx.
\end{split}
\end{align}
Here $f(x,v,t)$ is the number density of charged particles on the phase point $(x,v)\in\mathbb{T}^{d_x}\times \mathbb{R}^{d_v}$ at time $t>0$ where the space $\mathbb{T}^{d_x}$ and $\mathbb{R}^{d_v}$ denote the $d_x$-dimensional torus $\mathbb{R}^{d_x}/(L\mathbb{Z}^{d_x})$ and ${d_v}$-dimensional real numbers, respectively. The notation $\mathbb{E}(x,t)$ is the electric field, and $\phi(x,t)$ is the self-consistent electrostatic potential. We use $\rho$ to denote charge density and $m$ is the ion density assumed to be uniformly distributed on the background. 

In the system, the distribution function $f$ satisfies the following properties:
\begin{enumerate}
	\item Maximum principle:
	\[
	0 \leq f(x,v,t) \leq M, \quad \forall t>0,
	\]
	provided $0 \leq f_0(x,v) \leq M$. 
	
	\item Conservation of total mass:
	\[
	\frac{d}{dt}\int_{\mathbb{R}^{d_v}}\int_{\mathbb{T}^{d_x}}f\,dx\,dv=0.
	\]
	\item Conservation of $L^p$ norm $1\leq p<\infty$:
	\begin{equation*}
	\|f\|_p=\left(\int_{\mathbb{R}^{d_v}}\int_{\mathbb{T}^{d_x}}|f|^p\,dx\,dv\right)^\frac1p.
	\end{equation*}
	\item Conservation of energy:
	\begin{equation}\label{Energy}
	\text{Energy}=\int_{\mathbb{R}^{d_v}}\int_{\mathbb{T}^{d_x}}\frac{|v|^2}{2}f\,dx\,dv + \int_{\mathbb{T}^{d_x}} \frac{|\mathbb{E}|^2}{2}\,dx.
	\end{equation}
	\item Conservation entropy:
	\begin{equation*}
	\text{Entropy}=\int_{\mathbb{R}^{d_v}}\int_{\mathbb{T}^{d_x}}f\log \frac{1}{f}\,dx\,dv.
	\end{equation*}
\end{enumerate}


\subsection{BGK model for the Boltmann equation}

The BGK model is a well-known approximation model for the Boltzmann equation of rarefied gas dynamics. The main feature of the BGK model is the replacement of the collision integral of the Boltzmann equation with a relaxation process. The BGK model reads:
\begin{equation}\label{bgk}
\frac{\partial{f}}{\partial{t}} + v \cdot \nabla_x{f} = \frac{1}{\kappa}\left(\mathcal{M}(f)-f\right),
\end{equation}
where $f(x,v,t)$ is the number density of monatomic particles on the phase space $\mathbb{R}^{d_x} \times \mathbb{R}^{d_v}$ at time $0 \leq t \in \mathbb{R}_+$. Here $\mathbb{R}$ is the real number, and $d_x$ and $d_v$ denote the dimension for space and velocity variables, respectively. We consider the fixed collision frequency with the Knudsen number $\kappa>0$ defined by a ratio between the mean free path of particle and the physical length scale under consideration. 
In \eqref{bgk}, 
the collision integral of the Boltzmann equation is replaced with the relaxation process towards the local thermodynamical equilibrium, so called the local Maxwellian $\mathcal{M}(f)$:
\[
\mathcal{M}(f)(x,v,t):=\frac{\rho(x,t)}{\sqrt{\left(2 \pi T(x,t) \right)^{d_v}}}e^{-\frac{|v-U(x,t)|^2}{2T}}.
\]
Here macroscopic densities for mass $\rho$, momentum $U$, energy $E$ and temperature $T$ are defined by
\begin{align*}
\begin{split}
\bigg(\rho,\, \rho U,\, 2E,\, d_v \rho T\bigg)^\top& = \int_{\mathbb{R}^{d_v}} \left(1,\,v,\,|v|^2,|v-U|^2\right) f  dv,
\end{split}
\end{align*}
Note that the relaxation operator still maintains qualitative features of the Boltzmann collision operator:
\begin{itemize}
	\item Collision invariance $ 1,v,|v|^2$:
	\begin{align}\label{cancellation property}
	\int_{\mathbb{R}^{d_v}} (\mathcal{M}(f)-f)\begin{pmatrix}
	1\\v\\ |v|^2
	\end{pmatrix}dv=0.
	\end{align}
	\item Conservation laws for mass, momentum and energy:
	\[
	\frac{d}{dt}  \int_{\mathbb{R}^{d_x} \times\mathbb{R}^{d_v} } f \phi(v) \,dx\,dv=0.
	\]
	\item The H-theorem:
	\[
	\frac{d}{dt}  \int_{\mathbb{R}^{d_x}\times\mathbb{R}^{d_v}} f \log{\frac{1}{f}} \,dx\,dv   \geq 0.
	\]
	\item In the fluid regime $\kappa \rightarrow 0$, the solution $f$ tends to $\mathcal{M}$. Then, the macroscopic moments of $f$ satisfy the compressible Euler system:
	\begin{align*}
	\begin{split}
	&\rho_t + \nabla_x \cdot(\rho u)=0,\cr
	&(\rho u)_t + \nabla_x \cdot(\rho u \otimes u+\rho T Id)=0,\cr
	& E_t + \nabla_x \cdot\left[\left(E+\rho T\right)u\right] =0, 
	\end{split}
	\end{align*}
\end{itemize}

\section{Conservatie Reconstruction}\label{sec: recon}
Here we briefly summarize the one-dimensional conservative reconstruction introduced in \cite{BCRY1}. Since extension to two-dimensional case is straightforward, we refer to \cite{BCRY1} for its detail. 

Now, we provide a description of the conservative reconstruction in the point-wise framework. Let us consider uniform mesh of size $\Delta x$ with grid points $x_i\equiv x_{min} + i \Delta x$, {of the computational domain $[x_{min}, x_{max}]$}. The set of index $i$ will be denoted by $\mathcal{I}$ . Suppose that $u(x)= \frac{1}{\Delta x} \int_{x-\Delta x/2}^{x-\Delta x/2} \hat{u}(y)dy$, the procedure for conservative reconstruction is given as follows:
\begin{enumerate}
	\item Given point-wise values $\{u_i\}_{i \in \mathcal{I}}$  for each $i \in \mathcal{I}$, reconstruct a polynomial of even degree $k$:$$R_i(x)= \sum_{\ell=0}^k \frac{R_i^{(\ell)}}{\ell !}(x-x_i)^{\ell}$$ such that
	\begin{itemize}
		\item High order accurate in the approximation of smooth $\hat{u}(x)$:
		\begin{itemize}
			\item If $\ell$ is an even integer such that $0 \leq \ell\leq k$,
			\begin{align*}
			\begin{split}
			\hat{u}_i^{(\ell)}&=R_i^{(\ell)}+ \mathcal{O}(\Delta x^{k+2-\ell}). 
			\end{split}
			\end{align*}
			\item If $\ell$ is an odd integer such that $0\leq \ell < k$,
			\begin{align*}
			\begin{split}
			\hat{u}_i^{(\ell)}-\hat{u}_{i+1}^{(\ell)}&=R_i^{(\ell)}-R_{i+1}^{(\ell)}+ \mathcal{O}(\Delta x^{k+2-\ell}).
			\end{split}
			\end{align*}
		\end{itemize}
		\item Essentially non-oscillatory. 
		\item Positive preserving. 
		\item Conservative in the sense of cell averages:
		\[
		\frac{1}{\Delta x}\int_{x_{i-\frac{1}{2}}}^{x_{i+\frac{1}{2}}}R_i(x)\,dx = u_{i}.
		\] 
	\end{itemize}    
	\item Using the obtained values $R_i^{(\ell)}$ for $0 \leq \ell\leq k$, approximate $u(x_{i+\theta})$ with
	\begin{align*}
	Q(x_i+\theta \Delta x)= \sum_{\ell=0}^k (\Delta x)^{\ell} \left(\alpha_{\ell}(\theta)R_{i}^{(\ell)} +\beta_\ell(\theta)R_{i+1}^{(\ell)}\right),
	\end{align*}
	where $\alpha_\ell(\theta)$ and $\beta_\ell(\theta)$ are given by 
	\begin{align*}
	\alpha_\ell(\theta) = \frac{1- (2\theta -1)^{\ell+1}}{2^{\ell+1}(\ell+1)!}, \quad 
	\beta_\ell(\theta) = \frac{(2\theta -1)^{\ell+1} - (-1)^{\ell+1} }  {2^{\ell+1}(\ell+1)!}.
	\end{align*}
	for $\theta \in [0,1)$.
\end{enumerate}

\begin{remark}
	In the previous paper, we adopt CWENO reconstructions  as basic reconstructions $R_i(x)$, which leads to high order conservative non-oscillatory reconstruction. Here, we can further impose the maximum principle preserving property on $Q(x)$ by using a linear scaling technique \cite{zhang2010maximum,zhang2011maximum,friedrich2019maximum}. 
	
	The procedure is following: we first compute linear scaling parameter:
	\begin{align}\label{limiter}
	\xi_i= \max\left\{1, \frac{|M-u_i|}{|M'_i-u_i|},\frac{|m-u_i|}{|m_i'-u_i|}\right\},
	\end{align} 
	where $M$ and $m$ are maximum and minimum values of a known function $u$ on the whole computational domain, and $M_i'$ and $m_i'$ are maximum and minimum values of local reconstruction $R_i(x)$ on the cell $I_i$. Next, we use $\xi_i$ to CWENO reconstructions $R_i(x)$ as follows:
	\begin{align*}
	\tilde{R}_i(x)= u_i + \xi_i( R_i(x) - u_i),
	\end{align*}
	where $\tilde{R}_i(x)$ is the rescaled polynomial defined on $I_i$:
	:$$\tilde{R}_i(x)= \sum_{\ell=0}^k \frac{\tilde{R}_i^{(\ell)}}{\ell !}(x-x_i)^{\ell}$$
	Consequently, 
	\begin{align*}
	\tilde{R}_i^{(0)}= u_i + \xi_i( R_i^{(0)} - u_i), \quad 	\tilde{R}_i^{(\ell)}= \xi_i R_i^{(\ell)}, \quad 1 \leq \ell \leq k,
	\end{align*}
	from which we obtain a positive reconstruction $Q(x_i + \theta \Delta x)$:
	\begin{align*}
	Q(x_i+\theta \Delta x)= \sum_{\ell=0}^k (\Delta x)^{\ell} \left(\alpha_{\ell}(\theta)\tilde{R}_{i}^{(\ell)} +\beta_\ell(\theta)\tilde{R}_{i+1}^{(\ell)}\right).
	\end{align*}
	
	We remark that maximum preserving property still holds for $Q(x)$ for any $x$, because $Q(x)$ is obtained by taking sliding average of piecewise reconstructions satisfying the maximum preserving property:
	\[
	m \leq Q(x) = \frac{1}{\Delta x} \int_{x- \Delta x/2}^{x+ \Delta x/2} R(y)\,dy \leq M.
	\]	
	
	In some applications like the Vlasov-Ppisson model, it is not necessary to impose non-oscillatory property.
	In this case, $R_i(x)$ will be polynomials interpolating $u_i$ in terms of cell average, i.e. the CWENO polynomial with linear weights.
\end{remark}


\section{Semi-Lagrangian methods}\label{SLmethods}
A semi-Lagrangian method can be obtained by integrating the equation of interest along its characteristic curve. As the simplest case, we consider a one-dimensional transport equation:
\begin{align}\label{linear}
\partial_t f + a \partial_x f &= 0, \quad t\geq 0,\, x\in \mathbb{R},
\end{align}
where $a$ is a fixed constant in $\mathbb{R}$. 
Let us denote the $n$-th time step by $t^n:=n\Delta t$, $n \in \mathbb{N}$ for $\Delta t>0$. Given the distribution function $f$ at time $t^n$, the exact solution on $x_i$ at $t^{n+1}$ is given by $$f(x_i,t^{n+1})= f(X(t^n;x_i,t^{n+1}),t^n),$$ 
where $X(t;x_i,t^{n+1})$ satisfies
$$\frac{dX}{dt}=a, \quad X(t^{n+1})=x_i.$$
That is, we have $X(t^n;x_i,t^{n+1})= x_i-a(t^{n+1}-t^n)$, and this gives
$$f(x_i,t^{n+1})= f(x_i-a\Delta t,t^n).$$
Considering that the characteristic foot $X(t^n;x_i,t^{n+1})=x_i-a\Delta t$ can be located on off-grid points, one can update solutions after reconstructing the solution $f(x-a\Delta t,t^n)$ from given solutions on grid points. For the semi-Lagrangian schemes proposed in this paper, we will use the conservative reconstruction described in section \ref{sec: recon}. In the following we will discuss the application of semi-Lagrangian methods to a basic test problem: {\em the rigid body rotation} and to two kinetic equations: {\em Vlasov-Poisson system} and {\em BGK model}.


\subsection{Semi-Lagrangian methods for Rigid body rotation}
Rigid body rotation of objects can be described by a linear model:
\begin{align} \label{rigid}
u_t -yu_x +xu_y = 0,\quad  (x,y) \in \mathbb{R}^2.
\end{align}
One can design a semi-Lagrangian method using the information that the solution is constant along the characteristic curves: $$u(t) \equiv u(X(t;x,y,t^{n+1}),Y(t;x,y,t^{n+1})$$
where $\left(X(t), Y(t)\right)\equiv \left(X(t;x,v,t^{n+1}), Y(t;x,v,t^{n+1})\right)$ satisfies:
\begin{align*}
\frac{dX}{dt}= -Y, \quad \frac{dY}{dt}= X.
\end{align*}
In this case the backward characteristics can be solved exactly as 
\begin{equation}
X(t)= x\cos(t^{n+1}-t) + y\sin(t^{n+1}-t), \quad Y(t)= -x\sin(t^{n+1}-t) + y\cos(t^{n+1}-t),
\label{eq:exact_char}
\end{equation}
however, we shall not make use of such exact solution. Furthermore, the resulting scheme would not be exactly  conservative.

A general way to construct conservative schemes is based on splitting methods.  A brief introduction of splitting methods is postponed to the Section \ref{SplittingMethods}. 

\begin{remark}
	Since the analytical solution to \eqref{rigid} satisfies the maximum principle preserving property, it should be taken into account for spatial reconstructions. Furthermore, when initial solutions are discontinuous, it is also important to capture the position of discontinuity without spurious oscillations. For these reasons, we adopt CWENO reconstructions together with an MPP limiter \eqref{limiter} for conservative reconstructions (see section \ref{sec rigid}.)
\end{remark}

\subsection{Semi-Lagrangian methods for Vlasov-Poisson system}
\label{sec:VPS}
Following the notation of Sect. (\ref{SLmethods}), we set $a \equiv (v, \mathbb{E}(t,x,v))$ and $\nabla \equiv (\nabla_x, \nabla_v)$. Then we can rewrite (\ref{VP}) in the form (\ref{linear}). If $a$ is sufficiently smooth (Lipschitz continuous), we can define unique characteristic curves of the first order differential operator,
$$
\frac{\partial}{\partial t} + a \cdot \nabla, 
$$
which solve the following system
\begin{align}\label{CharC}
\begin{split}
\frac{d X}{dt}(t;x,v,t^{n+1}) &= V(t;x,v,t^{n+1}),\cr
\frac{d V}{dt}(t;x,v,t^{n+1}) &= \mathbb{E}(t,X(t;x,v,t^{n+1})).
\end{split} 
\end{align}	
Here $(X(t;x,v,t^{n+1}),V(t;x,v,t^{n+1}))$ denotes the position in phase space at the time $t$, which passes $(x, v)$ at time $t^{n+1}$. Since the distribution function 
of the VP system is constant along the particle
trajectories, we have $$f(x,v,t^{n+1})=f(X(t^n;x,v,t^{n+1}),V(t^n;x,v,t^{n+1})).$$
To update solution in this way, however, it is necessary to trace all characteristic curves by solving solving (\ref{CharC}). 
We note that this procedure ca be simplified a lot using splitting methods where we are treating two simple transport-type equations with an Poisson equation. High-order time splitting methods used for VP system in the numerical tests are listed in Sect. \ref{SplittingMethods}.

\begin{remark} Hereafter, let us denote velocity grids by $v_j=v_{min}+(j-1)\Delta v$ with mesh size $\Delta v >0$ and the set of indices $j$ by $\mathcal{J}$. This notation will be also used for the description  of numerical methods for BGK model. In \eqref{VP}, considering that $\mathbb{E} = -\nabla_x \phi$, one should find the potential function $\phi_i^n$ in advance. For this, we solve the Poisson equation $-\Delta_x \phi= \rho - m$ in \eqref{VP} with a Fast Fourier Transform using $\{\rho_i^n\}_{i \in \mathcal{I}}$ and $m^n$ obtained by the rectangular rule: 
	\begin{align}\label{m}
	\rho_i^n=\sum_j f_{i,j}^n (\Delta v)^{d_v},\quad  m^n=\frac{1}{L^{d_x}}\sum_{i,j} f_{i,j}^n (\Delta x)^{d_x}(\Delta v)^{d_v}.
	\end{align} 
	We remark that the initial value of $m^n$ is maintained during time evolution of numerical solutions obtained by mass conservative schemes, i.e. $m^n=m^0$.
	
	In one space dimension, we calculate the discrete electric field $E_i^n=-\partial_x \phi\big|_{(x,t)}$ as follows:
	\begin{itemize}
		\item 4th order approximation
		\begin{align}\label{phi 4}
		\partial_x \phi\big|_{(x,t)=(x_i,t^n)} = \frac{8(\phi_{i+1}^n - \phi_{i-1}^n)  -  (\phi_{i+2}^n - \phi_{i-2}^n)}{12\Delta x}+ \mathcal{O}(\Delta x)^4.
		\end{align}	
		\item 6th order approximation
		\begin{align}\label{phi 6}
		\partial_x \phi\big|_{(x,t)=(x_i,t^n)} = \frac{45(\phi_{i+1}^n - \phi_{i-1}^n)  -  9(\phi_{i+2}^n - \phi_{i-2}^n)  + (\phi_{i+3}^n - \phi_{i-3}^n)}{60\Delta x}+ \mathcal{O}(\Delta x)^6.
		\end{align}
		\item 8th order approximation
		\begin{align}\label{phi 8}
		\partial_x \phi\big|_{(x,t)=(x_i,t^n)} = \frac{672(\phi_{i+1}^n - \phi_{i-1}^n)  -  168(\phi_{i+2}^n - \phi_{i-2}^n)  + 32(\phi_{i+3}^n - \phi_{i-3}^n) - 3(\phi_{i+4}^n - \phi_{i-4}^n)}{840\Delta x}+ \mathcal{O}(\Delta x)^6.
		\end{align}
	\end{itemize}
	The approximation of $\nabla_x \phi$ in multi-dimension can be done by dimension-wise computations.
\end{remark}

\begin{remark}
	Considering the maximum principle of the Vlasov-Poisson system, it is natural to adopt a basic reconstruction which maintains the maximum and minimum of initial data together with conservation property. For these reasons, we use the MPP limiter \eqref{limiter} on polynomial interpolation in the sense of cell average.
	
\end{remark}


\subsection{A semi-Lagrangian method for BGK model}
Here, we review a class of non-splitting implicit semi-Lagrangian methods for BGK model \cite{GRS,RS}. Note that we will apply the conservative reconstruction in section \ref{sec: recon} to the reconstruction of distribution on off-grid points.

\subsubsection{First order semi-Lagrangian method for BGK model}
A first order SL method for the BGK model \eqref{bgk} is obtained by applying implicit Euler method to its Lagrangian framework:
\begin{align}\label{B-3}
\begin{split}
\frac{df}{dt}(X(t;x_i,v_j,t^{n+1}),v_j,t)&=\frac{1}{\kappa}\left(\mathcal{M}(f)-f\right)(X(t;x_i,v_j,t^{n+1}),v_j,t)\cr
X(t;x_i,v_j,t^{n+1})&=x_i-v_j(t^{n+1}-t),
\end{split}
\end{align}
which gives
\begin{align}\label{2D IE}
f_{i,j}^{n+1}=\tilde{f}_{i,j}^n + \frac{\Delta t}{\kappa} \left(\mathcal{M}[f]_{i,j}^{n+1}-f_{i,j}^{n+1}\right),
\end{align}
where $\tilde{f}_{i,j}^n$ is an approximation of $f(x_i - v_j \Delta t,v_j,t^n)$ which can be reconstructed from  $\{f_{i,j}^n\}_{i \in \mathcal{I}}$. The quantity $\mathcal{M}[f]_{i,j}^{n+1}$ is the local Maxwellian computed by
\begin{align*}
\mathcal{M}[f]_{i,j}^{n+1} =\frac{\rho_{i}^{n+1}}{\sqrt{\left(2 \pi T_{i}^{n+1} \right)^2}}\exp\left(-\frac{|v_{j}-U_{i}^{n+1}|^2}{2T_{i}^{n+1}}\right),
\end{align*}
with discrete macroscopic quantities:
\begin{align*}
\left(\rho_{i}^{n+1},\, \rho_{i}^{n+1}U_{i}^{n+1},\, d_v \rho_{i}^{n+1} T_{i}^{n+1}\right) &:= \sum_{j \in \mathcal{J}} f_{i,j}^{n+1} \left(1,v_j,\big| v_j- U_{i}^{n+1} \big|^2\right)(\Delta v)^{d_v}.
\end{align*}
Using the structural feature of the relaxation term of \eqref{bgk}, one can compute \eqref{2D IE} explicitly. Multiplying the collision invariant $\phi_{j}:=(1,v_{j},\frac{|v_{j}|^2}{2})$ to both sides of \eqref{2D IE} and taking summation over $j \in \mathcal{J}$, one can obtain
\begin{align}\label{approximation}
\sum_{j \in \mathcal{J}} \big(f_{i,j}^{n+1} - \tilde{f}_{i,j}^n\big)\phi_{j} (\Delta v)^{d_v} = \frac{\Delta t}{\kappa}\sum_{j \in \mathcal{J}} \left(\mathcal{M}[f]_{i,j}^{n+1}-f_{i,j}^{n+1}\right)\phi_{j} (\Delta v)^{d_v}.
\end{align}
Here the right hand side vanishes if the number of velocity girds is enough and the velocity domain is large enough. Then, the discrete macroscopic quantities $\rho_i^n$, $U_i^n$ and $T_i^n$ can be computed as follows: 
\begin{align}
\left(\rho_{i}^{n+1},\, \rho_{i}^{n+1}U_{i}^{n+1},\, E_{i}^{n+1}\right) &=\sum_{j \in \mathcal{J}} f_{i,j}^{n+1}\phi_{j} (\Delta v)^{d_v} = \sum_{j \in \mathcal{J}} \tilde{f}_{i,j}^n \phi_{j} (\Delta v)^{d_v}=:\left(\tilde{\rho}_{i}^{n+1},\, \tilde{\rho}_{i}^{n+1}\tilde{U}_{i}^{n+1},\, \tilde{E}_{i}^{n+1}\right),
\end{align}
this further gives
\begin{align}
d_v \rho_{i}^{n+1}T_{i}^{n+1} &=\sum_{j \in \mathcal{J}} f_{i,j}^{n+1} \big|v_j-U_i^{n+1}\big|^2(\Delta v)^{d_v} = \sum_{j \in \mathcal{J}} \tilde{f}_{i,j}^n \big|v_j-\tilde{U}_i^{n+1}\big|^2 (\Delta v)^{d_v} = d_v\tilde{\rho}_{i}^{n+1}\tilde{T}_{i}^{n+1},
\end{align}
Consequently, one can update solution as follows:
\begin{align}\label{first order scheme}
f_{i,j}^{n+1}=\tilde{f}_{i,j}^n + \frac{\Delta t}{\kappa} \left(\mathcal{M}[\tilde{f}]_{i,j}^{n}-f_{i,j}^{n+1}\right),
\end{align}
with
\begin{align*}
\mathcal{M}[\tilde{f}]_{i,j}^{n} =\frac{\tilde{\rho}_{i}^{n}}{\sqrt{\left(2 \pi \tilde{T}_{i}^{n} \right)^{d_v}}}\exp\left(-\frac{|v_{j}-\tilde{U}_{i}^{n}|^2}{2\tilde{T}_{i}^{n}}\right).
\end{align*}


\subsubsection{Conservative approximation of Maxwellian}
\label{sec:CM}
The DIRK and BDF based semi-Lagrangian schemes are based on discrete velocity model (DVM), so conservation errors coming from the computation of the Maxwellian is crucial as much as the errors coming from the reconstruction of the solutions at characteristic foots. (See \cite{BCRY}.) 

Consider $d_v$-dimensional velocity domain with grid points $v_j$ for $j \in \mathcal{J}$. Given macroscopic moments $\mathcal{U}:=(\rho,\rho U, E)^{\top} \in \mathbb{R}^{d_v+2}$, one can compute a Maxwellian $\left\{\mathcal{M}_j\right\}_{j \in \mathcal{J}}$ as follows:
\[
\mathcal{M}_j:=\frac{\rho}{\sqrt{\left(2 \pi T \right)^{d_v}}}\exp\left(-\frac{\left|v_j-U\right|^2}{2T}\right), \quad j \in \mathcal{J},
\]
with the relation $d_v \rho T = 2E-\rho U^2$. Since we only use finite number of velocity grids in a truncated domain, the discrete moments of $\left\{\mathcal{M}_j\right\}_{j \in \mathcal{J}}$ cannot exactly reproduce the given moments $\mathcal{U}$, i.e.
$ \sum_{j \in \mathcal{J}} \mathcal{M}_j \phi_{j} (\Delta v)^{d_v} \ne \mathcal{U} $, for $\phi_j=(1,v_j,\frac{|v_j|^2}{2})$. This issue can be treated by finding an Maxwellian $\left\{g_j\right\}_{j \in \mathcal{J}}$ such that
\begin{align}\label{same momemnt}
\sum_{j \in \mathcal{J}} g_j \phi_{j} (\Delta v)^{d_v} = \mathcal{U}.
\end{align}
Since the set $\mathcal{J}$ contains $(N_v+1)^{d_v}$ meshes, and $\left\{\mathcal{M}_j\right\}_{j \in \mathcal{J}}$ and $\left\{g_j\right\}_{j \in \mathcal{J}}$ are vectors in $\mathbb{R}^{(N_v+1)^{d_v}}$. We will write  $$\mathcal{M} \equiv (\mathcal{M}_1,\mathcal{M}_2,\cdots,\mathcal{M}_{(N_v+1)^{d_v}})^\top, \quad g =(g_1,g_2,\cdots,g_{(N_v+1)^{d_v}})^\top .$$

We now briefly review two ways introduced in \cite{M,gamba2009spectral} which reduce the conservation errors coming from the computation of the Maxwellian.

\paragraph{$L^2$ minimization ($dM_1$)}
The first approach is to solve the following minimization problem:
\begin{align}\label{min prb}
\min \big\|\mathcal{M} - g\big\|_2^2 \quad\text{ s.t }\quad Cg=\mathcal{U}
\end{align}
where $$\displaystyle C_{(d_v+2) \times (N_v+1)^{d_v}}:=\begin{pmatrix}
(\Delta v)^{d_v}\\ v_j (\Delta v)^{d_v} \\ \frac{|v_j|^2}{2} (\Delta v)^{d_v}
\end{pmatrix}, \quad \mathcal{U}_{(d_v+2) \times 1} = (\rho,\rho U, E)^{\top}.$$
To solve the minimization problem \eqref{min prb}, we use the Lagrange multiplier method with the following Lagrangian $\mathcal{L}(M,\lambda)$:
$$\mathcal{L}(g,\lambda)= \big\|\mathcal{M} - g\big\|_2^2 + \lambda^{\top}(Cg-\mathcal{U}). $$
We first compute the gradient of $\mathcal{L}$ with respect to $M$ and $\lambda$:
\begin{align*}
\begin{array}{rlrl}
\nabla_g \mathcal{L} =0 &<=>& g&= \mathcal{M}+ \frac{1}{2}C^{\top} \lambda\cr
\nabla_\lambda \mathcal{L} = 0 &<=>& Cg&=\mathcal{U}.
\end{array}
\end{align*}
Then, we use these to derive
\begin{align*}
Cg= C\mathcal{M}+ \frac{1}{2}CC^{\top} \lambda &<=> \mathcal{U}= C\mathcal{M}+ \frac{1}{2}CC^{\top} \lambda \cr
&<=>  \lambda= 2(CC^{\top})^{-1}(\mathcal{U}-C\mathcal{M}).
\end{align*}
Note that the matrix $CC^{\top}$ is invertible because it is symmetric and positive definite. Consequently,
\begin{align}\label{first dm}
g= \mathcal{M}+C^{\top}(CC^{\top})^{-1}(\mathcal{U}-C\mathcal{M}),
\end{align}
Once the matrix $C^{\top}(CC^{\top})^{-1}$ is computed, it can be stored and reused for other Maxwellian computations. In view of this, the minimization approach \eqref{min prb} only requires additional matrix computations in determining $g$.

\paragraph{Discrete Maxwellian (Entropy minimization, $dM2$)}
The second approach is to determine the value of $g$ is to consider the following entropy minimization problem:
\begin{align}\label{min entro}
\min \left\{ \sum_{j\in \mathcal{J}}g_j\log{g_j},\,  g_j \geq 0 \quad \text{s.t.} \quad  \sum_{j \in \mathcal{J}} g_j \phi_{j} (\Delta v)^{d_v} = \mathcal{U} \right\}
\end{align}
which is equivalent to solving the following equation:
\begin{align}\label{root prb}
\sum_{j\in\mathcal{J}} \left(\exp(a \cdot \phi_j) \right)\phi_j (\Delta v)^{d_v}=\mathcal{U} \quad\text{where}\quad  a \in \mathbb{R}^{d_v+2}.
\end{align}
Then, the root $a \in \mathbb{R}^{d_v+2}$ of \eqref{root prb} gives the solution to \eqref{min entro}:
\begin{align}\label{second dm}
g_j=\exp(a \cdot \phi_j) \text{ for } j\in \mathcal{J}.
\end{align}

Compared to the first method \eqref{first dm}, this method can be more complicate in that one should implement an root find algorithm such as Newton's method to solve \eqref{root prb}. However, this approach always gives non-negative Maxwellian, hence it can be more helpful to preserve the positivity of numerical solution to the BGK model \eqref{bgk}.

\section{High order time discretization}
In this section we present high order integrators for the time discretization that we use in our numerical tests. In the context of semi-Lagrangian methods for rigid body rotation \eqref{rigid} and VP system (\ref{VP}), splitting methods can be a more suitable approach in terms of conservation. On the other hand, for BGK model we consider diagonally implicit R-K method and BDF multistep methods \cite{HW} for the treatment of stiffness arising in the small values of Knudsen number.
\subsection{High-order splitting methods}\label{SplittingMethods}
We consider an arbitrary system $\dot{y} = f_1(y) + f_2(y)$ in $\mathbb{R}^n$. The idea of splitting methods is to decompose this system into two parts over a time: $\dot{y} = f_1(y)$ and $\dot{y} = f_2(y)$, and compute solutions of the two systems successively to yield a solution to the original system. For more details on splitting methods, we refer to \cite{mclachlan2002splitting}. Now, we apply the idea of splitting methods in the framework of semi-Lagrangian methods for the rigid body rotation \eqref{rigid} and VP system (\ref{VP}).

\subsubsection*{High-order splitting methods for rigid body rotation.}
In splitting methods for \eqref{rigid}, we are to solve two problems: \begin{align} \label{rigid split}
u_t -yu_x= 0,\quad u_t + xu_y = 0.
\end{align}
To employ semi-Lagrangian methods, we need to deal with two problems:
\begin{itemize}
	\item Shift solutions along $x$-axis over a time step $\tau>0$: In the first equation of \eqref{rigid split}, we solve 
	\begin{align} \label{rigid split1}
	\frac{d}{dt}u\left(X(t;x,y,t+\tau),Y(t;x,y,t+\tau),t\right)= 0, \quad \frac{dX}{dt}=-Y, \quad \frac{dY}{dt}=0.
	\end{align}
	The solution to this system is given by
	\[  u(t+\tau)\equiv u\left(x,y,t+\tau\right) = u\left(X(t),Y(t),t\right) = u\left(x+y\tau,y,t\right) \equiv u(t).
	\]
	For later use, we denote this flow by
	\begin{align*}
	u(t+\tau)&=e^{\tau\mathcal{X}}u(t).
	\end{align*}
	\item Shift along $y$-axis over a time step $\tau>0$. For the second equation of \eqref{rigid split}, we consider
	\begin{align} \label{rigid split2}
	\frac{d}{dt}u\left(X(t;x,y,t+\tau),Y(t;x,y,t+\tau),t\right) = 0, \quad \frac{dX}{dt}=0, \quad \frac{dY}{dt}=X.
	\end{align}
	In this case, we have
	\[  u(t+\tau)\equiv u\left(x,y,t+\tau\right) = u\left(X(t),Y(t),t\right) = u\left(x,y-x\tau,t\right) \equiv u(t),
	\]
	and the corresponding flow will be represented as follows:
	\begin{align*}
	u(t+\tau)&=e^{\tau\mathcal{Y}}u(t).
	\end{align*}
\end{itemize}
Note that each problem can be dealt with as in the linear transport equation case \eqref{linear}. To achieve high order accuracy in the time, one should successively solve (\ref{rigid split1}) and (\ref{rigid split2}) with suitable time steps. For this, we here adopt splitting methods used in \cite{yoshida1990construction}. 
In the following descriptions for high order splitting methods, we denote the solution $u$ at time $t^n$ by $u^n$.
\begin{itemize}
	\item 2nd order Strang splitting method: 
	\begin{enumerate}
		\item Shift along $y$-axis: $u^*=
		e^{\frac{\Delta t}{2}\mathcal{Y}}u^n$.
		\item Shift along $x$-axis: $u^{**}=
		e^{\Delta t\mathcal{X}}u^*$.
		\item Shift along $y$-axis: $u^{n+1}=e^{\frac{\Delta t}{2}\mathcal{Y}}u^{**}$.
	\end{enumerate}
	We can write this compactly as 
	\begin{align}\label{strang}
	u^{n+1}=\mathcal{S}_2^{\Delta t}u^n:=
	e^{\frac{\Delta t}{2}\mathcal{Y}}e^{\Delta t\mathcal{X}}e^{\frac{\Delta t}{2}\mathcal{Y}}u^n.
	\end{align}
	\item 4th order splitting method \cite{yoshida1990construction}:  
	\begin{align}\label{yoshida 4}
	u^{n+1}=\mathcal{S}_4^{\Delta t}u^n:=\mathcal{S}_2^{x_1\Delta t}\mathcal{S}_2^{x_0\Delta t}\mathcal{S}_2^{x_1\Delta t}u^n,
	\end{align}
	where $x_0=-\frac{2^{1/3}}{2-2^{1/3}}$, $x_1=\frac{1}{2-2^{1/3}}$.
	\item 6th order splitting method \cite{yoshida1990construction}:
	\begin{align}\label{yoshida 6}
	u^{n+1}=\mathcal{S}_6^{\Delta t}u^n:=\mathcal{S}_4^{y_1\Delta t}\mathcal{S}_4^{y_0\Delta t}\mathcal{S}_4^{y_1\Delta t}u^n,
	\end{align}
	where $y_0=-\frac{2^{1/5}}{2-2^{1/5}}$, $y_1=\frac{1}{2-2^{1/5}}$.
\end{itemize}


\subsubsection*{High-order splitting methods for VP system}
For splitting methods to VP system \eqref{VP}, we separate it into the following two parts:
\begin{align}\label{VP2}
\begin{split}
\frac{\partial{f}}{\partial{t}} + v \cdot \nabla_x{f}=0,
\end{split}
\end{align}
\begin{align}\label{VP3}
\begin{split}
\frac{\partial{f}}{\partial{t}} + \mathbb{E}(x,t) \cdot \nabla_v{f}=0.
\end{split}
\end{align}

For the semi-Lagrangian treatment of two systems, we rewrite them as follows: 
\begin{itemize}
	\item Advection step over a time step $\tau>0$: The equation \eqref{VP2} can be solved as in \eqref{linear}. Consider
	\begin{align} \label{VP split1}
	\frac{d}{dt}f(X(t;x,v,t+\tau),V(t;x,v,t+\tau),t)= 0, \quad \frac{dX}{dt}=V, \quad \frac{dY}{dt}=0,
	\end{align}
	which gives 
	\[  f(t+\tau)\equiv f\left(x,v,t+\tau\right) = f\left(X(t),V(t),t\right) = f\left(x-v\tau,v,t\right) \equiv f(t).
	\]
	We define the associated flow by
	\begin{align*}
	f(t+\tau)&=e^{\tau\mathcal{T}}f(t).
	\end{align*}
	\item Drift step over a time step $\tau>0$.
	\begin{align} \label{VP split2}
	\frac{d}{dt}f(X(t;x,v,t+\tau),V(t;x,v,t+\tau),t)= 0, \quad \frac{dX}{dt}=0, \quad \frac{dV}{dt}=\mathbb{E}(X(t;x,v,t+\tau),t)
	\end{align}
	In this case, we have
	\[  f(t+\tau)\equiv f\left(x,v,t+\tau\right) = f\left(X(t),V(t),t\right) = f\left(x,v-\mathbb{E}(x,t)\tau,t\right) \equiv f(t).
	\]
	Here we can compute the electric field $\mathbb{E}(x,t)$ using $\{f(t)\}$. Depending on the desired accuracy of spatial reconstructions, we will use one of \eqref{phi 4}, \eqref{phi 6} or \eqref{phi 8}. This flows will be denoted with the following notation: 
	\begin{align*}
	f(t+\tau)&=e^{\tau\mathcal{U}}f(t).
	\end{align*}
\end{itemize}
Note that both characteristic curves are given by straight lines, which enables one to compute the location of characteristic foots analytically.\\ 
The extension to high order methods for VP system is straightforward. For second and fourth order splitting methods, we will use (\ref{strang}) and (\ref{yoshida 4}), respectively. To attain sixth order in time, we use a 11-stage splitting method constructed in \cite{CCFM}.
Let us denote the solution $f$ at time $t^n$ by $f^n$. Then, the 6th order 11-stage splitting method in \cite{CCFM} can be described as follows:	 
\begin{align}\label{meren}
f^{n+1}=\psi_{6}^{\Delta t}f^n:=e^{\Delta t\mathcal{D}_1}e^{a_1\Delta t\mathcal{T}}e^{\Delta t\mathcal{D}_2}e^{a_2\Delta t\mathcal{T}}e^{\Delta t\mathcal{D}_3}e^{a_3\Delta t\mathcal{T}}e^{\Delta t\mathcal{D}_3}e^{a_2\Delta t\mathcal{T}}e^{\Delta t\mathcal{D}_2}e^{a_1\Delta t\mathcal{T}}e^{\Delta t\mathcal{D}_1}f^n,
\end{align}
where $\mathcal{D}_i=(b_i + 2c_im (\Delta t)^2 + 4d_im^2 (\Delta t)^4 - 8e_im^3 (\Delta t)^6)\mathcal{U}$, $i=1,2,3$. The value of $m$ is computable as in \eqref{m}, and the coefficients in $D_i$ are given as follows:
\begin{align*}
a_1&=  0.168735950563437422448196\cr
a_2&=  0.377851589220928303880766\cr
a_3&=  -0.093175079568731452657924\cr
b_1&=  0.049086460976116245491441\cr
b_2&=  0.264177609888976700200146\cr
b_3&=  0.186735929134907054308413\cr
c_1&=  -0.000069728715055305084099\cr
c_2&=  -0.000625704827430047189169\cr
c_3&=  -0.002213085124045325561636\cr
d_2&=  -2.916600457689847816445691 \cdot 10^{-6}\cr
d_3&=  3.048480261700038788680723 \cdot 10^{-5}\cr
e_3&=  4.985549387875068121593988 \cdot 10^{-7}\cr
d_1&= e_1 =e_2 =0.
\end{align*}	

\subsection{High-order Runge-Kutta methods for BGK model}
In \cite{GRS,RS}, considering the stiffness arising in the fluid limit $\kappa \to 0$ of BGK model \eqref{bgk}, high order semi-Lagrangian methods are obtained by combining high order L-stable diagonally implicit Runge-Kutta methods (DIRK) or backward-difference formulas (BDF) with high order spatial reconstructions. 
In this paper, we use a second order L-stable diagonally implicit Runge-Kutta (DIRK) scheme \cite{HWN} and a third order DIRK method introduced in \cite{BCRY}. We represent the methods using the representation of Butcher's table:
\begin{align*}
\begin{array}{c|c}
c & A \\
\hline
& b^{\top} 
\end{array}
\end{align*}
where $A=[a_{ij}]$ is a $s \times s$ matrix such that $a_{ij}=0$ for $i<j$, $c= (c_1,...,c_s)^{\top}$ and
$b=(b_1,...,b_s)^{\top}$ are coefficients vectors. The corresponding Butcher's tables are given as follows:
\begin{align}\label{Butcher}
\text{DIRK2}=\begin{array}{c|c c}
\alpha & \alpha & 0 \\
1 & 1-\alpha & \alpha \\
\hline
& 1-\alpha & \alpha
\end{array}, \quad 
\text{DIRK3}=\begin{array}{c|c c c}
\gamma & \gamma & 0 & 0\\
c_2 & c_2-\gamma_2 & \gamma_2 & 0\\
1 & 1-b_2-\gamma & b_2 & \gamma\\
\hline
& 1-b_2-\gamma & b_2 & \gamma_3
\end{array}
\end{align}
where $\alpha= 1 -\frac{\sqrt{2}}{2}$, $\gamma=0.3$, $\gamma_2 =13/3$, $b_2=-3/710$ and $c_2 = 8/3$. 


In the description of $s$-stage DIRK method, we use the following notation:
\begin{itemize}
	\item The location of characteristic foot on the $\ell$-th stage along the $k$-th backward-characteristic starting from $x_i$ with $v_j$ at time $t^n + c_k \Delta t$:
	\[
	x_{i,j}^{k, \ell}:= x_{i}-(c_k-c_\ell) v_{j}\Delta t
	\]
	\item The stage value of $f$ on $x_{i,j}^{k, \ell}$:
	\[F_{i,j}^{(k,\ell)} := f(x_{i}-v_{j}(c_k-c_\ell)\Delta t, v_{j}, t^n + c_\ell\Delta t).\]
	Setting $c_0=0$, we define
	\[F_{i,j}^{(k,0)} := f(x_{i}-v_{j}c_k\Delta t,v_{j}, t^n).\]
	\item RK fluxes:
	\[
	K(x,v,t) = \frac{1}{\kappa}\left(\mathcal{M}(x,v,t)-f(x,v,t)\right).
	\]
	and its discrete values on $x_{i,j}^{k, \ell}$
	\[
	K_{i,j}^{(k,\ell)} = \frac{1}{\kappa}\big( \mathcal{M}(F_{i,j}^{(k,\ell)})-F_{i,j}^{(k,\ell)} \big).
	\]
\end{itemize}

Now, we illustrate DIRK methods:
\subsubsection*{Algorithm for s-order DIRK based methods} For $k=1,\dots,s$, 
\begin{figure}[t]
	\centering
	\includegraphics[width=0.6\linewidth]{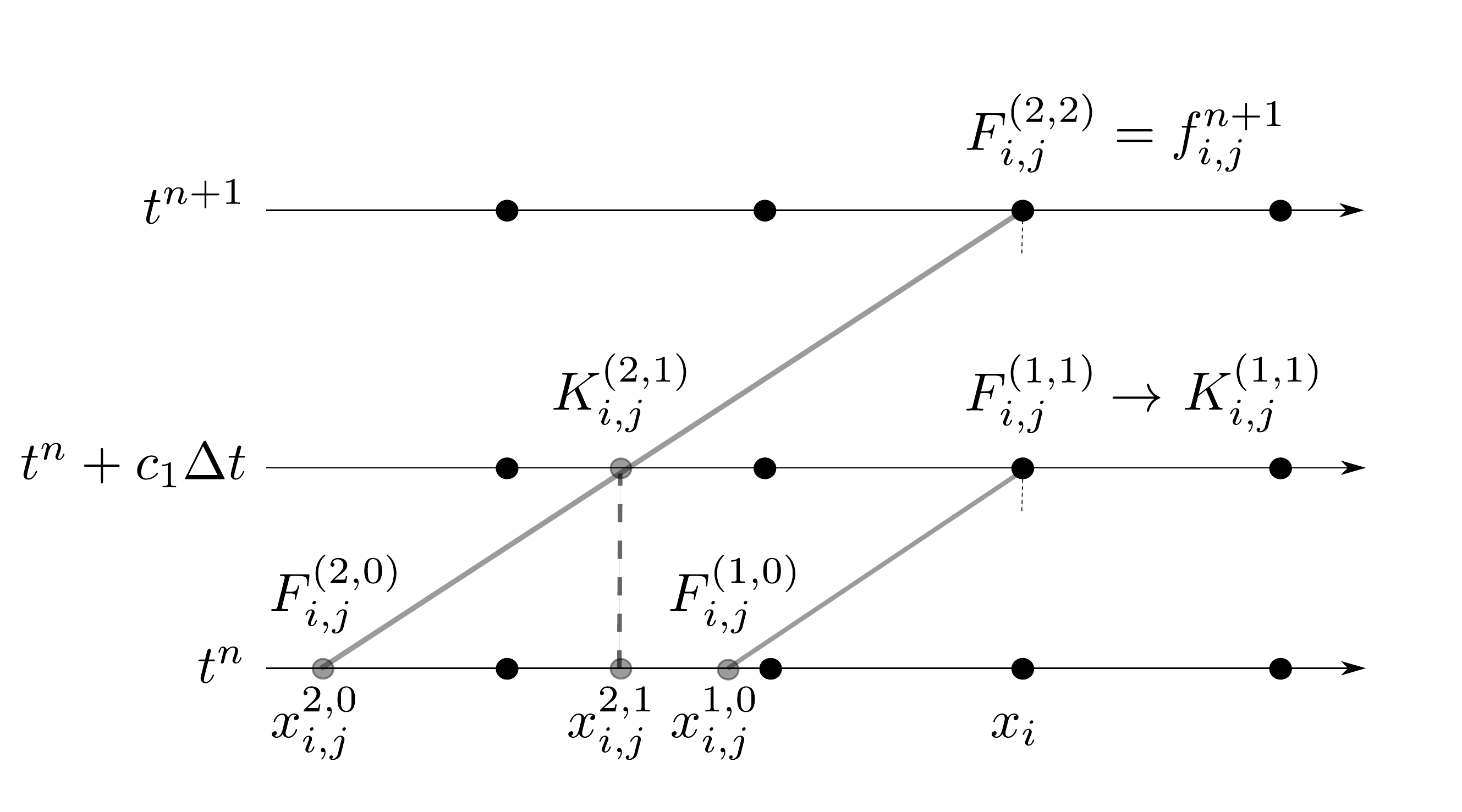}			\caption{Representation of DIRK2 scheme. Black circles: grid nodes, grey circles: points where interpolation is needed.}\label{FigRK2}	
\end{figure}
\begin{enumerate}
	\item Interpolate $F_{i,j}^{(k,0)}$ on $x_{i,j}^{k,0}$ from $\{F_{i,j}^{(k,0)}\}_{i\in \mathcal{I}}$.
	\item For $\ell=1\dots,k-1$, interpolate $K_{i,j}^{(k,\ell)}$ on $x_{i,j}^{k, \ell}$ from $\{K_{i,j}^{(\ell,\ell)}\}_{i\in \mathcal{I}}$. (Skip for $k=1$)
	\item
	Compute: 
	\[
	F_{i,j}^{(k,k)} = F_{i,j}^{(k,0)} + \Delta t \sum_{\ell=1}^{k-1} a_{k \ell}K_{i,j}^{(k,\ell)} + \frac{\Delta t}{\kappa}a_{kk}\left( \mathcal{M}\left[F_{i,j}^{(k,k)}\right] - F_{i,j}^{(k,k)} \right),
	\]
	where we use $\mathcal{M}\left[F_{i,j}^{(k,k)}\right]=\mathcal{M}\left[F_{i,j}^{(k,0)} + \Delta t \sum_{\ell=1}^{k-1} a_{k \ell}K_{i,j}^{(k,\ell)}\right].$
\end{enumerate}
Due to the property that the adopted DIRK methods are stiffly accurate (SA),
one can update solutions by putting $f^{n+1}_{i,j} = F_{i,j}^{(s,s)}$.

See Fig. \ref{FigRK2} for the schematic of DIRK2 based scheme.

\subsubsection*{BDF methods.} 
The backward differentiation formula (BDF) is a linear multi step method (see \cite{HW}), which gives uniformly accurate solution for all values of $\kappa$. The $s$ order BDF method is represented by
\begin{align}\label{BDFscheme}
\begin{array}{l}
\displaystyle	BDF: y^{n+1} = \sum_{k=1}^{s} \alpha_k y^{n+1-k} +  \beta_s \Delta t g(y^{n+1},t_{n+1}).
\end{array}
\end{align} 
where $\alpha_k$ and $\beta_s$ are constants depending on $s$. In particular, we here consider BDF2 $(s=2)$ and BDF3 $(s=3)$ methods:
\begin{align}\label{BDFschemes}
\begin{array}{l}
\displaystyle	BDF2: y^{n+1}= \frac{4}{3}y^n -\frac{1}{3}y^{n-1} + \frac{2}{3} \Delta t g(y^{n+1},t_{n+1}).
,\\[3mm]
BDF3:
\displaystyle	y^{n+1}= \frac{18}{11}y^n - \frac{9}{11}y^{n-1} + \frac{2}{11}y^{n-2} + \frac{6}{11} \Delta t g(y^{n+1},t_{n+1}).
\end{array}
\end{align} 
The main advantage of $s$ order BDF method against a $s$-stage DIRK method is that it enables to design semi-Lagrangian schemes with less interpolations and computations of Maxwellian. Now, we illustrate BDF based semi-Lagrangian methods.
We say $x_{i,j}^k:=x_i-k v_j \Delta{t}$ in the following algorithm.

\subsubsection*{Algorithm for s-order BDF based methods} 
\begin{enumerate}
	\item For $k=1,\dots,s$, interpolate $f_{i,j}^{n,k}= f(x_{i,j}^k, v_j, t^{n+1-k})$  from $\left\{f_{i,j}^{n+1-k}\right\}_{i\in \mathcal{I}}$.
	\item
	Compute: 
	\[
	f_{i,j}^{n+1}= \sum_{k=1}^{s} \alpha_k f_{i,j}^{n,k} +  \beta_s \frac{\Delta t}{\kappa} \left( \mathcal{M}\left[f_{i,j}^{n+1}\right] - f_{i,j}^{n+1} \right),
	\]
	where we use $\mathcal{M}\left[f_{i,j}^{n+1}\right]=\mathcal{M}\left[\sum_{k=1}^{s} \alpha_k f_{i,j}^{n,k}\right].$
\end{enumerate}
\begin{figure}[t]
	\centering
	\includegraphics[width=0.4\linewidth]{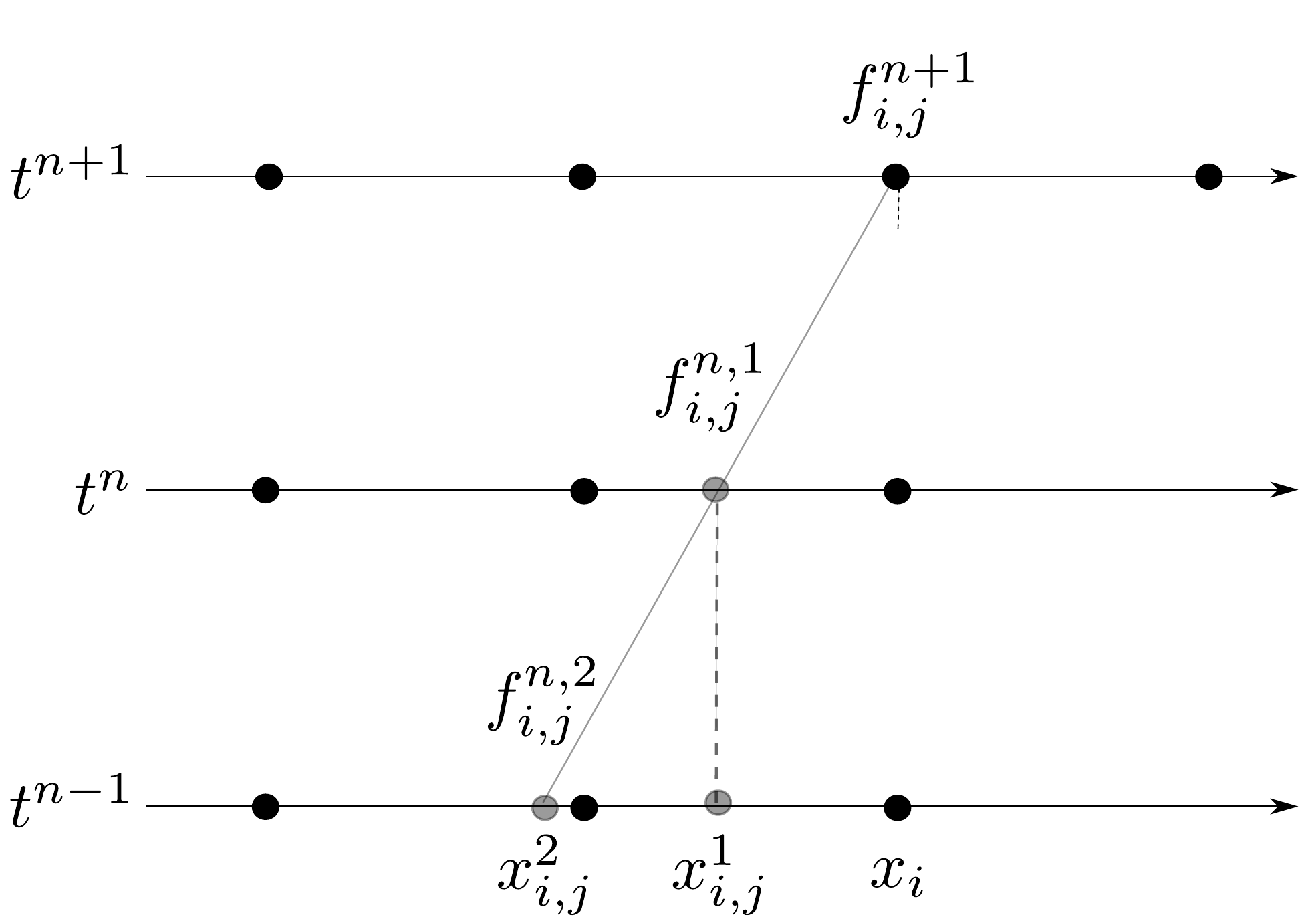}			\caption{Representation of BDF2 scheme. Black circles: grid nodes, grey circles: points where interpolation is needed.}	
	\label{Figbdf2}
\end{figure}
%

%
See Fig. \ref{Figbdf2} for the schematic of BDF2 based scheme.

\begin{remark}
	We remark that the positive preserving property of semi-Lagrangian schemes for BGK model is not straightforward even if we impose the positiveness for spatial reconstructions. This is mainly due to the use of high order time discretization such as DIRK and BDF methods. In particular, for DIRK methods, the values of $K_{i,j}^{(k,\ell)} = \frac{1}{\kappa}\big( \mathcal{M}(F_{i,j}^{(k,\ell)})-F_{i,j}^{(k,\ell)} \big)$ can be negative. In case of BDF schemes, the negative values of $\alpha_k$ in \eqref{BDFscheme} may lead to negative solutions. 
	For these reasons, in this paper, we only impose conservation and non-oscillatory properties. (For a positive scheme to the BGK model in the Eulerian framework, we refer to \cite{hu2018asymptotic}.)
\end{remark}

\subsubsection{Asymptotic preserving property of SL schemes for the BGK model}

In \cite{BCRY}, we show that the first order scheme SL \eqref{first order scheme} with time discretization gives
$f=\mathcal{M}(f)+\mathcal{O}(\kappa)$. The result can be easily extended to fully discretized high order SL schemes. Consider the BDF3 based scheme:
\begin{align}\label{f star}
f_{i,j}^{n+1}= f_{i,j}^* + \frac{6}{11} \frac{\Delta t}{\kappa} \left( \mathcal{M}[f]_{i,j}^* - f_{i,j}^{n+1} \right),
\end{align}
where $f_{i,j}^*=\frac{18}{11} f_{i,j}^{n,1} -\frac{9}{11} f_{i,j}^{n,2} +\frac{2}{11} f_{i,j}^{n,3}$ and  $\mathcal{M}[f]_{i,j}^*$ is the Maxwellian constructed from $\{f_{i,j}^*\}_{i \in \mathcal{I}}$. Then, we can rewrite \eqref{f star} as
\[
f_{i,j}^{n+1}= \frac{\kappa}{\kappa + \frac{6}{11}\Delta t}f_{i,j}^* + \frac{\frac{6}{11}\Delta t}{\kappa + \frac{6}{11}\Delta t}\mathcal{M}[f]_{i,j}^*.
\]
Under the assumption that $f_{i,j}^*,\,\mathcal{M}[f]_{i,j}^* = \mathcal{O}(1)$,
for a fixed time step $\Delta t>0$, the limit $\kappa \rightarrow 0$ implies that
\begin{align}\label{f mstar}
f_{i,j}^{n+1}= \mathcal{M}[f]_{i,j}^* + \mathcal{O}(\kappa).
\end{align}
for any $n\geq 0$ regardless of initial data. We further note that the following identity holds:
\begin{align}\label{m mstar}
\mathcal{M}[f]_{i,j}^{n+1} = \mathcal{M}[f]_{i,j}^*,
\end{align} 
which comes from the fact that $f_{i,j}^{n+1}$ and $f_{i,j}^*$ have the same moments.
Consequently, combining \eqref{f mstar} and \eqref{m mstar}, we obtain 
\begin{align}\label{f to M}
f_{i,j}^{n+1}= \mathcal{M}[f]_{i,j}^{n+1} + \mathcal{O}(\kappa).
\end{align}
In the similar manner, DIRK based high order semi-Lagrangian methods satisfy the relation \eqref{f to M}, which is referred to \textit{Asymptotic preserving property}. 
This related numerical tests will be  in sects \ref{sec 1d riemann} and \ref{sec 2d riemann}.

\begin{section}{Numerical tests}
	
	In this section, we test the performance of the proposed conservative high order semi-Lagrangian methods on  rigid body rotation \eqref{rigid}, Vlasov-Poisson system \eqref{VP} and BGK model \eqref{bgk}.
	
	For accuracy tests, we measure numerical errors and corresponding convergence rate using a relative $L^1$-norm of the numerical solutions. For example, in case of 1D a space-time dependent solution $u(x,t)$, let $u_i^n(h)$ be the solution obtained by using $h=\Delta x$ with $N_x$ grid points. Similarly, we denote by $u_i^n(h/2)$ the solution obtained with $h/2$ with $2N_x$ grid points.
	\begin{itemize}
		\item Relative $L^1$-norm of $u^n$ using $N_x$  
		\begin{align}\label{rel norm}
		\text{error$(N_x,2N_x)$} = \frac{\sum_{i=1}^{N_x^{d_x}} |u_i^n(h) - u_{2i-1}^{n}(h/2)|}{\sum_{i=1}^{N_x^{d_x}} |u_{2i-1}^{n}(h/2)|}
		\end{align}
		We measure two-dimensional relative $L^1$-norm in the same way.
		\item Convergence rate
		\begin{align}\label{rate}
		\text{Rate}(N_x)= \log_2\left(\frac{\text{error$(N_x,2N_x)$}}{\text{error$(2N_x,4N_x)$}}\right)
		\end{align}
	\end{itemize}
	
	For two-dimensional problems, we assume uniform square meshes both in space and velocity domains, i.e., $\Delta x=\Delta y=h$ for space and $\Delta v^1=\Delta v^2=\Delta v$ for velocity. In two-dimensional accuracy tests, we will compute errors and convergence rates using relative $L^1$-norm as in \eqref{rel norm} and \eqref{rate}.


	\subsection{Rigid body rotation}\label{sec rigid}
	
	Since the analytical solution for rigid body rotation \eqref{rigid} preserves the shape, we designed high order conservative splitting methods to have both non-oscillatory and MPP properties.
	
	For numerical tests, we take CWENO23, CWENO35 and CWENO47 reconstructions as basic reconstructions with a MPP limiter \eqref{limiter}. {(For details of CWENO, we refer to \cite{LPR1,LPR2,C-2008,CPSV}.)} For comparison, as counterparts, we consider conservative polynomials of degree 3, 5 and 7 as basic reconstructions with a MPP limiter.
	Note that such polynomials of degree 3, 5 and 7 are the same optimal polynomials used for CWENO23, CWENO35 and CWENO47 reconstructions.
	
	The schemes compared in this test are reported in Table \ref{tab rigid name}. 
	\begin{center}
		\begin{table}[ht]
			\centering
			{\begin{tabular}{|ccc|}
					\hline
					\multicolumn{1}{ |c }{}&
					\multicolumn{1}{ |c| }{Order of time splitting method}& \multicolumn{1}{ c|  }{Basic reconstruction} \\ \hline
					\multicolumn{1}{ |c| }{2-P3-MPP}&
					\multirow{2}{*}{2nd order Strang splitting \eqref{strang}}
					&\multicolumn{1}{ |c|  }{P3}    
					\\
					\cline{1-1}
					\cline{3-3}
					\multicolumn{1}{ |c| }{2-QCWENO23-MPP}&
					&\multicolumn{1}{ |c|  }{CWENO23}
					\\
					\hline
					\multicolumn{1}{ |c |}{4-P5-MPP}&
					\multirow{2}{*}{4th order Yoshida splitting \eqref{yoshida 4}}
					&\multicolumn{1}{ |c|  }{P5}    
					\\
					\cline{1-1}
					\cline{3-3}
					
					\multicolumn{1}{ |c| }{4-QCWENO35-MPP}&
					&\multicolumn{1}{ |c|  }{CWENO35}
					\\
					\hline
					\multicolumn{1}{ |c| }{6-P7-MPP}&
					\multirow{2}{*}{6th order Yoshida splitting \eqref{yoshida 6}}
					&\multicolumn{1}{ |c|  }{P7}    
					\\
					\cline{1-1}
					\cline{3-3}
					\multicolumn{1}{ |c| }{6-QCWENO47-MPP}&
					&\multicolumn{1}{ |c|  }{CWENO47}
					\\
					\hline
					\multicolumn{1}{ |c }{2D-P3-non-splitting}&
					\multicolumn{1}{ |c|  }{Exact characteristics (Eq.\ref{eq:exact_char})}     
					&\multicolumn{1}{ |c|  }{2D-P3}
					\\
					\hline
			\end{tabular}}
			\caption{Numerical methods used for rigid body rotation.}\label{tab rigid name}
		\end{table}
	\end{center}
	Here 2D-P3 is obtained by adopting, as a basic reconstruction, the two dimensional optimal polynomial of degree 2 used for 2D CWENO23 reconstruction \cite{LPR1}. Note that only 2D-P3-non-splitting does not use the MPP limiter.
	
	Here the CFL number is defined as 			
	\begin{align}\label{CFL rigid}
	\text{CFL}= \Delta t \left(\frac{|y|_{max}}{\Delta x} + \frac{|x|_{max}}{\Delta y}\right).
	\end{align}
	In the following examples, the computational domain is the square  $\Omega = [-\pi,\pi]^2$, therefore 
	$|x|_{max}=|y|_{max}=\pi$.
	
	\subsubsection{Accuracy test}
	To confirm the accuracy of the proposed schemes, we consider the following smooth initial data
	\begin{equation}
	u_0(x,y) = f\left(\sqrt{x^2 + (y-1.5)^2}\right), \quad x\in[-\pi,\pi],\quad y\in[-\pi,\pi].
	\label{eq:u0}
	\end{equation}
	where  
	\begin{align}\label{initial rigid acc}
	f(r) &= 
	\begin{cases}
	1, \quad\qquad  r < \frac{1}{4}\\
	p(r), \quad \frac{1}{4}\leq r \leq \frac{5}{4}\\
	0, \quad\qquad \text{otherwise} 
	\end{cases}, 	
	\end{align}
	Here $p:\mathbb{R} \rightarrow \mathbb{R}$
	is a polynomial of degree 16 determined by imposing 17 conditions, so that the piecewise polynomial function $f$ is of class $C^7(\mathbb{R})$.
	In particular, we impose that the first 8 derivatives of $p(r)$ are continuous at $r=\frac{1}{4}$ and that the first 7 derivatives are  continuous at $r=\frac{5}{4}$. We plot initial data in 
	Fig.~\ref{rigid_order_initial}. 
	
	\begin{figure}[ht]
		\centering
		\begin{subfigure}[b]{0.45\linewidth}
			\includegraphics[width=1\linewidth]{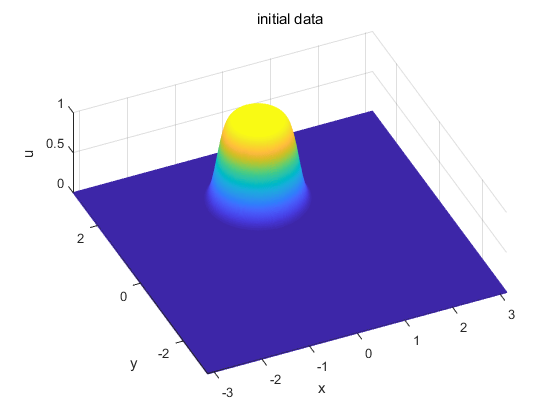}
		\end{subfigure}
		\caption{Accuracy test for 2D Rigid body motion. Initial data given in \eqref{initial rigid acc}.}\label{rigid_order_initial}
	\end{figure}

	For space, we assumed zero-boundary condition, and compute numerical solutions at $T_f=1$ with $N_x=N_y=40, 80, 160, 320, 640$. We choose the time step imposing that the CFL number  \eqref{CFL rigid} has a prescribed value. To obtain fast convergence to the desired accuracy, in the CWENO basic reconstructions we adopt $\varepsilon=\Delta x$. (See \cite{BCRY1}, Eq.~(18), and \cite{CPSV,kolb2014full} for the choice of $\varepsilon$.)
	
	In Table \ref{tab rigid body}, we report the relative $L^1$-norm \eqref{rel norm} and the corresponding convergence rates \eqref{rate} of the numerical solution $u$. The results imply that the use of MPP limiter \eqref{limiter} does not introduce any order reduction. For small CFL=$1.05$, spatial errors are dominant, which gives the desired accuracy. On the other hand, for large CFL$=8.4$, we observe the orders for time splitting methods which implies that time relevant errors are dominant.

	\begin{center}
		\begin{table}[ht]
			\centering
			{\begin{tabular}{|cccccc|}
					\hline
					\multicolumn{1}{ |c }{}&
					\multicolumn{1}{ |c| }{}& \multicolumn{2}{ c|  }{CFL$=1.05$} & \multicolumn{2}{ c|  }{CFL$=8.4$}  \\ \hline
					\multicolumn{1}{ |c }{}&
					\multicolumn{1}{ |c|  }{$(N_x,2N_x)$} &
					\multicolumn{1}{ c  }{error} &
					\multicolumn{1}{ c|  }{rate} &
					\multicolumn{1}{ c  }{error} &
					\multicolumn{1}{ c|  }{rate}   \\ 
					\hline
					\hline	
					\multicolumn{1}{ |c }{}&
					\multicolumn{1}{ |c|  }{$(40,80)$}&4.3025e-02      
					&2.12
					&2.4377e-01
					
					&1.98 
					
					\\
					\multicolumn{1}{ |c }{2-QCWENO23}&
					\multicolumn{1}{ |c|  }{$(80,160)$}&9.8567e-03      
					&2.30 
					&6.1939e-02
					
					&1.92 
					
					\\
					\multicolumn{1}{ |c }{MPP}&
					\multicolumn{1}{ |c|  }{$(160,320)$}&2.0039e-03      
					&2.97 
					&1.6316e-02
					
					&1.98 
					\\
					\multicolumn{1}{ |c }{}&
					\multicolumn{1}{ |c|  }{$(320,640)$}&
					2.5531e-04      
					&    &4.1422e-03
					&
					\\
					\hline	
					\multicolumn{1}{ |c }{}&
					\multicolumn{1}{ |c|  }{$(40,80)$}&4.2412e-02      
					&2.92 
					&1.4000e-01
					
					&4.26 
					
					\\
					\multicolumn{1}{ |c }{4-QCWENO35}&
					\multicolumn{1}{ |c|  }{$(80,160)$}&5.5900e-03      
					&4.06 
					&7.3184e-03
					
					&4.06 
					
					\\
					\multicolumn{1}{ |c }{MPP}&
					\multicolumn{1}{ |c|  }{$(160,320)$}&3.3613e-04      
					&4.94 
					&4.3762e-04
					
					&3.97 
					\\
					\multicolumn{1}{ |c }{}&
					\multicolumn{1}{ |c|  }{$(320,640)$}&
					1.0976e-05      
					&    &2.7920e-05
					&
					\\
					\hline	
					\multicolumn{1}{ |c }{}&
					\multicolumn{1}{ |c|  }{$(40,80)$}&2.3053e-02      
					&3.08 
					&3.0354e-01
					
					&8.10 
					
					\\
					\multicolumn{1}{ |c }{6-QCWENO47}&
					\multicolumn{1}{ |c|  }{$(80,160)$}&2.7278e-03      
					&5.31 
					&1.1058e-03
					
					&5.74 
					
					\\
					\multicolumn{1}{ |c }{MPP}&
					\multicolumn{1}{ |c|  }{$(160,320)$}&6.8834e-05      
					&6.92 
					&2.0731e-05
					
					&6.17 
					\\
					\multicolumn{1}{ |c }{}&
					\multicolumn{1}{ |c|  }{$(320,640)$}&
					5.6901e-07      
					&    &2.8779e-07
					&
					\\
					\hline
					\hline
			\end{tabular}}
			\caption{Accuracy test for splitting methods of rigid body rotation with initial data \eqref{initial rigid acc}.}\label{tab rigid body}
		\end{table}
	\end{center}

	\begin{center}
		\begin{table}[ht]
			\centering
			{\begin{tabular}{|cccc|}
					\hline
					\multicolumn{1}{ |c }{}&
					\multicolumn{1}{ |c| }{}& \multicolumn{1}{ c|  }{CFL$=1.05$} & \multicolumn{1}{ c|  }{CFL$=8.4$}  \\ \hline
					\multicolumn{1}{ |c }{}&
					\multicolumn{1}{ |c|  }{$N_x$} &
					\multicolumn{1}{ c|  }{$L^1$ error} &
					\multicolumn{1}{ c|  }{$L^1$ error}   \\ 
					\hline
					\hline	
					%
					\multicolumn{1}{ |c }{}&
					\multicolumn{1}{ |c|  }{$80$}&0.0000e-00
					&3.5927e-16

					\\
					\multicolumn{1}{ |c }{2-QCWENO23}&
					\multicolumn{1}{ |c|  }{$160$}&3.5927e-16           
					&1.7964e-16
					
					\\
					\multicolumn{1}{ |c }{MPP}&
					\multicolumn{1}{ |c|  }{$320$}      
					&1.7964e-16       &7.1854e-16
					
					\\\multicolumn{1}{ |c }{}&
					\multicolumn{1}{ |c|  }{$640$}      
					&3.5927e-16   
					&5.3891e-16
					\\
					\hline	
					%
					%
					%
					\multicolumn{1}{ |c }{}&
					\multicolumn{1}{ |c|  }{$80$}&1.9034e-12   
					&1.6167e-15

					\\
					\multicolumn{1}{ |c }{4-QCWENO35}&
					\multicolumn{1}{ |c|  }{$160$}&2.6406e-14              
					&3.0538e-15

					\\
					\multicolumn{1}{ |c }{MPP}&
					\multicolumn{1}{ |c|  }{$320$}      
					&5.2992e-14          &7.1854e-15

					\\\multicolumn{1}{ |c }{}&
					\multicolumn{1}{ |c|  }{$640$}      
					&1.0706e-13      
					&1.3473e-14
					\\
					\hline	
					%
					%
					%
					\multicolumn{1}{ |c }{}&
					\multicolumn{1}{ |c|  }{$80$}&1.4654e-11   
					&1.1856e-14

					\\
					\multicolumn{1}{ |c }{6-QCWENO47}&
					\multicolumn{1}{ |c|  }{$160$}&1.7964e-16              
					&7.1854e-16

					\\
					\multicolumn{1}{ |c }{MPP}&
					\multicolumn{1}{ |c|  }{$320$}      
					&1.4371e-15          &8.9818e-16

					\\\multicolumn{1}{ |c }{}&
					\multicolumn{1}{ |c|  }{$640$}      
					&3.5927e-15               
					&0.0000e-00
					\\
					\hline
					\hline
			\end{tabular}}
			\caption{Conservation errors obtained by splitting methods for rigid body rotation with initial data \eqref{initial rigid acc}.}\label{tab rigid body con}
		\end{table}
	\end{center}
	As a remark, we additionally present numerical results obtained by the non-splitting semi-Lagrangian method for rigid body rotation which makes use of exact characteristics. Here the only error (besides round off) is due to interpolation. In Tab.\ref{tab rigid body non-splitting}, we report numerical errors obtained with different time steps: (1) $\Delta t$ is determined by CFL=$1.05,\,8.4$ in \eqref{CFL rigid} (2) $\Delta t= T_f = 1$.
	For CFL$=1.05,/, 8.4$, the global accuracy is degraded by one as expected, hence it is close to 3. On the other hand, one-step based method ($\Delta t=1$ case) does not have an order reduction, so it gives order 4 which is the accuracy of QCWENO23 reconstruction for smooth functions.
	\begin{center}
		\begin{table}[ht]
			\centering
			{\begin{tabular}{|cccccccc|}
					\hline
					\multicolumn{1}{ |c }{}&
					\multicolumn{1}{ |c| }{}& \multicolumn{2}{ c|  }{CFL$=1.05$} & \multicolumn{2}{ c|  }{CFL$=8.4$} & \multicolumn{2}{ c|  }{$\Delta t=1$}  \\ \hline
					\multicolumn{1}{ |c }{}&
					\multicolumn{1}{ |c|  }{$(N_x,2N_x)$} &
					\multicolumn{1}{ c  }{error} &
					\multicolumn{1}{ c|  }{rate} &
					\multicolumn{1}{ c  }{error} &
					\multicolumn{1}{ c|  }{rate}&
					\multicolumn{1}{ c  }{error} &
					\multicolumn{1}{ c|  }{rate}   \\ 
					\hline
					\hline	
					\multicolumn{1}{ |c }{}&
					\multicolumn{1}{ |c|  }{$(40,80)$}&4.4767e-02        &2.1738            
					&1.5350e-02
					
					&2.9408
					&8.7760e-03 &3.5933
					\\
					\multicolumn{1}{ |c }{2D-P3}&
					\multicolumn{1}{ |c|  }{$(80,160)$}&9.9211e-03         
					&							2.7599            
					&1.9991e-03
					&3.2646
					&7.2713e-04 & 3.9797
					\\
					\multicolumn{1}{ |c }{non-splitting}&
					\multicolumn{1}{ |c|  }{$(160,320)$}&							1.4647e-03        
					&							2.9500    
					
					& 2.0802e-04

					&2.9769
					&4.6088e-05 & 3.9874
					\\
					\multicolumn{1}{ |c }{}&
					\multicolumn{1}{ |c|  }{$(320,640)$}&
					1.8955e-04   &    &2.6423e-05
					&
					&2.9058e-06&
					\\
					\hline
					\hline
			\end{tabular}}
			\caption{Accuracy test for 2D-P3-non-splitting method of rigid body rotation with initial data \eqref{initial rigid acc}.}\label{tab rigid body non-splitting}
		\end{table}
	\end{center}

	\begin{center}
		\begin{table}[ht]
			\centering
			{\begin{tabular}{|ccccc|}
					\hline
					\multicolumn{1}{ |c }{}&
					\multicolumn{1}{ |c| }{}& \multicolumn{1}{ c|  }{CFL$=1.05$} & \multicolumn{1}{ c|  }{CFL$=8.4$}  & \multicolumn{1}{ c|  }{$\Delta t=1$}   
					\\ \hline
					\multicolumn{1}{ |c }{}&
					\multicolumn{1}{ |c|  }{$N_x$} &
					\multicolumn{1}{ c|  }{mass error} &
					\multicolumn{1}{ c|  }{mass error} &
					\multicolumn{1}{ c|  }{mass error}   \\ 
					\hline
					\hline	
					\multicolumn{1}{ |c }{}&
					\multicolumn{1}{ |c|  }{$80$}&9.9358e-07   
					&2.7724e-06 & 5.7758e-07 
					\\
					\multicolumn{1}{ |c }{2D-P3}&
					\multicolumn{1}{ |c|  }{$160$}&2.5417e-08       	&1.2029e-07
					&9.5474e-10 
					\\
					\multicolumn{1}{ |c }{non-splitting}&
					\multicolumn{1}{ |c|  }{$320$}      
					&5.8763e-11          &8.3774e-10
					&8.1131e-11 
					\\\multicolumn{1}{ |c }{}&
					\multicolumn{1}{ |c|  }{$640$}      
					&5.6213e-12      
					&4.9656e-11
					&2.3963e-13
					\\
					\hline	
					
					\hline
					\hline
			\end{tabular}}
			\caption{Conservation errors in mass obtained by 2D-P3-non-splitting for rigid body rotation with initial data \eqref{initial rigid acc}.}\label{tab rigid body con non-splitting}
		\end{table}
	\end{center}

	\subsubsection{Maximum principle preserving and non-oscillatory properties}

	In this example, we show that the combined use of CWENO reconstructions and MPP limiter for the basic reconstruction enables to achieve both
	non-oscillatory property and the maximum principle preserving property at the same time.
	
	For this, we consider initial data composed of a slotted disk, a cone and a smooth hump, as in  \cite{qiu2011positivity}, together with a double layered disk. We plot initial data in 
	Fig.~\ref{rigid_sharp_inintial} and compute numerical solutions with CFL$=6.3$.
	\begin{figure}[ht]
		\centering
		\begin{subfigure}[b]{0.45\linewidth}
			\includegraphics[width=1\linewidth]{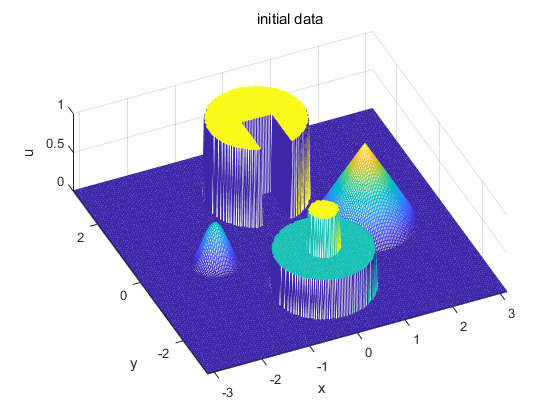}
		\end{subfigure}
		\caption{2D Rigid body motion. Initial data used for results in Figs. \ref{rigid QCWENO23 all}, \ref{rigid QCWENO35 all}, \ref{rigid QCWENO47 all}. }\label{rigid_sharp_inintial}
	\end{figure}
	In Figs. \ref{rigid QCWENO23 all}, \ref{rigid QCWENO35 all}, \ref{rigid QCWENO47 all}, we compare numerical methods listed in Tab. \ref{tab rigid name}. In particular, we plot contours and slides of numerical solutions at different locations at final time $T_f=2\pi$. We first confirm that the use of MPP limiter makes numerical solutions to be bounded by $[0,1]$ in Figs. \ref{rigid QCWENO23 Cylinder}, \ref{rigid QCWENO35 Cylinder}, \ref{rigid QCWENO47 Cylinder}. In Figs. \ref{rigid QCWENO23 Cone}, \ref{rigid QCWENO23 hump}, \ref{rigid QCWENO35 Cone}, \ref{rigid QCWENO35 hump}, \ref{rigid QCWENO47 Cone}, \ref{rigid QCWENO47 hump}, the use of non-linear weights make CWENO based methods to be more dissipative. Although such dissipative tendency can be improved for high order CWENO based methods, the use of CWENO reconstruction does not seem to be necessary for the numerical treatment of cone and hump. 
	
	On the other hand, in Figs. \ref{rigid QCWENO23 double layer}, \ref{rigid QCWENO35 double layer}, \ref{rigid QCWENO47 double layer}, one can see that the MPP limiter is not enough to prevent oscillations in a middle layer. However, the use CWENO as basic reconstruction maintains the  non-oscillatory behaviour of the solution in the middle layer. Consequently, the technique based on CWENO with MPP limiter makes provides a high order conservative reconstruction which avoids spurious oscillation and maintains the maximum principle preserving property at the same time.
	
	We remark that the careful choice of $\varepsilon$ or non-linear weights for CWENO is necessary to inherit its non-oscillatory property. In this numerical test, we simply set
	$\varepsilon=10^{-4}$ for all CWENO reconstructions. For too large value of $\varepsilon$, non-linear weights reduce to linear weights, which do not prevent oscillations,  while if  $\varepsilon$ is too small,  the resulting reconstruction can be too dissipative. For the optimal choice of $\varepsilon$ and non-linear weights for CWENO reconstruction, we refer to \cite{CPSV-2017,CPSV,kolb2014full}.
	\begin{figure}[htbp]
		\centering
		\begin{subfigure}[b]{0.45\linewidth}
			\includegraphics[width=1\linewidth]{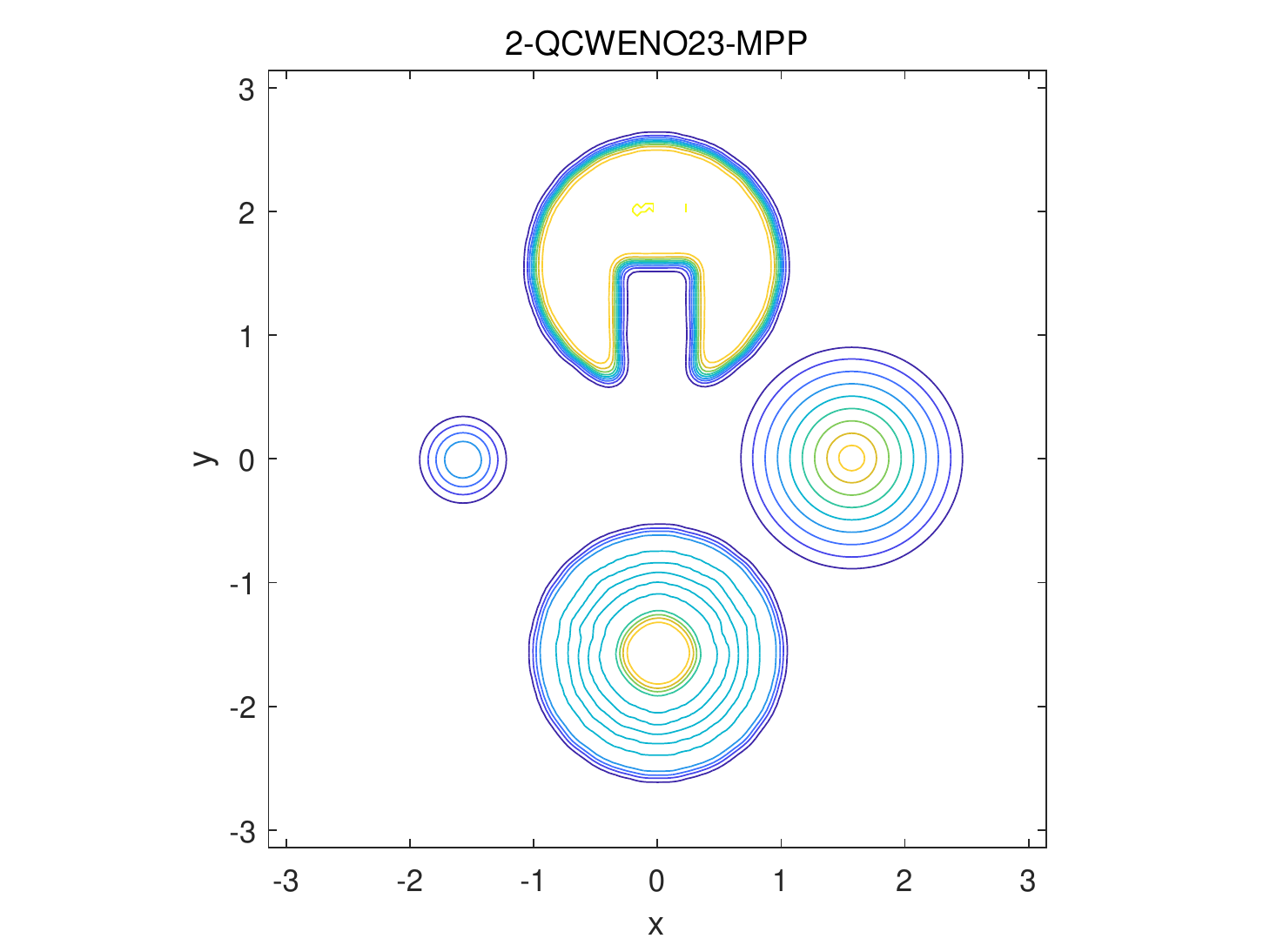}
			\subcaption{Contour}
		\end{subfigure}
		\vspace*{5mm}
		\begin{subfigure}[b]{0.45\linewidth}
			\includegraphics[width=1\linewidth]{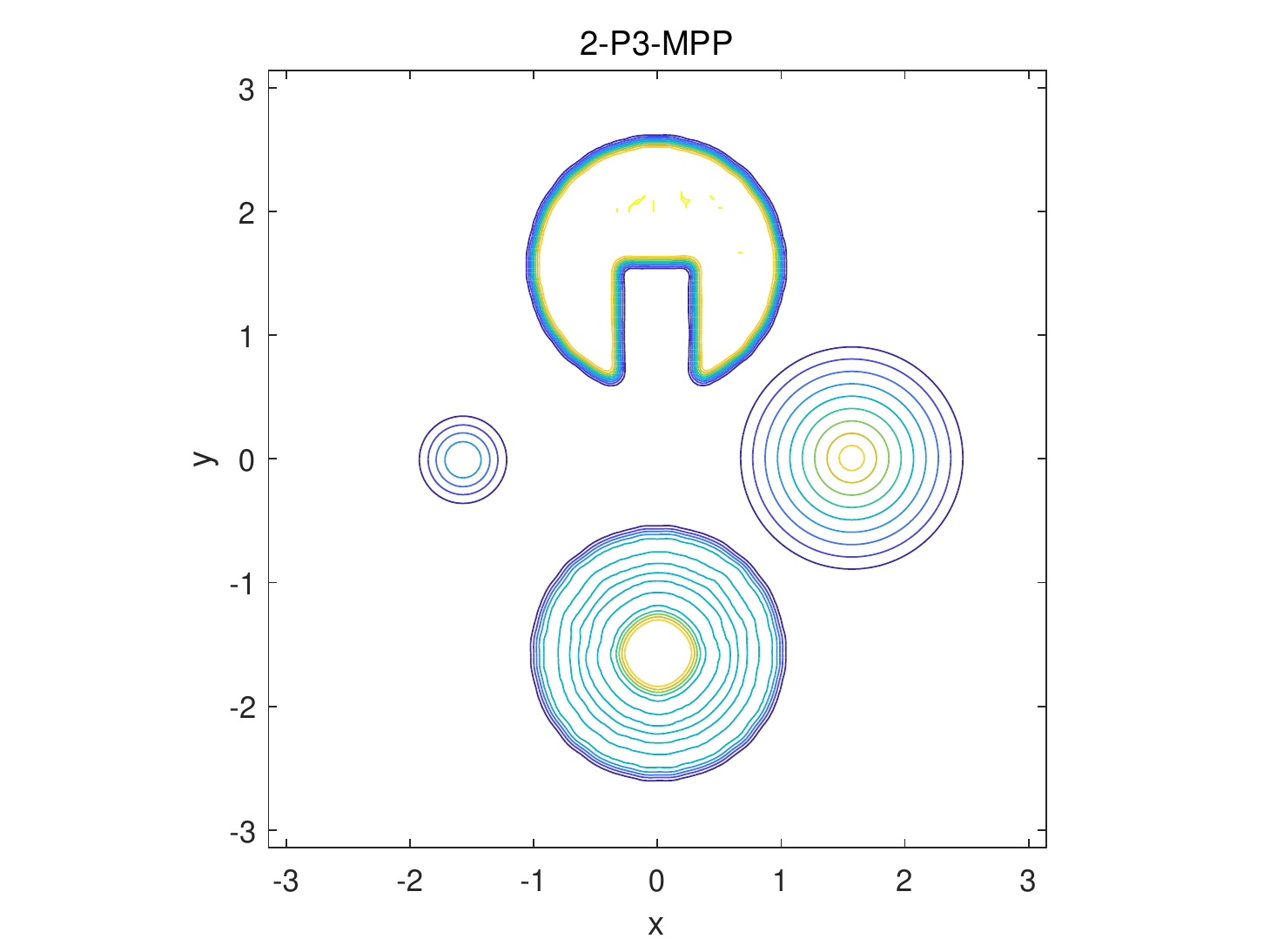}
			\subcaption{Contour}
		\end{subfigure}
		\vspace*{5mm}
		\begin{subfigure}[b]{0.45\linewidth}
			\includegraphics[width=1\linewidth]{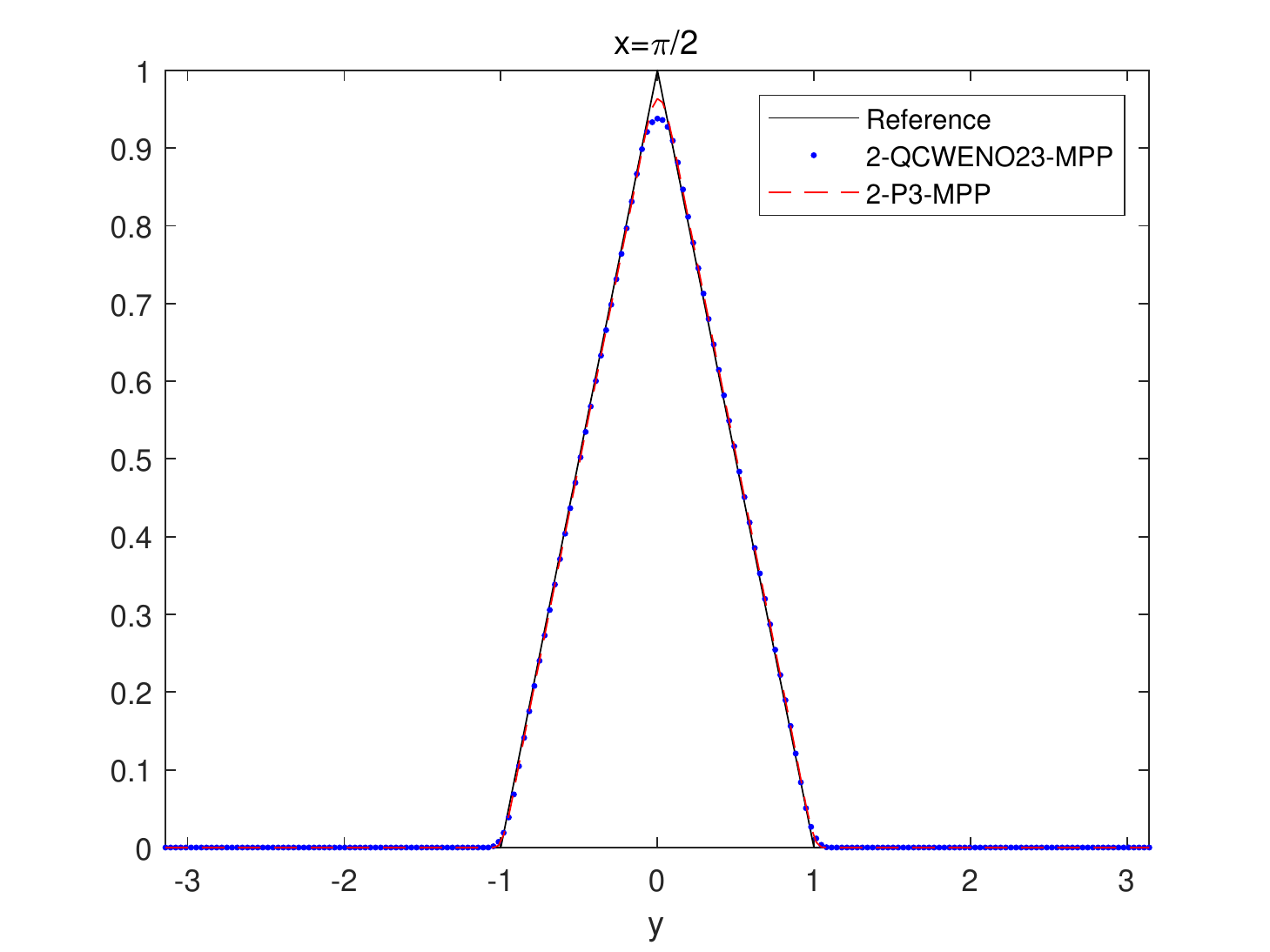}
			\subcaption{Cone}\label{rigid QCWENO23 Cone}
		\end{subfigure}	
		\begin{subfigure}[b]{0.45\linewidth}
			\includegraphics[width=1\linewidth]{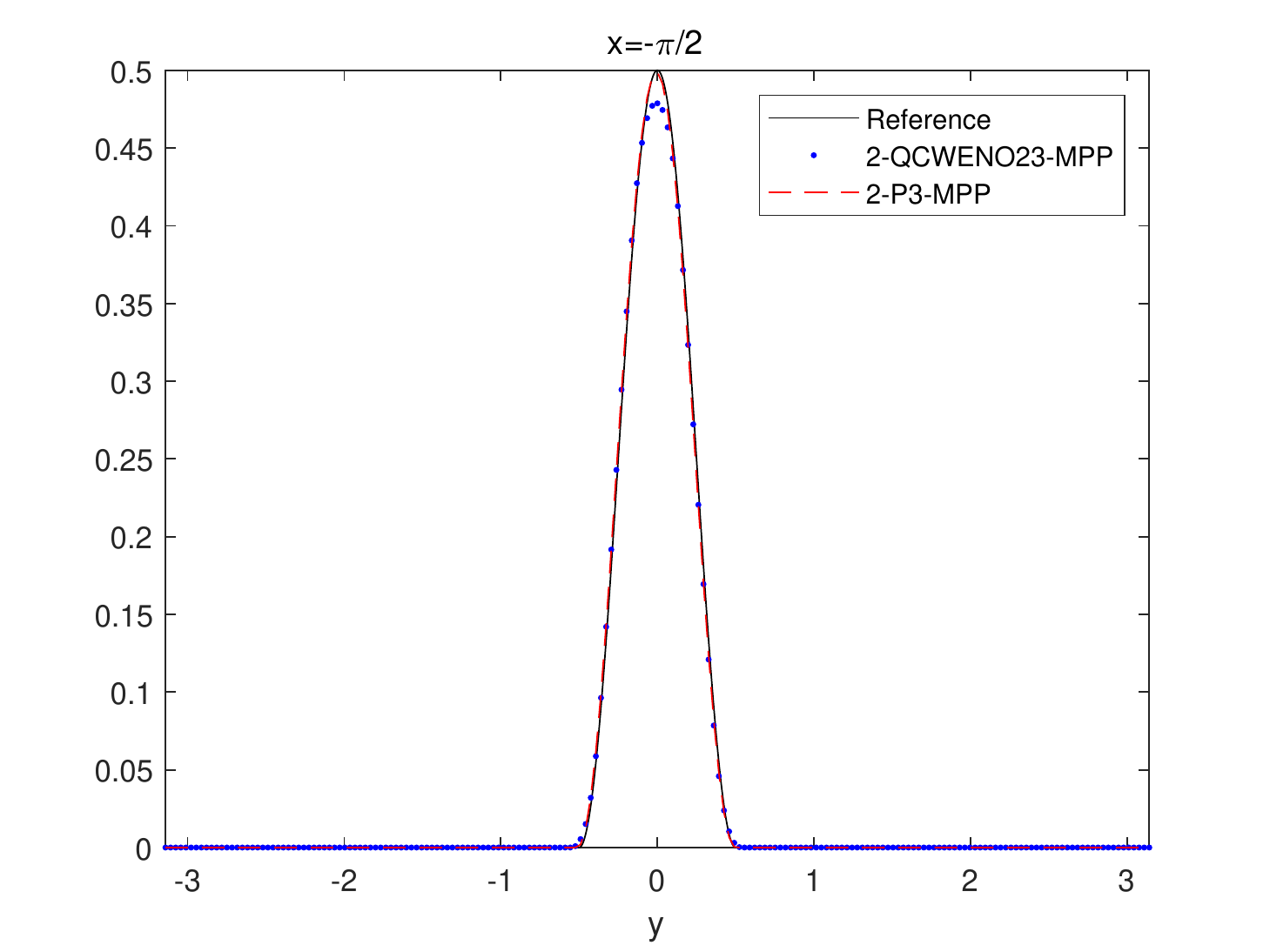}
			\subcaption{Humb}\label{rigid QCWENO23 hump}
		\end{subfigure}	
		\begin{subfigure}[b]{0.45\linewidth}
			\includegraphics[width=1\linewidth]{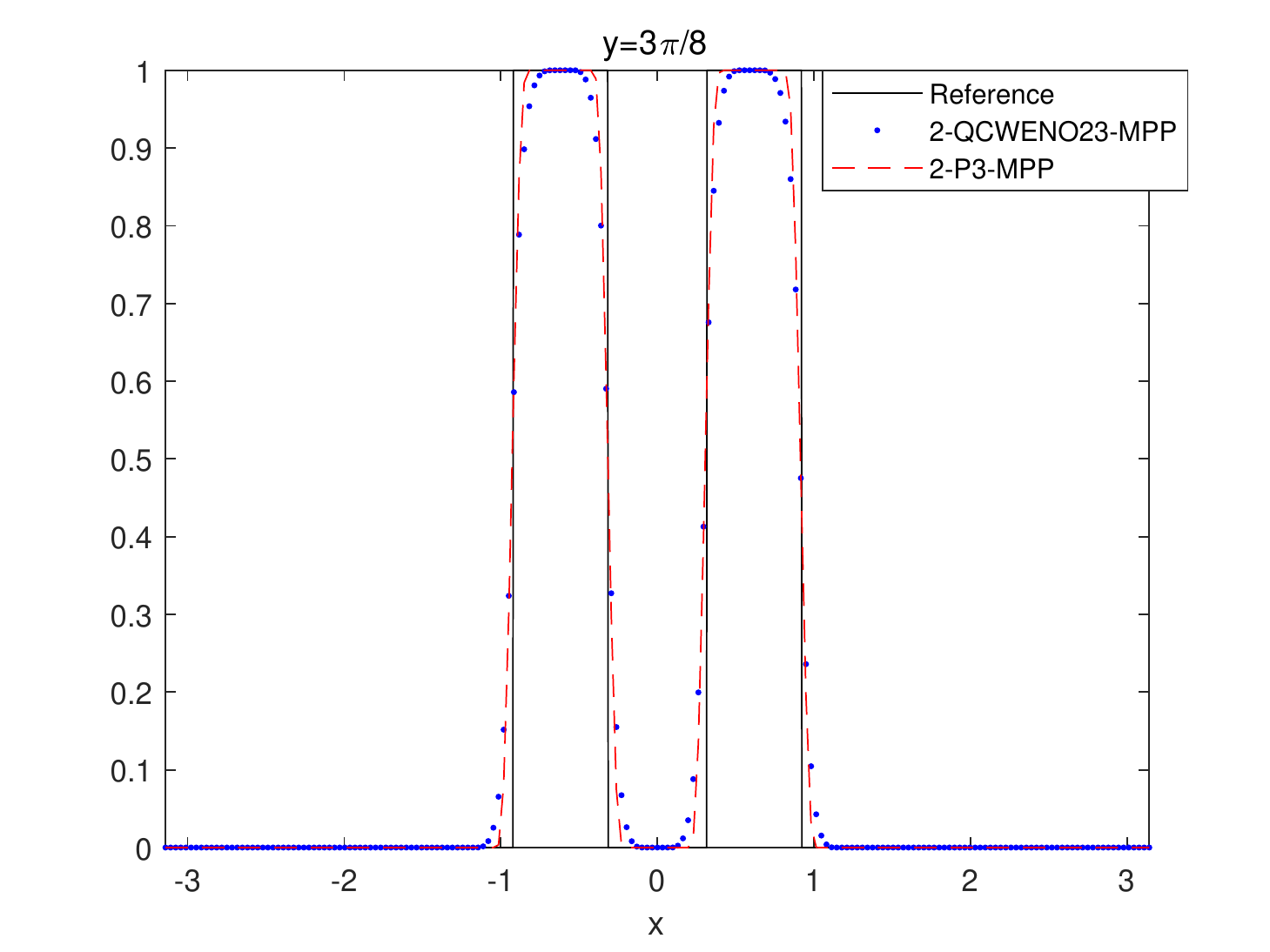}
			\subcaption{Cylinder}\label{rigid QCWENO23 Cylinder}
		\end{subfigure}	
		\begin{subfigure}[b]{0.45\linewidth}
			\includegraphics[width=1\linewidth]{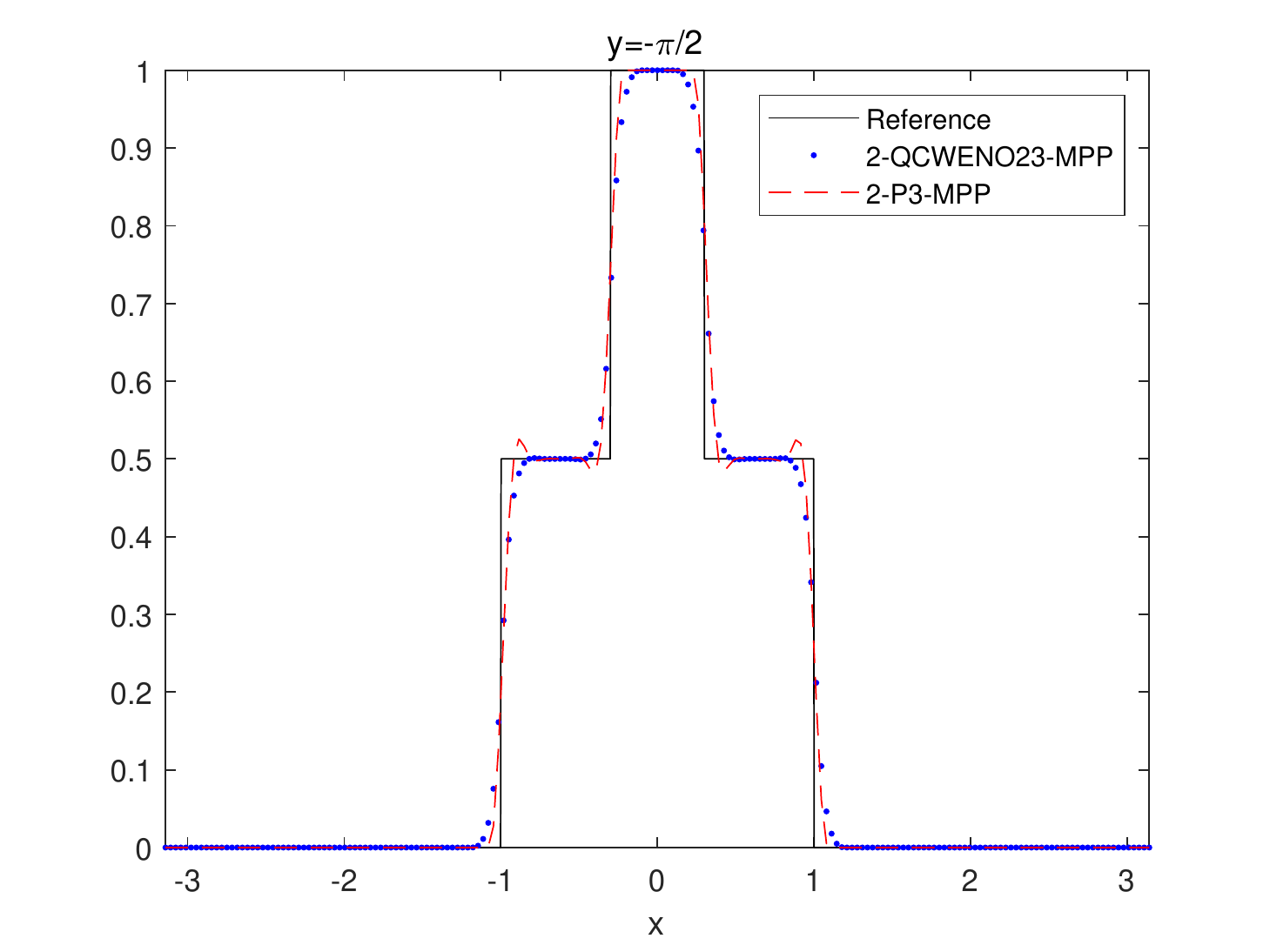}
			\subcaption{Double layer}\label{rigid QCWENO23 double layer}
		\end{subfigure}					
		\caption{2D Rigid body motion. Numerical solutions obtained for $(N_x,N_y)=192 \times 192$. For CWENO23, we use $\varepsilon=10^{-4}$}\label{rigid QCWENO23 all}
	\end{figure}
	\begin{figure}[htbp]
		\centering
		\begin{subfigure}[b]{0.45\linewidth}
			\includegraphics[width=1\linewidth]{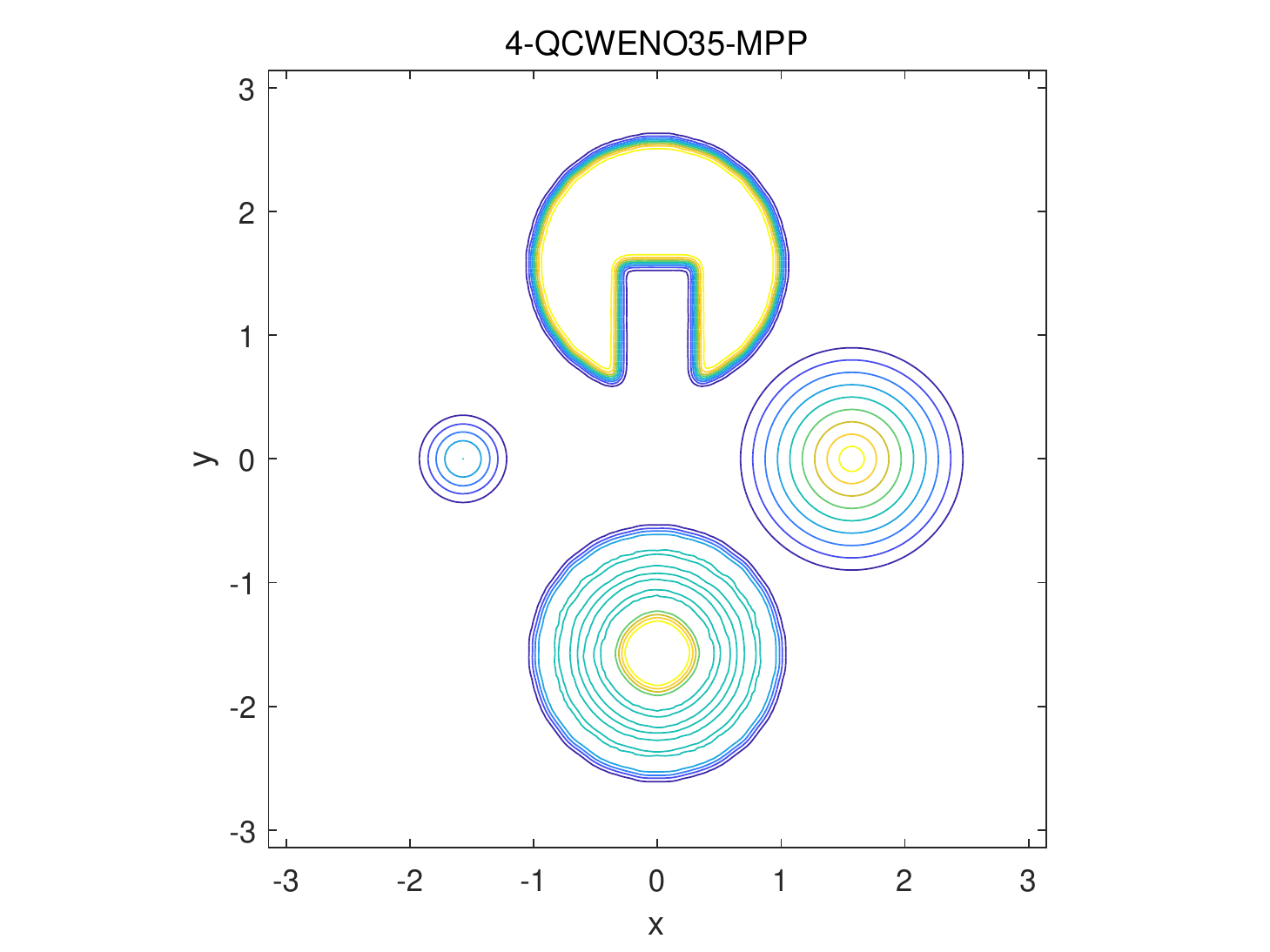}
			\subcaption{Contour}
		\end{subfigure}
		\vspace*{5mm}
		\begin{subfigure}[b]{0.45\linewidth}
			\includegraphics[width=1\linewidth]{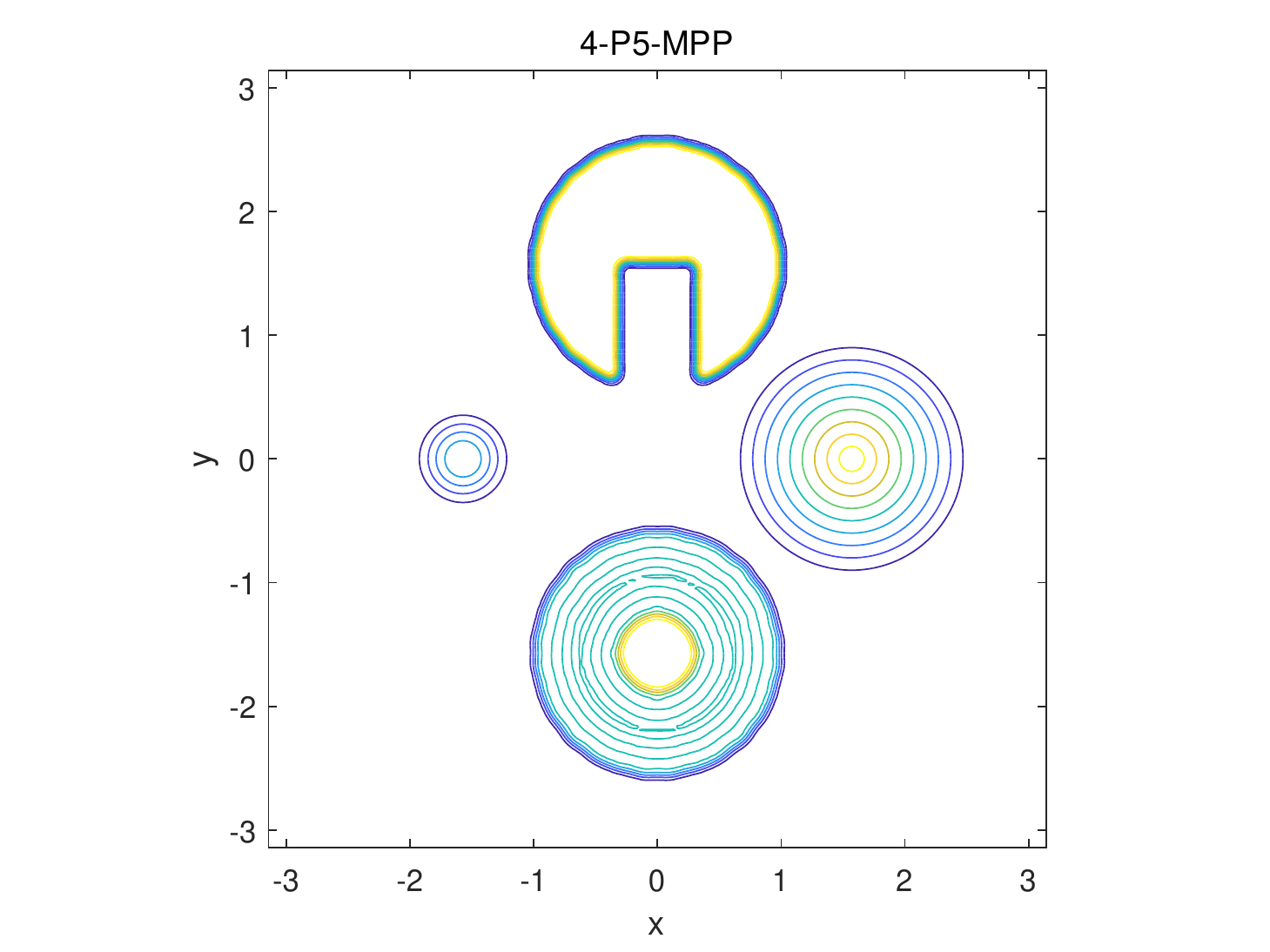}
			\subcaption{Contour}
		\end{subfigure}
		\vspace*{5mm}
		\begin{subfigure}[b]{0.45\linewidth}
			\includegraphics[width=1\linewidth]{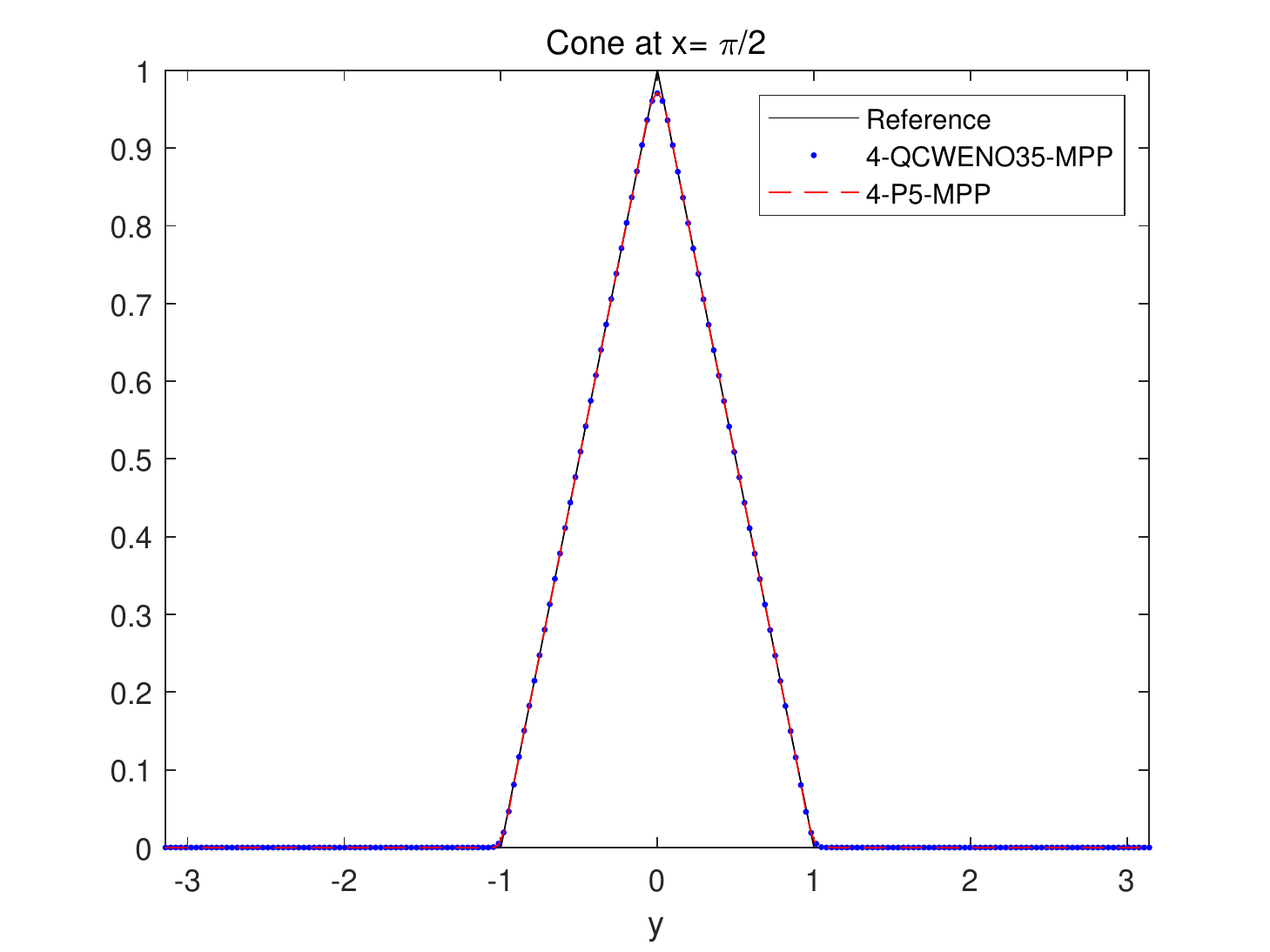}
			\subcaption{Cone}\label{rigid QCWENO35 Cone}
		\end{subfigure}	
		\begin{subfigure}[b]{0.45\linewidth}
			\includegraphics[width=1\linewidth]{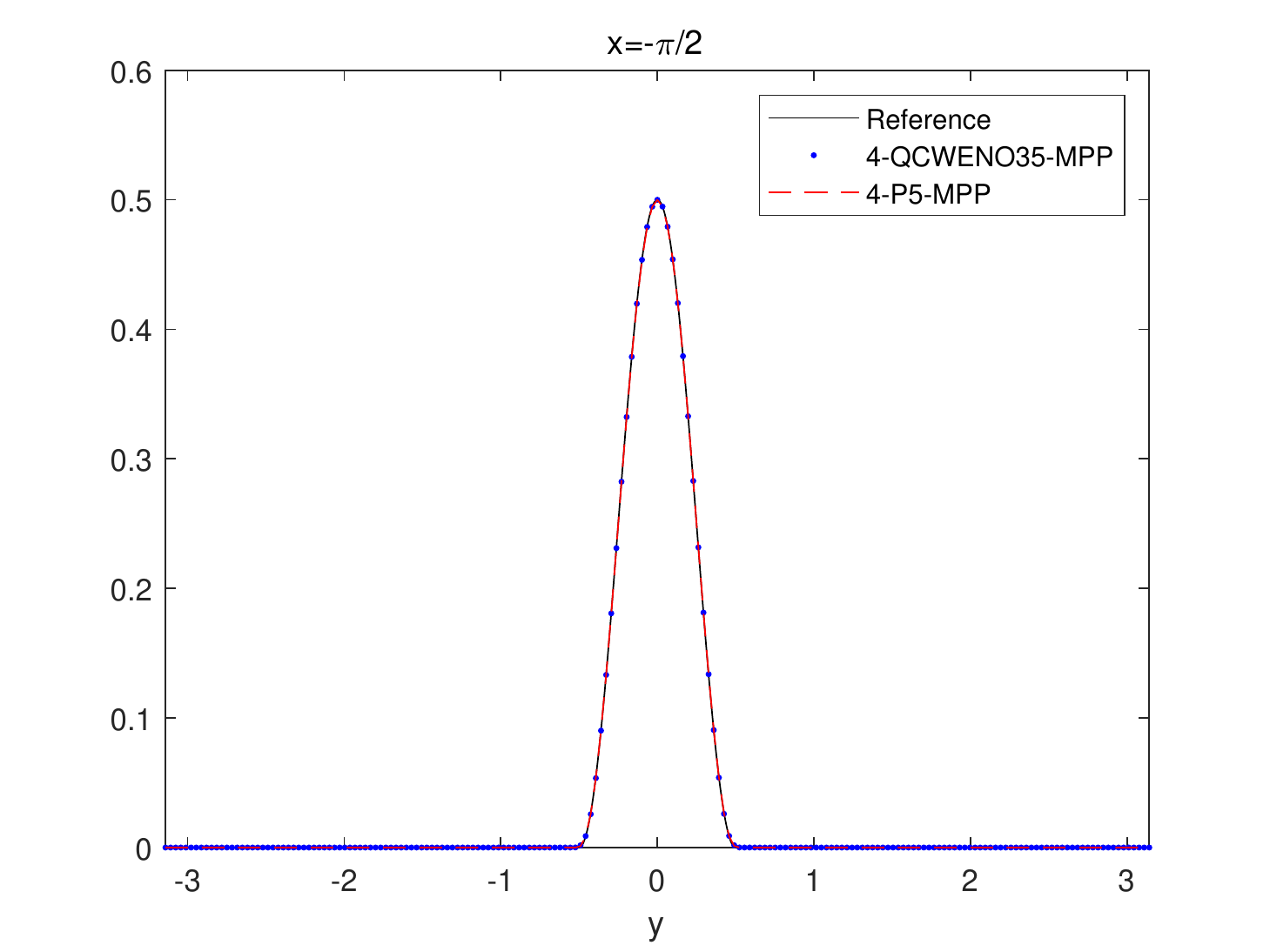}
			\subcaption{Humb}\label{rigid QCWENO35 hump}
		\end{subfigure}	
		\begin{subfigure}[b]{0.45\linewidth}
			\includegraphics[width=1\linewidth]{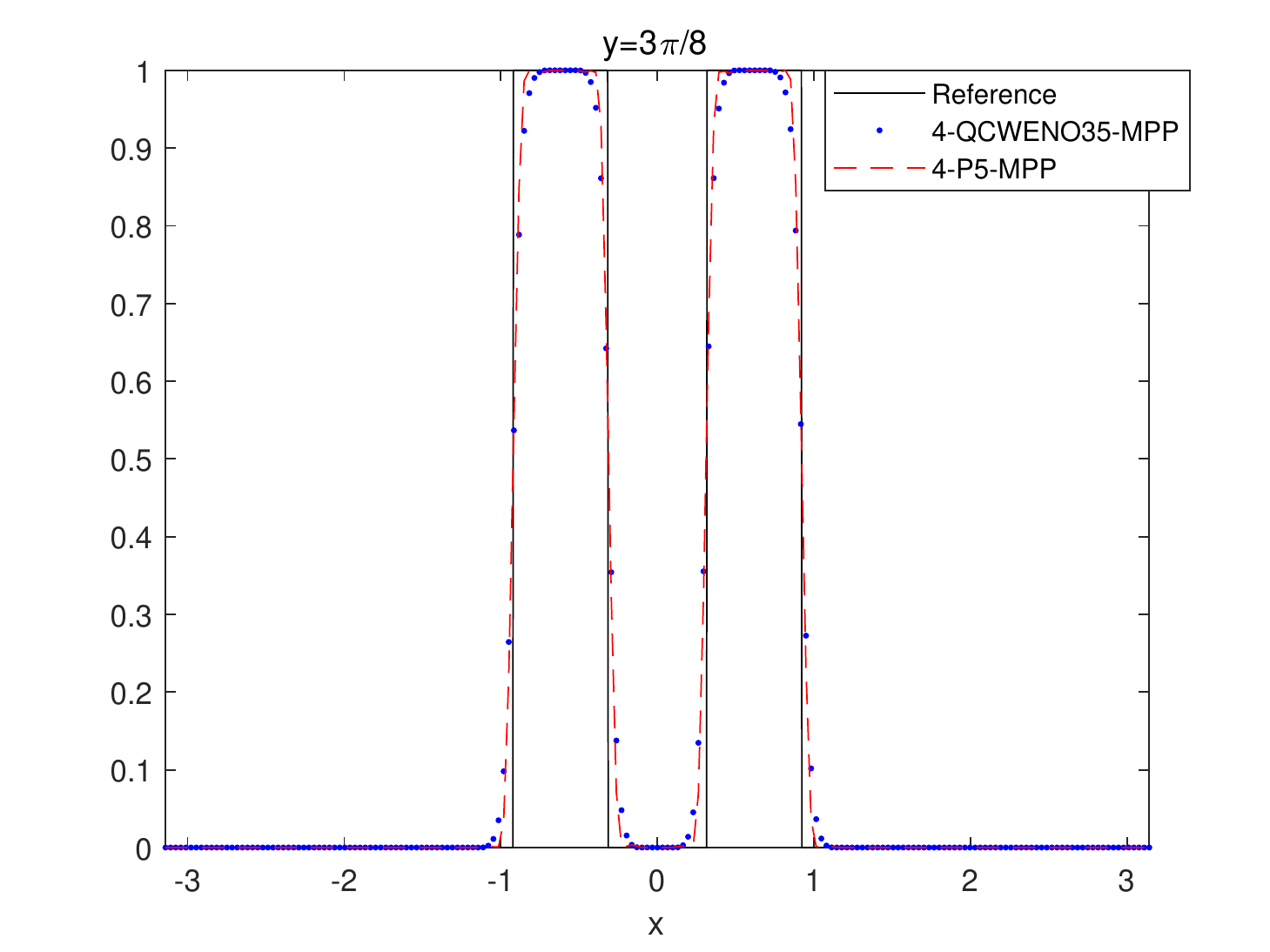}
			\subcaption{Cylinder}\label{rigid QCWENO35 Cylinder}
		\end{subfigure}	
		\begin{subfigure}[b]{0.45\linewidth}
			\includegraphics[width=1\linewidth]{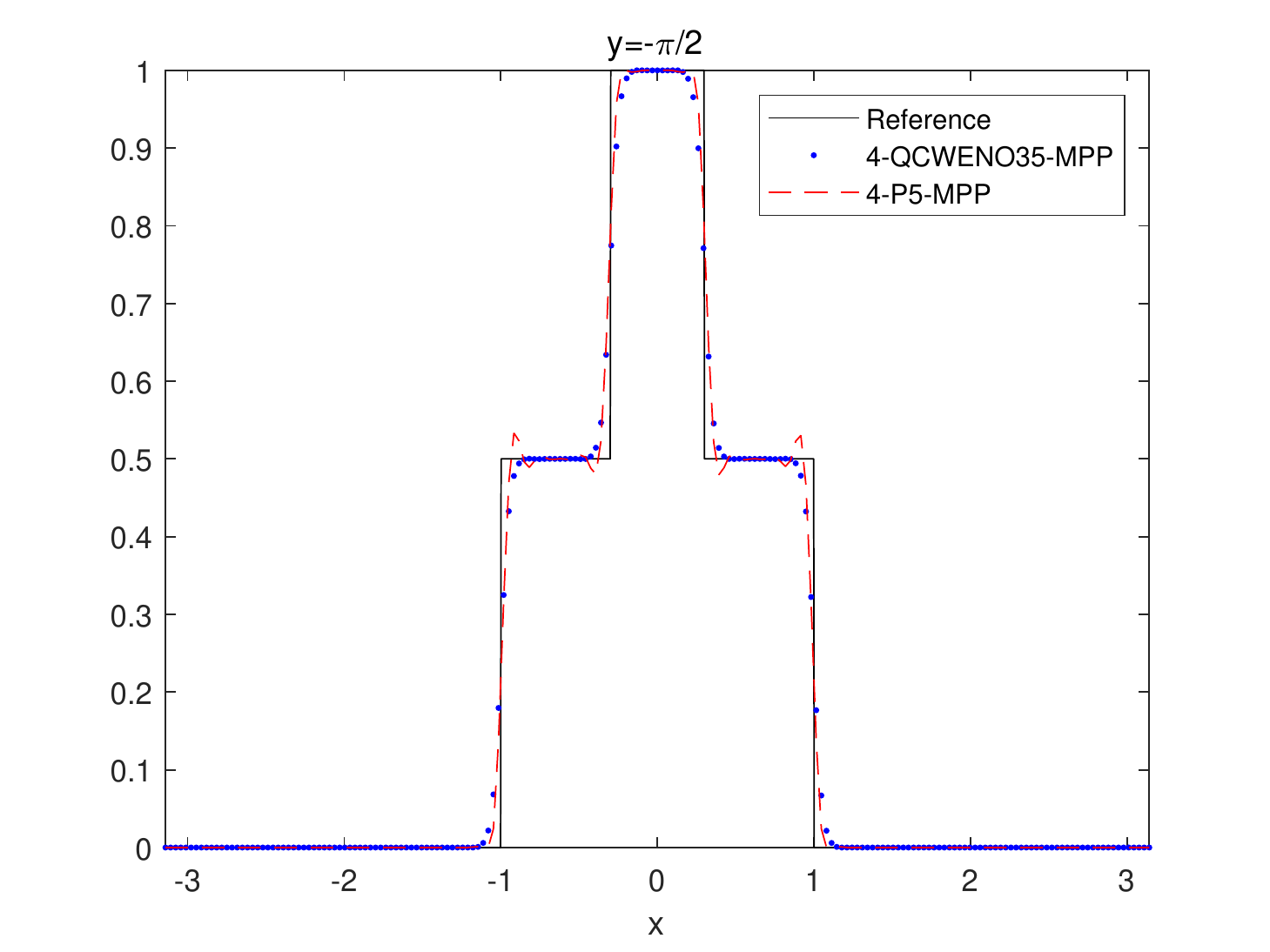}
			\subcaption{Double layer}\label{rigid QCWENO35 double layer}
		\end{subfigure}					
		\caption{2D Rigid body motion. Numerical solutions obtained for $(N_x,N_y)=192 \times 192$. For CWENO35, we use $\varepsilon=10^{-4}$}\label{rigid QCWENO35 all}
	\end{figure}
	\begin{figure}[htbp]
		\centering
		\begin{subfigure}[b]{0.45\linewidth}
			\includegraphics[width=1\linewidth]{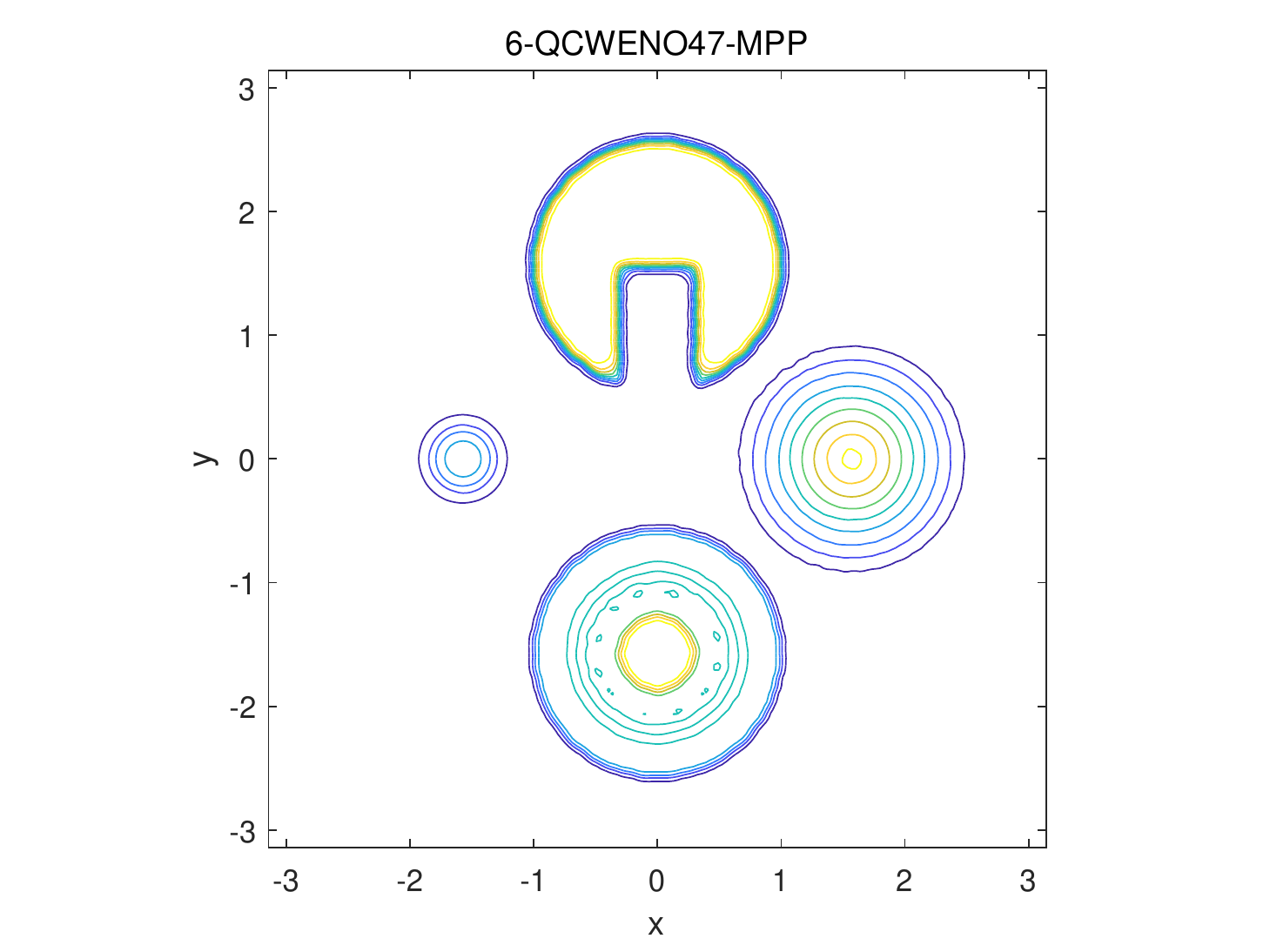}
			\subcaption{Contour}
		\end{subfigure}
		\vspace*{5mm}
		\begin{subfigure}[b]{0.45\linewidth}
			\includegraphics[width=1\linewidth]{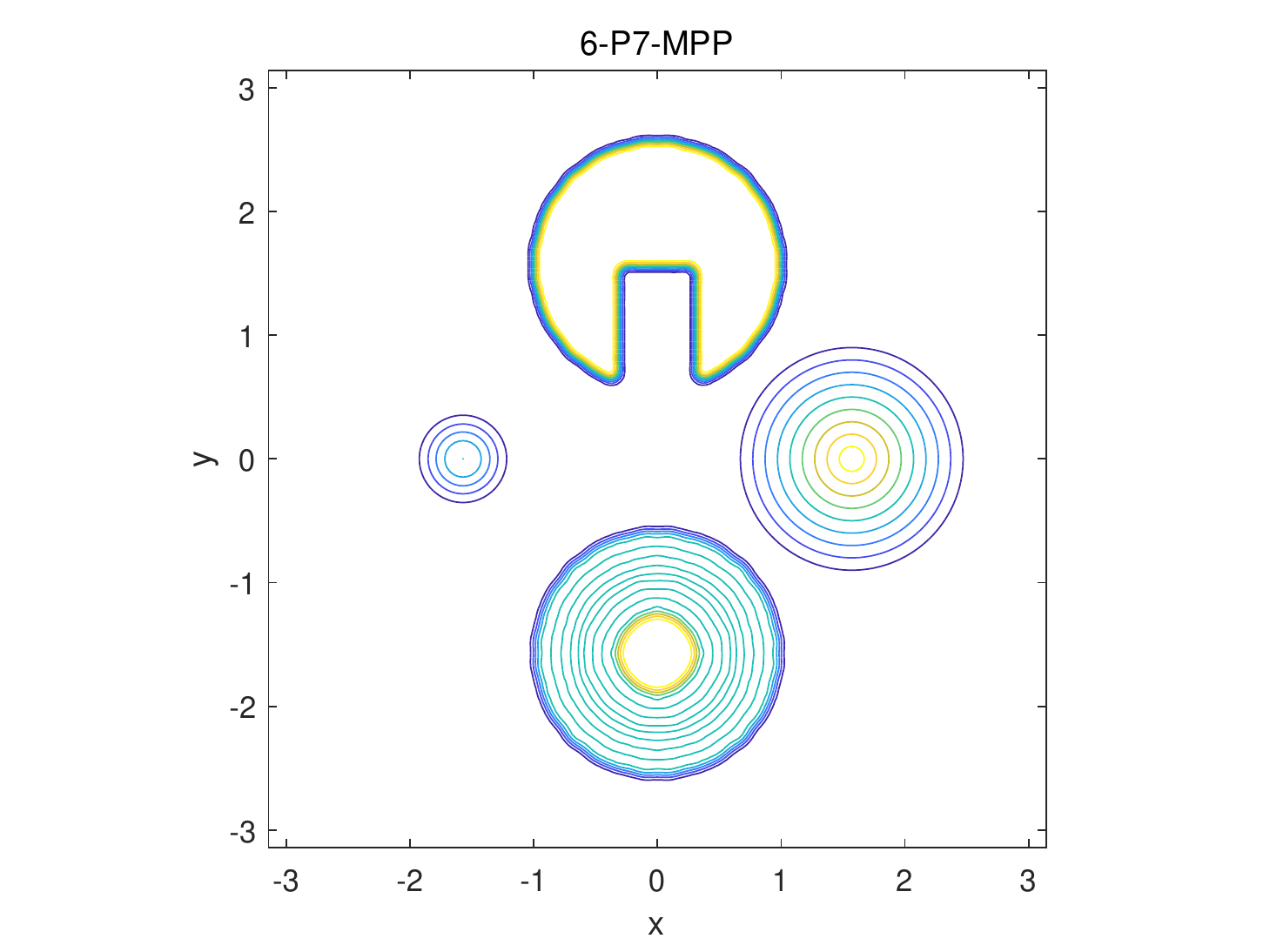}
			\subcaption{Contour}
		\end{subfigure}
		\vspace*{5mm}
		\begin{subfigure}[b]{0.45\linewidth}
			\includegraphics[width=1\linewidth]{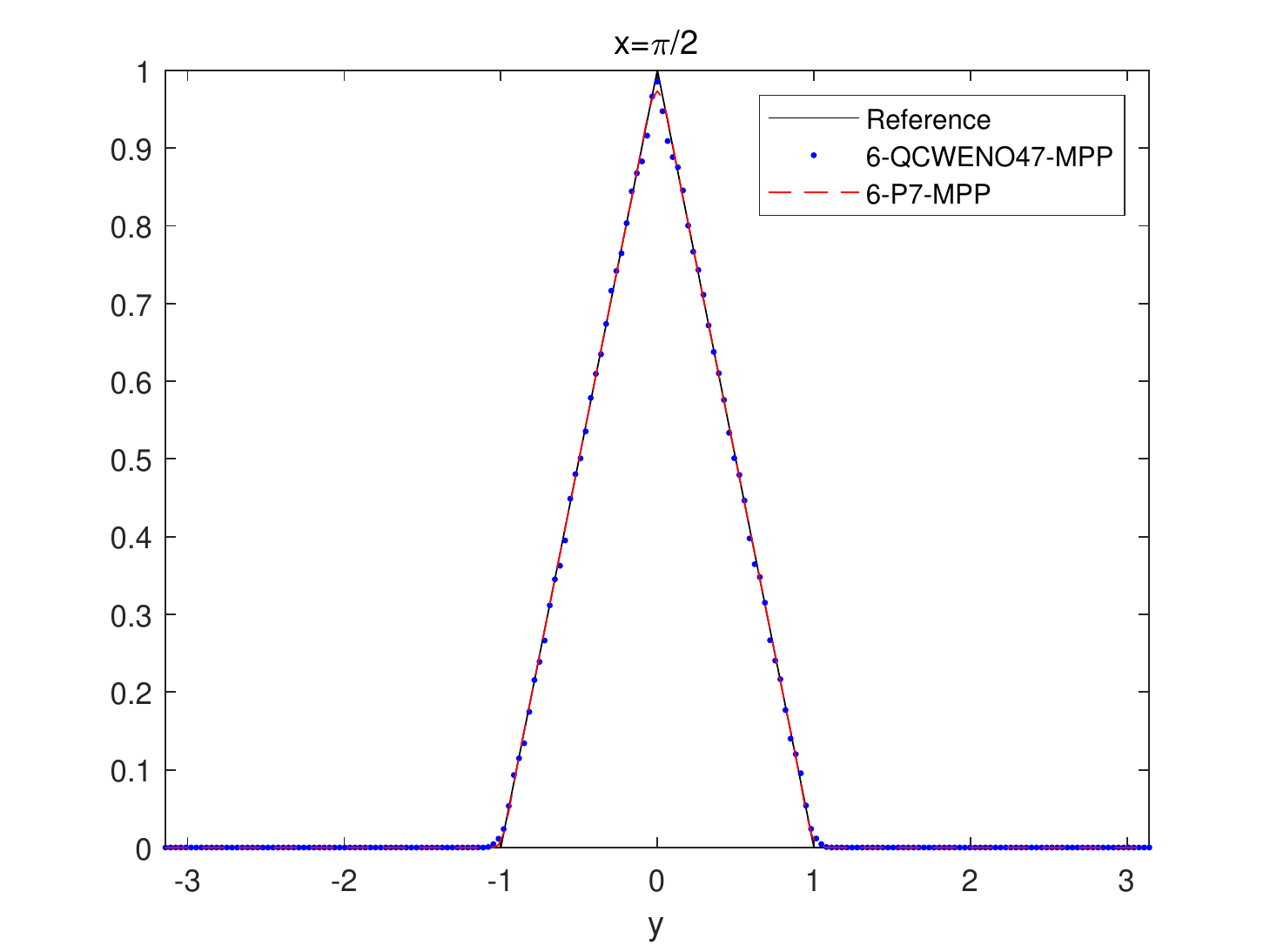}
			\subcaption{Cone}\label{rigid QCWENO47 Cone}
		\end{subfigure}	
		\begin{subfigure}[b]{0.45\linewidth}
			\includegraphics[width=1\linewidth]{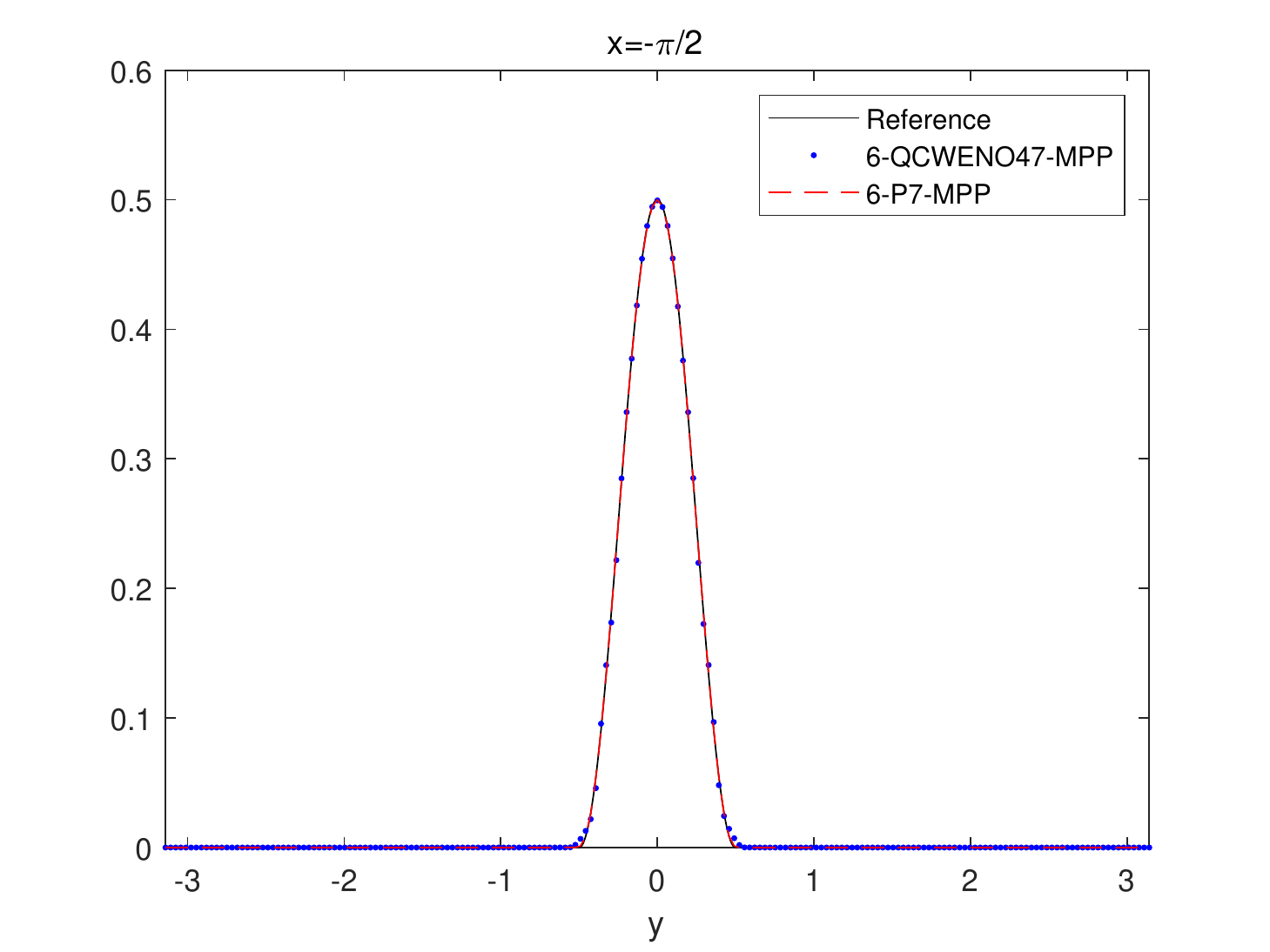}
			\subcaption{Humb}\label{rigid QCWENO47 hump}
		\end{subfigure}	
		\begin{subfigure}[b]{0.45\linewidth}
			\includegraphics[width=1\linewidth]{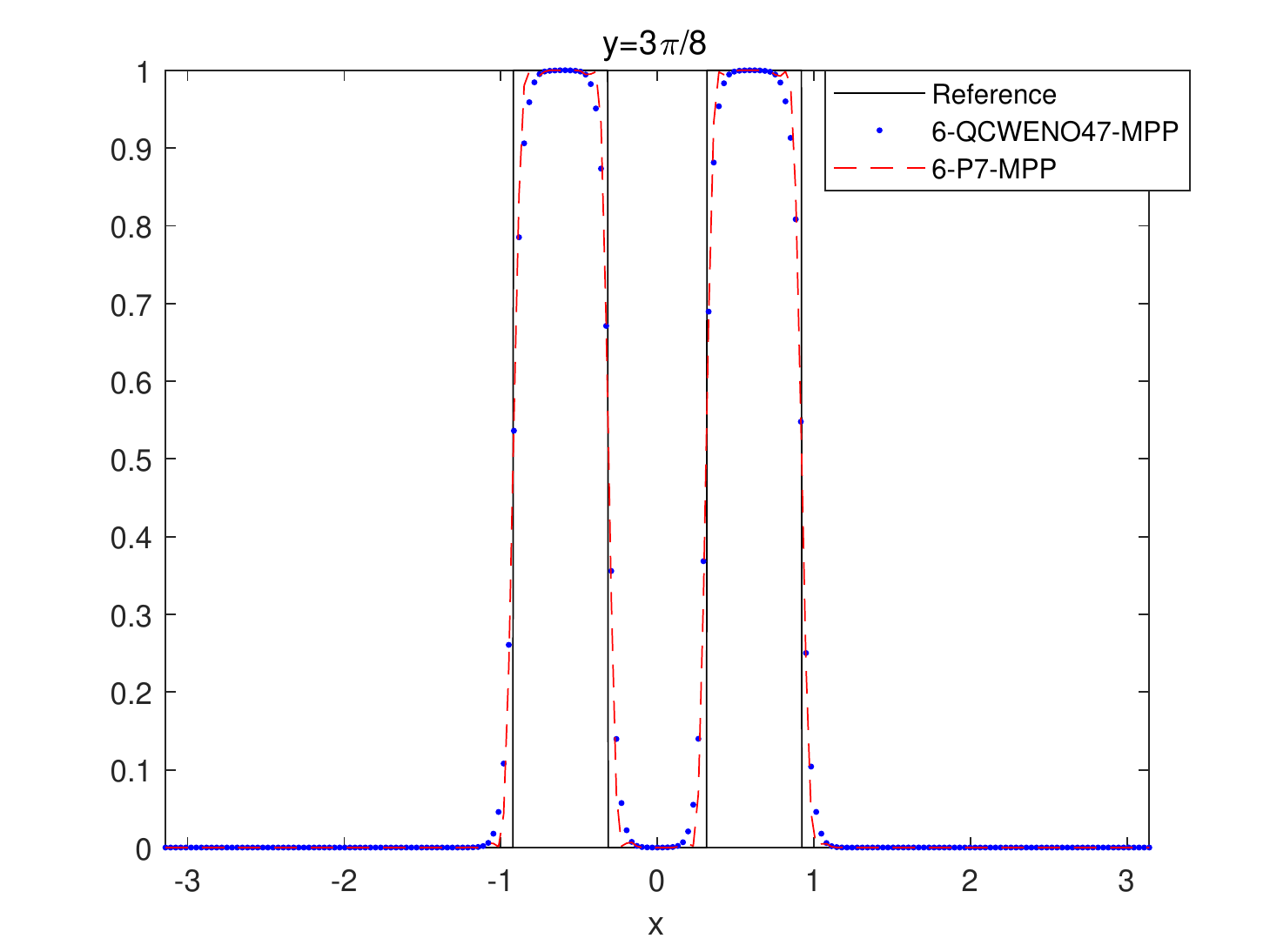}
			\subcaption{Cylinder}\label{rigid QCWENO47 Cylinder}
		\end{subfigure}	
		\begin{subfigure}[b]{0.45\linewidth}
			\includegraphics[width=1\linewidth]{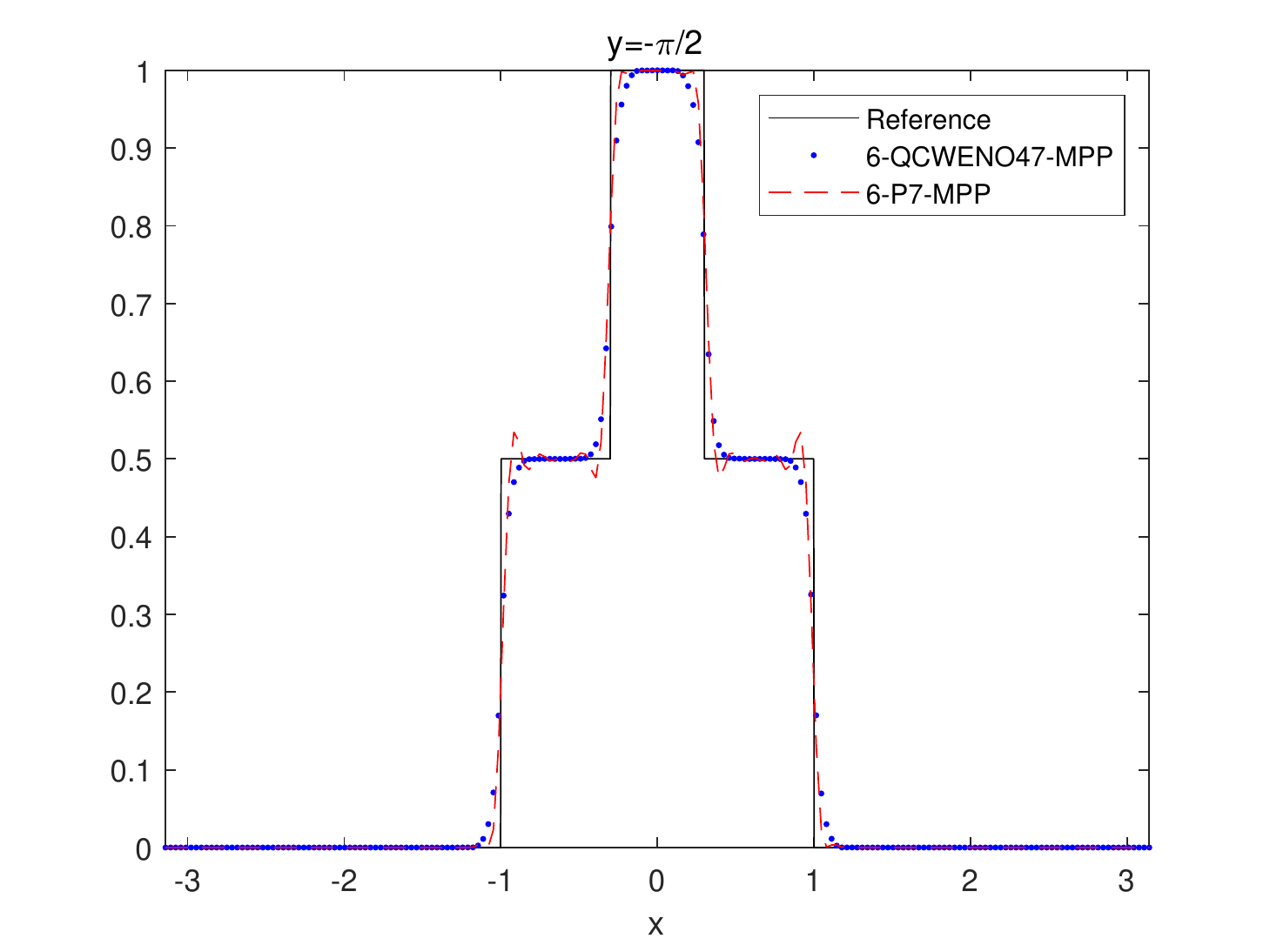}
			\subcaption{Double layer}\label{rigid QCWENO47 double layer}
		\end{subfigure}					
		\caption{2D Rigid body motion. Numerical solutions obtained for $(N_x,N_y)=192 \times 192$. For CWENO47, we use $\varepsilon=10^{-4}$}\label{rigid QCWENO47 all}
	\end{figure}


	\subsection{1D Vlasov-Poisson model}
	
	In this section we apply the high order splitting schemes described in Section \ref{SplittingMethods} to the Vlasov-Poisson system.
	
	To maintain conservation and MPP properties of analytical solutions for \eqref{VP}, we consider the numerical methods 
	in Table \ref{tab vp name}.
	\begin{center}
		\begin{table}[ht]
			\centering
			{\begin{tabular}{ccc}
					\hline
					\multicolumn{1}{ |c }{}&
					\multicolumn{1}{ |c| }{Order of time splitting method}& \multicolumn{1}{ c|  }{Basic reconstruction}   \\ \hline
					\multicolumn{1}{ |c }{2-P3-MPP}&
					\multicolumn{1}{ |c|  }{2nd order Strang splitting method \eqref{strang}}     
					&\multicolumn{1}{ |c|  }{P3}    
					\\
					\hline
					\multicolumn{1}{ |c }{4-P5-MPP}&
					\multicolumn{1}{ |c|  }{4th order Yoshida splitting \eqref{yoshida 4}}     
					&\multicolumn{1}{ |c|  }{P5}    
					\\
					\hline
					\multicolumn{1}{ |c }{6-P7-MPP}&
					\multicolumn{1}{ |c|  }{6th order method in \eqref{meren}}     
					&\multicolumn{1}{ |c|  }{P7}    
					\\
					\hline
			\end{tabular}}
			\caption{Numerical methods with MPP limiter used for 1D Vlasov-Poisson model.}\label{tab vp name}
		\end{table}
	\end{center}
	\vspace{-5mm}
	For 1D problems $(d_x, d_v)=(1,1)$, we adopted a  uniform mesh $\Delta x$ and $\Delta v$ both in space and velocity domain. The CFL number is defined as
	%
	
	\begin{align}\label{CFL VP 1d}
	\text{CFL}\displaystyle :=\max\left\{ \max_j|v_j|\frac{\Delta t}{\Delta x},\,\max_i|\mathbb{E}_i|\frac{\Delta t}{\Delta v} \right\}.
	\end{align}

	\subsubsection{Accuracy test}	
	
	%
	The aim of this section is to check the accuracy of high order splitting methods for the Vlasov-Poisson system. Here, we consider the following smooth initial data  \cite{filbet2001comparison}:
	\begin{align}\label{two stream instability}
	f(x,v,0)= \frac{2}{7\sqrt{2\pi}}(1+5v^2)\bigg(1 + \alpha\big( (\cos(2kx) + \cos(3kx))/1.2 + \cos(kx) \big)\bigg) \exp\left(-\frac{v^2}{2}\right), \quad \alpha=0.01,\quad k=0.5.
	\end{align}
	We assume periodic boundary condition on the interval $[0, 4\pi]$ and zero-boundary condition on velocity domain $[-8, 8]$. Using \eqref{CFL VP 1d}, We take a time step determined by \eqref{CFL VP 1d} with CFL$=4,\,80$. 
	
	From Table~\ref{tab vp name}, we note that each numerical scheme takes a spatial reconstruction whose order is higher than that of time splitting method. This means that for a fixed CFL number and sufficiently refined grid in all independent variables,  spatial errors become smaller than time errors, especially  for large CFL numbers and/or long integration time. 
	On the contrary, we can expect that convergence rates for reconstructions are easily confirmed for a small CFL numbers and/or  short integration time.
	
	For these reasons, we consider two cases for accuracy tests: (1) $N_x=N_v=40,80,160,320,640$ with CFL$=4$, (2) $N_x=N_v=160,320,640,1280,2560$ with CFL$=80$. Here we check the accuracy of all schemes based on the relative $L^1$-errors \eqref{rel norm} and convergence rates \eqref{rate} of the numerical solutions $f$ at a final time $T_f=2$.
	
	\begin{center}
		\begin{table}[ht]
			\centering
			{\begin{tabular}{|ccccccc|}
					\hline
					\multicolumn{1}{ |c }{}&
					\multicolumn{1}{ |c| }{}&  \multicolumn{2}{ |c }{CFL$=4$}& \multicolumn{1}{ |c| }{} &
					\multicolumn{2}{ |c| }{CFL$=80$} \\ \hline
					\multicolumn{1}{ |c }{}&
					\multicolumn{1}{ |c|  }{$(N_x,2N_x)$} &
					\multicolumn{1}{ |c  }{error} &
					\multicolumn{1}{ c  }{rate} &
					\multicolumn{1}{ |c|  }{$(N_x,2N_x)$} &
					\multicolumn{1}{ |c  }{error} &
					\multicolumn{1}{ c|  }{rate}     \\ 
					\hline
					\hline	
					\multicolumn{1}{ |c }{}&
					\multicolumn{1}{ |c|  }{$(40,80)$}                    
					&9.6577e-04
					&2.8269
					&\multicolumn{1}{ |c|  }{$(160,320)$}   
					&1.5142e-03
					&2.4288					
					\\
					\multicolumn{1}{ |c }{2-P3-MPP}&
					\multicolumn{1}{ |c|  }{$(80,160)$}              
					&1.3611e-04
					&2.9325
					&\multicolumn{1}{ |c|  }{$(320,640)$}                       
					&2.8122e-04
					&2.0837
					\\
					\multicolumn{1}{ |c }{}&
					\multicolumn{1}{ |c|  }{$(160,320)$}             
					&1.7829e-05&
					2.9189
					&\multicolumn{1}{ |c|  }{$(640,1280)$}                    &6.6341e-05
					&2.0194
					\\
					\multicolumn{1}{ |c }{}&
					\multicolumn{1}{ |c|  }{$(320,640)$}      
					&
					2.3575e-06      &   &\multicolumn{1}{ |c|  }{$(1280,2560)$} 
					&1.6364e-05&
					\\
					\hline
					\multicolumn{1}{ |c }{}&
					\multicolumn{1}{ |c|  }{$(40,80)$}          	    
					&1.4110e-02               
					&0.3910                      
					&\multicolumn{1}{ |c|  }{$(160,320)$}                   
					&1.9771e-03
					&2.7746
					\\
					\multicolumn{1}{ |c }{4-P5-MPP}&
					\multicolumn{1}{ |c|  }{$(80,160)$}                   
					&1.0760e-02               
					&13.0707
					&\multicolumn{1}{ |c|  }{$(320,640)$}                     &2.8893e-04
					&3.6838
					\\
					\multicolumn{1}{ |c }{}&
					\multicolumn{1}{ |c|  }{$(160,320)$}        
					&1.2507e-06
					&4.9704                      &\multicolumn{1}{ |c|  }{$(640,1280)$}                         &2.2483e-05
					&3.9038
					\\
					\multicolumn{1}{ |c }{}&
					\multicolumn{1}{ |c|  }{$(320,640)$}      
					&3.9895e-08               &   &\multicolumn{1}{ |c|  }{$(1280,2560)$} &1.5020e-06
					&
					\\
					\hline	
					\multicolumn{1}{ |c }{}&
					\multicolumn{1}{ |c|  }{$(40,80)$}
					&6.3266e-05         
					& 6.5135                                 
					&\multicolumn{1}{ |c|  }{$(160,320)$}        
					&7.0569e-07
					&8.7136
					\\
					\multicolumn{1}{ |c }{6-P7-MPP}&
					\multicolumn{1}{ |c|  }{$(80,160)$}
					&6.9250e-07         
					&6.8640
					&\multicolumn{1}{ |c|  }{$(320,640)$}                       
					&1.6810e-09					
					&5.9542
					\\
					\multicolumn{1}{ |c }{}&
					\multicolumn{1}{ |c|  }{$(160,320)$}
					&5.9451e-09
					&6.9704
					&\multicolumn{1}{ |c|  }{$(640,1280)$}     
					&2.7113e-11
					&5.9407
					\\
					\multicolumn{1}{ |c }{}&
					\multicolumn{1}{ |c|  }{$(320,640)$}      
					&4.7408e-11         &   &\multicolumn{1}{ |c|  }{$(1280,2560)$} &4.4141e-13
					&
					\\
					\hline
					\hline
			\end{tabular}}
			\caption{Accuracy test for the 1D Vlasov-Poisson system. Initial data is given in \eqref{two stream instability}.}\label{tab3}
		\end{table}
	\end{center}
	
	The final time is relatively short: a typical particle will move only a fraction of the domain before time $T_f$.   
	In Table \ref{tab3}, we report numerical results for all schemes. As expectated, for CFL$=80$ we obtain the convergence rates of time splitting methods. On the contrary, for relatively small CFL$=4$ we observe the accuracy of spatial reconstructions. In both cases, we note that the use of MPP limiter does not lead to any order reduction.
	
	\subsubsection{1D Long time simulation}
	Next, we move on to long time simulations for VP system. 
	
	The Vlasov-Poisson system admits several time invariants, such as $L^p$-norms of the solution total energy and entropy. These quantities can be used as  diagnostics for the numerical methods. Here we monitor the following 
	four quantities that should remain constant in time:
	
	\begin{itemize}
		\item $L^p$-norm of $f^n$, $p=1,2$
		\[
		\|f^n\|_p = \left(\sum_{i,j}|f_{i,j}^n|^p(\Delta x)^{d_x}(\Delta v)^{d_v}\right)^{\frac{1}{p}}.
		\]
		\item Total energy of $f^n$
		\[
		\text{Energy} = \sum_{i,j}\frac{|v_j|^2}{2}f_{i,j}^{n}(\Delta x)^{d_x}(\Delta v)^{d_v} + \sum_{i}\frac{|\mathbb{E}_i^{n}|^2}{2} (\Delta x)^{d_x}.
		\]
		\item Entropy of $f^n$
		\[
		\text{Entropy} = -\sum_{i,j}f_{i,j}^{n}\log{f_{i,j}^n}(\Delta x)^{d_x}(\Delta v)^{d_v}.
		\]
	\end{itemize}
	
	In addition, we monitor the behaviour of the electric field:
	\begin{itemize} 
		\item $L^2$-norm of $\mathbb{E}^n$
		\[
		\|\mathbb{E}^n\|_2 = \left(\sum_{i}|\mathbb{E}_i^n|^2(\Delta x)^{d_x}\right)^{\frac{1}{2}}.
		\]
		\item $L^\infty$-norm of $\mathbb{E}^n$
		\[
		\|\mathbb{E}^n\|_\infty = \max_i |\mathbb{E}_i^n|.
		\]
	\end{itemize}

	We consider two benchmark problems:
	\begin{enumerate}
		\item 	\textbf{Weak Landau damping}:
		We consider initial data in \cite{rossmanith2011positivity}: 
		\begin{align}\label{weak landau sec}
		f(x,v,0)= \frac{1}{\sqrt{2\pi}}\bigg(1+ \alpha \cos(k x)\bigg) \exp\left(-\frac{v^2}{2}\right), \quad \alpha=0.01,\quad k=0.5.
		\end{align}
		
		\item \textbf{Strong Landau damping}: We consider initial data in \cite{rossmanith2011positivity}
		\begin{align}\label{Strong Landau damping}
		f(x,v,0)= \frac{1}{\sqrt{2\pi}}\bigg(1+ \alpha \cos(k x)\bigg) \exp\left(-\frac{v^2}{2}\right), \quad \alpha=0.5,\quad k=0.5.
		\end{align}
		%
		
		%
	\end{enumerate}	
	
	For both problems, we impose periodic boundary condition on the physical domain $[-2\pi, 2\pi]$ and zero-boundary condition on velocity domain $[-2\pi, 2\pi]$. Numerical solutions are computed up to final time $T_f=60$ with a time step determined by CFL$=2$ using \eqref{CFL VP 1d}. The grid size is  $N_x=64$ and $N_v=128$.
	
	In Figure \ref{weak landau}, we plot numerical results obtained by high-order splitting methods listed in Table \ref{tab vp name}. The time evolution of $L^2$ and $L^{\infty}$-norm of electric field $E^n$ shows a theoretical damping rate for all schemes. After $t=50$, only scheme 6-P7-MPP maintains the correct damping rate. 
	
	Regarding the preservation of conserved quantities such as the $L^1$ and $L^2$-norm/energy/entropy of numerical solutions, all schemes produce small conservation errors within $10^{-6}$. Among them, the 6-P7-MPP scheme gives the most accurate profiles.
	These results are comparable with those obtained with other conservative finite difference methods available in the literature, such as, for example, \cite{RQT}.
	
	\begin{figure}[htbp]
		\centering
		\begin{subfigure}[b]{0.45\linewidth}
			\includegraphics[width=1\linewidth]{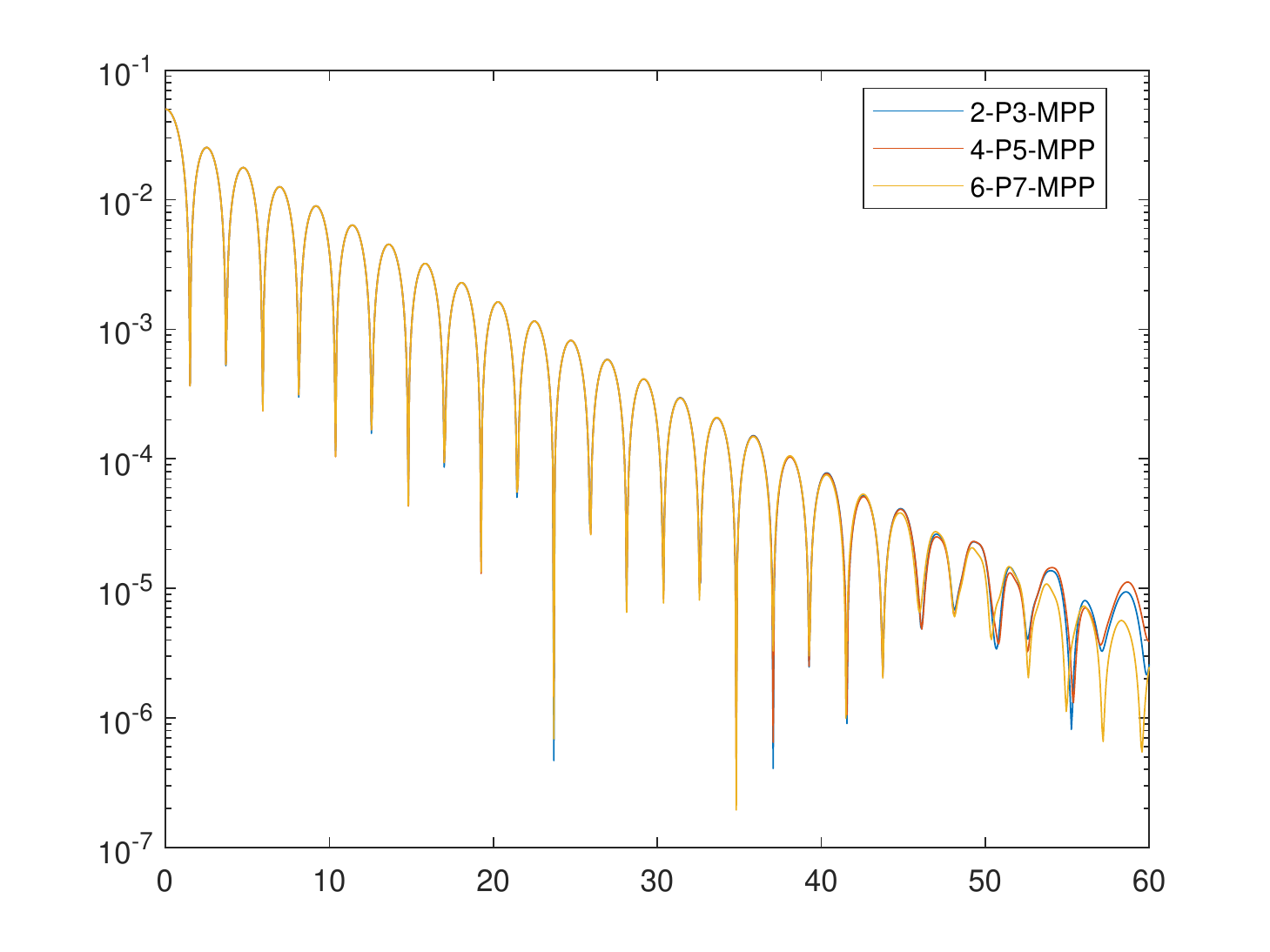}
			\subcaption{$L^2$-norm of $E$} \label{weak landau E}
		\end{subfigure}	
		\vspace*{5mm}
		\begin{subfigure}[b]{0.45\linewidth}
			\includegraphics[width=1\linewidth]{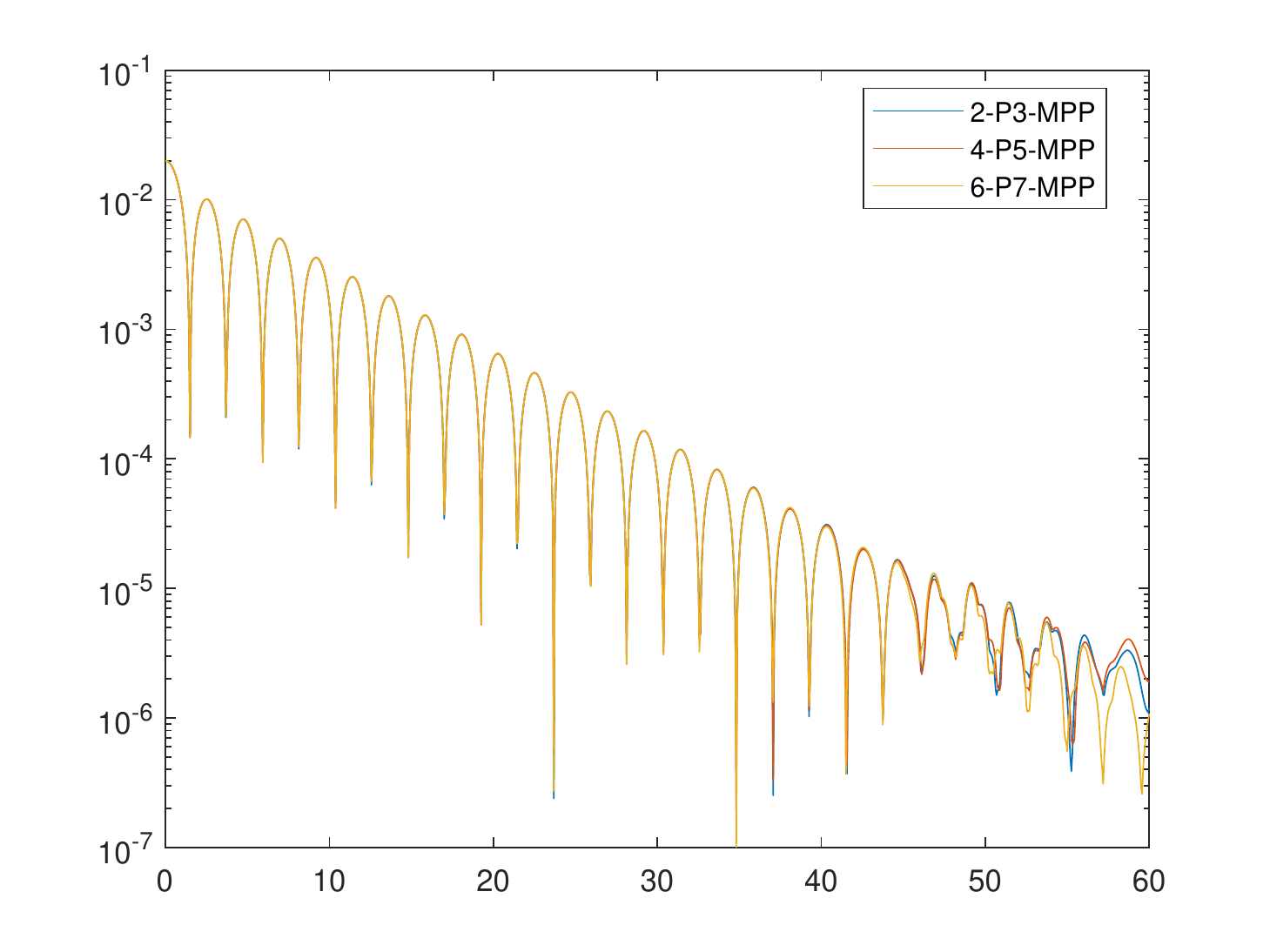}
			\subcaption{$L^{\infty}$-norm of $E$}
		\end{subfigure}	
		\vspace{5mm}
		%
		\begin{subfigure}[b]{0.45\linewidth}
			\includegraphics[width=1\linewidth]{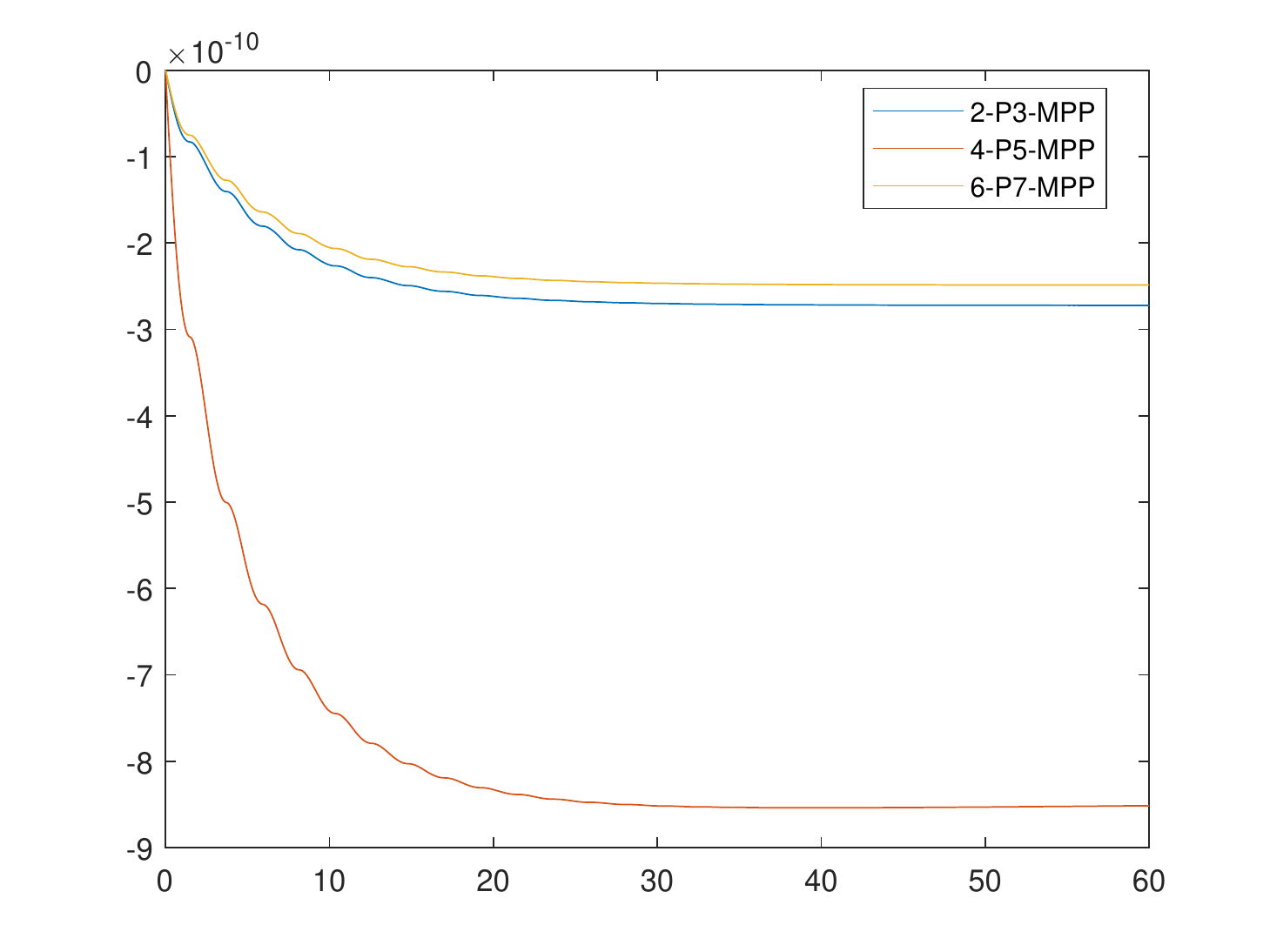}
			\subcaption{$\|f^{N_t}\|_1-\|f^0\|_1$}
		\end{subfigure}	
		\begin{subfigure}[b]{0.45\linewidth}
			\includegraphics[width=1\linewidth]{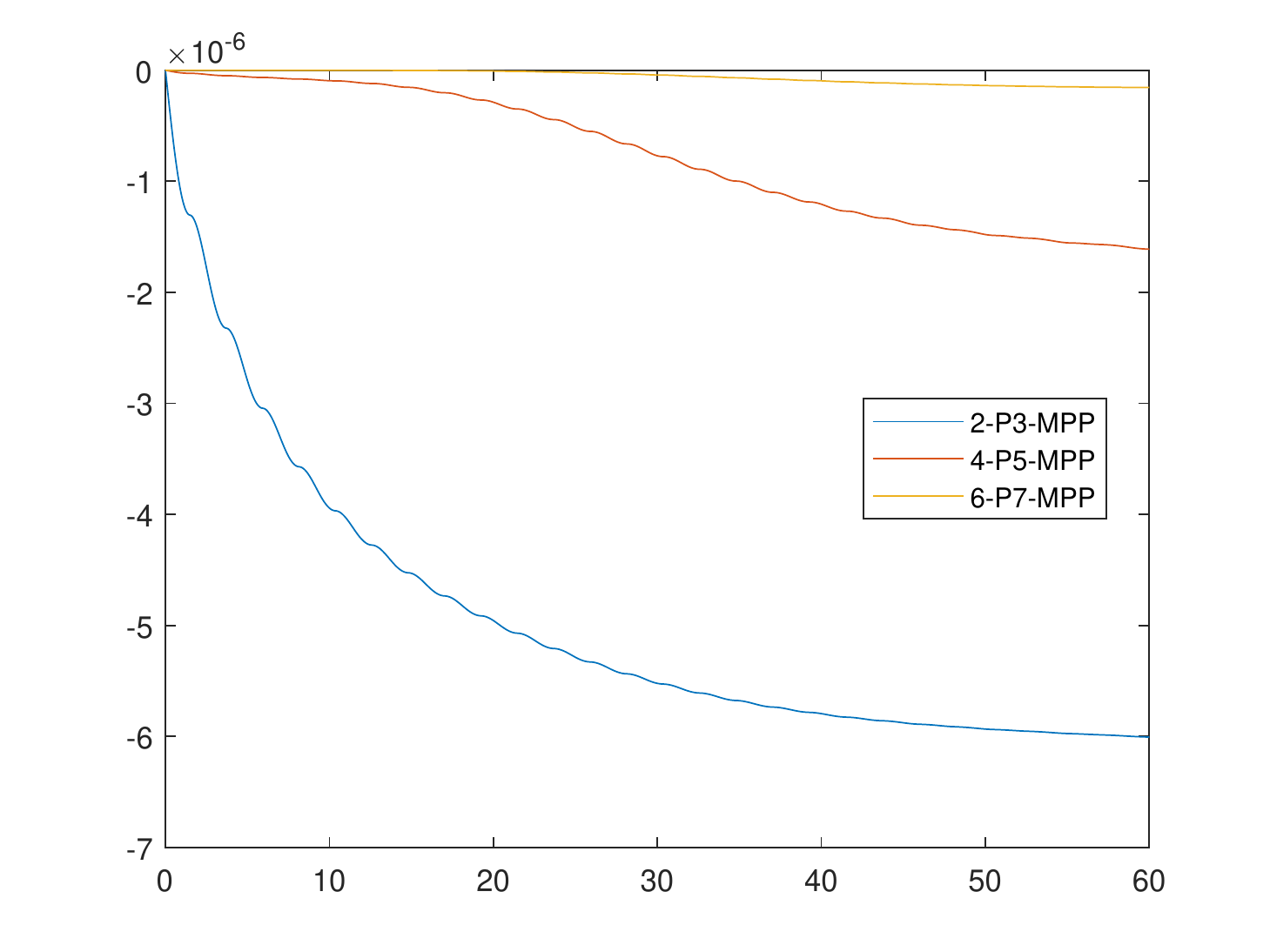}
			\subcaption{$\|f^{N_t}\|_1-\|f^0\|_2$}
		\end{subfigure}	
		\vspace{2mm}
		\begin{subfigure}[b]{0.45\linewidth}
			\includegraphics[width=1\linewidth]{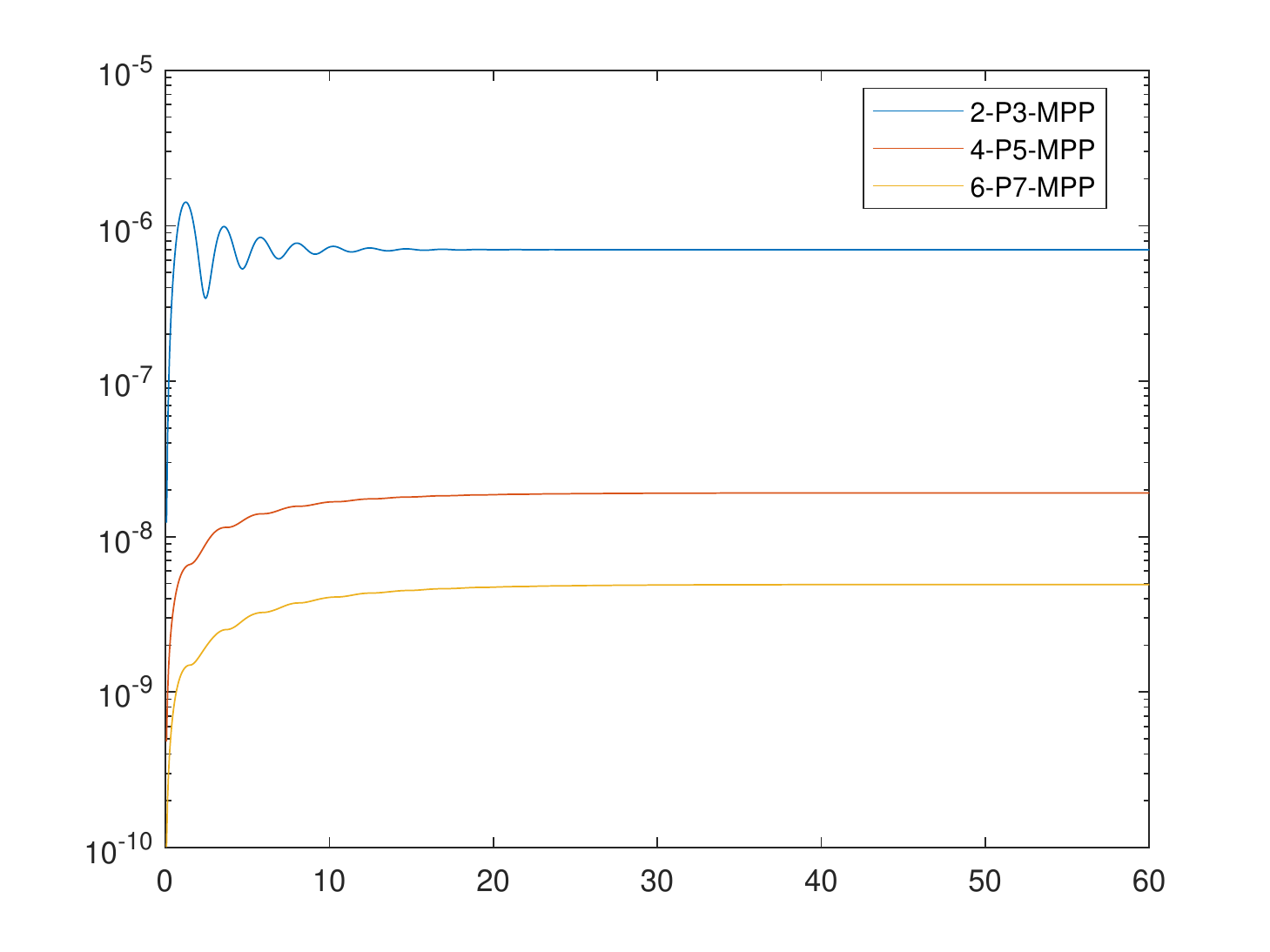}
			\subcaption{$\big|Energy(t=T_f)-Energy(t=0)\big|$}
		\end{subfigure}	
		\begin{subfigure}[b]{0.45\linewidth}
			\includegraphics[width=1\linewidth]{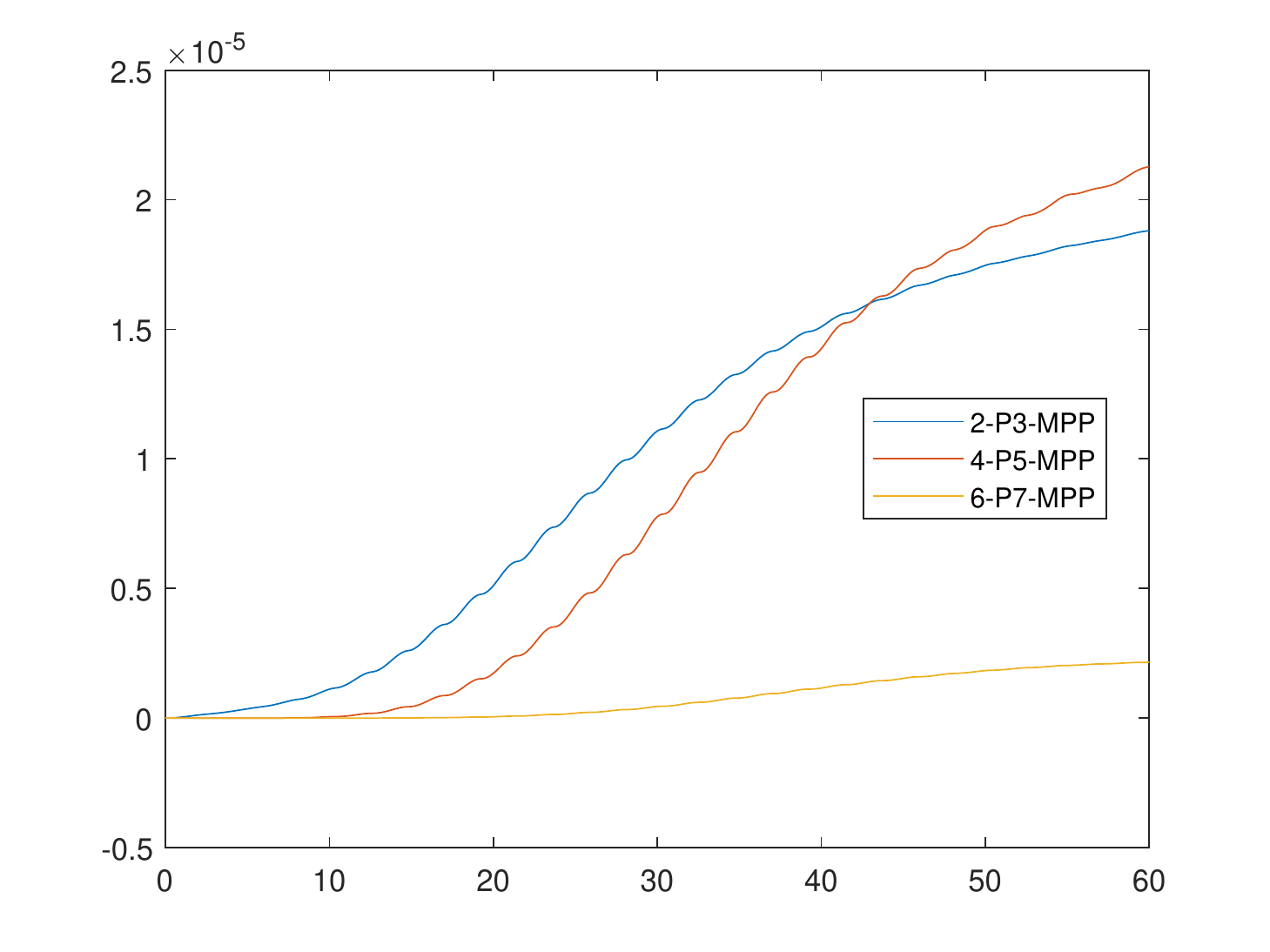}
			\subcaption{$Entropy(t=T_f)-Entropy(t=0)$}
		\end{subfigure}	
		
		\caption{1D Vlasov Poisson. Weak Landau damping with initial data \eqref{weak landau sec}. $N_x=64, N_v=128$}\label{weak landau}
	\end{figure}

	For strong Landau damping, Eq.~\eqref{Strong Landau damping}, we first compute numerical solutions taking the same domain and meshes used for the weak Landau damping. In Fig. \ref{strong landau}, numerical results are plotted for $N_x=64$, $N_v=128$ with velocity domain $[-2\pi,2\pi]$. Due to the increased perturbation amplitude $\alpha=0.5$, non-linear damping begins to appear around $t=10$. Such non-linear tendency of the electric field is well-captured by each scheme, but all conservation errors are non negligible. The small error in the $L^1$-norm is due to the loss of conservation caused by zero-boundary condition at the edge of the computational domain in velocity.
	For comparison with the results in \cite{rossmanith2011positivity}, we plot the distribution functions in phase space at different times in Fig. \ref{fig strong phase space}.
	
	{The effect of grid refinement is shown in  
		Fig.~\ref{strong landau 1}, where the number of grid points is double both in physical and velocity space. If we slightly enlarge the domain in velocity space, say from $-2.5\pi$ to $2.5\pi$, then conservation of $L^1$ norm improves by four orders of magnitude. }
	Although all schemes exhibit finite conservation errors in $L^2$-norm, entropy and energy, we confirm that scheme 6-P7-MPP has less conservation errors than schemes 2-P3-MPP and 4-P5-MPP.

	\begin{figure}[htbp]
		\centering
		\begin{subfigure}[b]{0.45\linewidth}
			\includegraphics[width=1\linewidth]{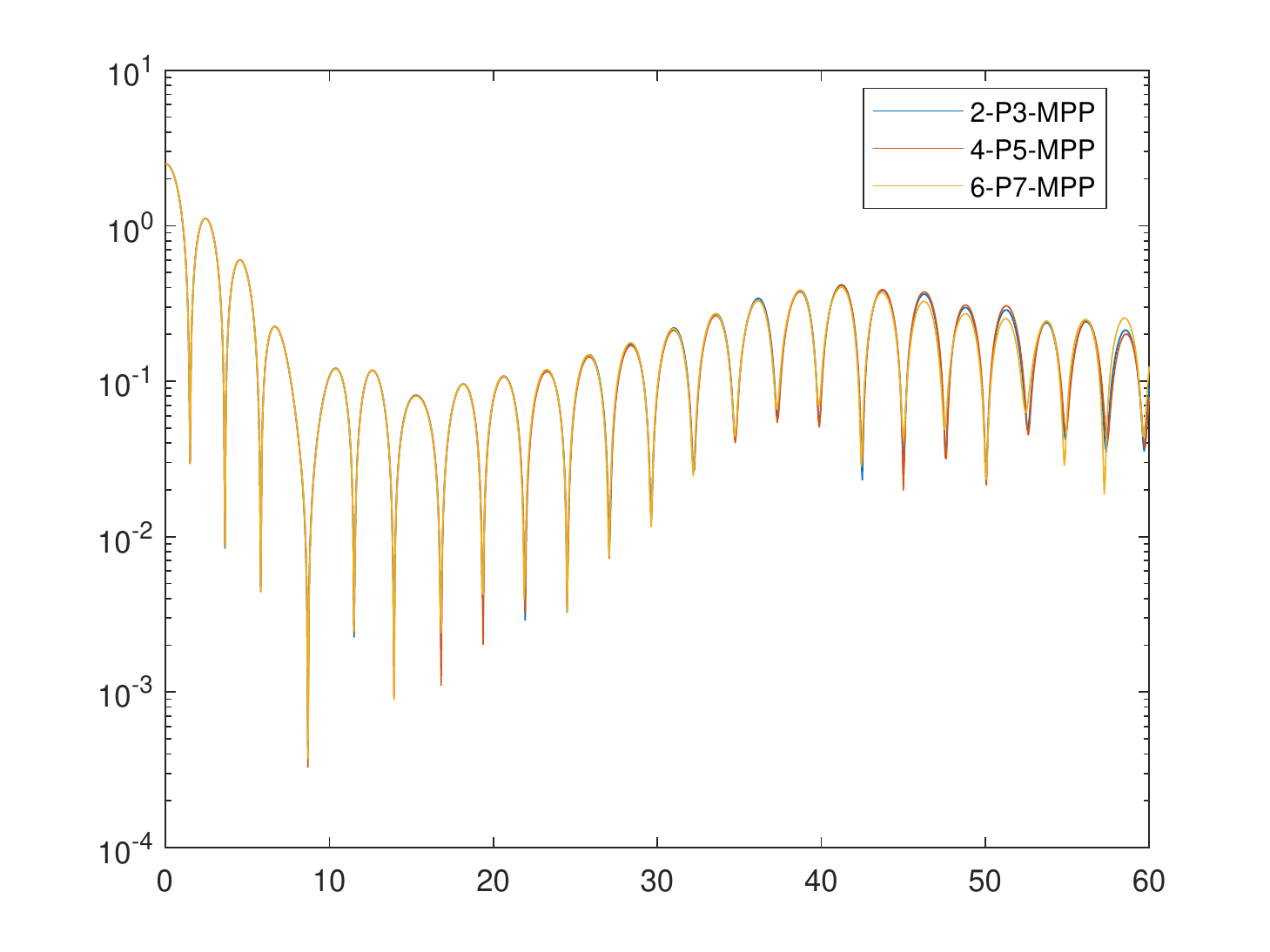}
			\subcaption{$L^2$-norm of $E$}
		\end{subfigure}	
		\vspace*{5mm}
		\begin{subfigure}[b]{0.45\linewidth}
			\includegraphics[width=1\linewidth]{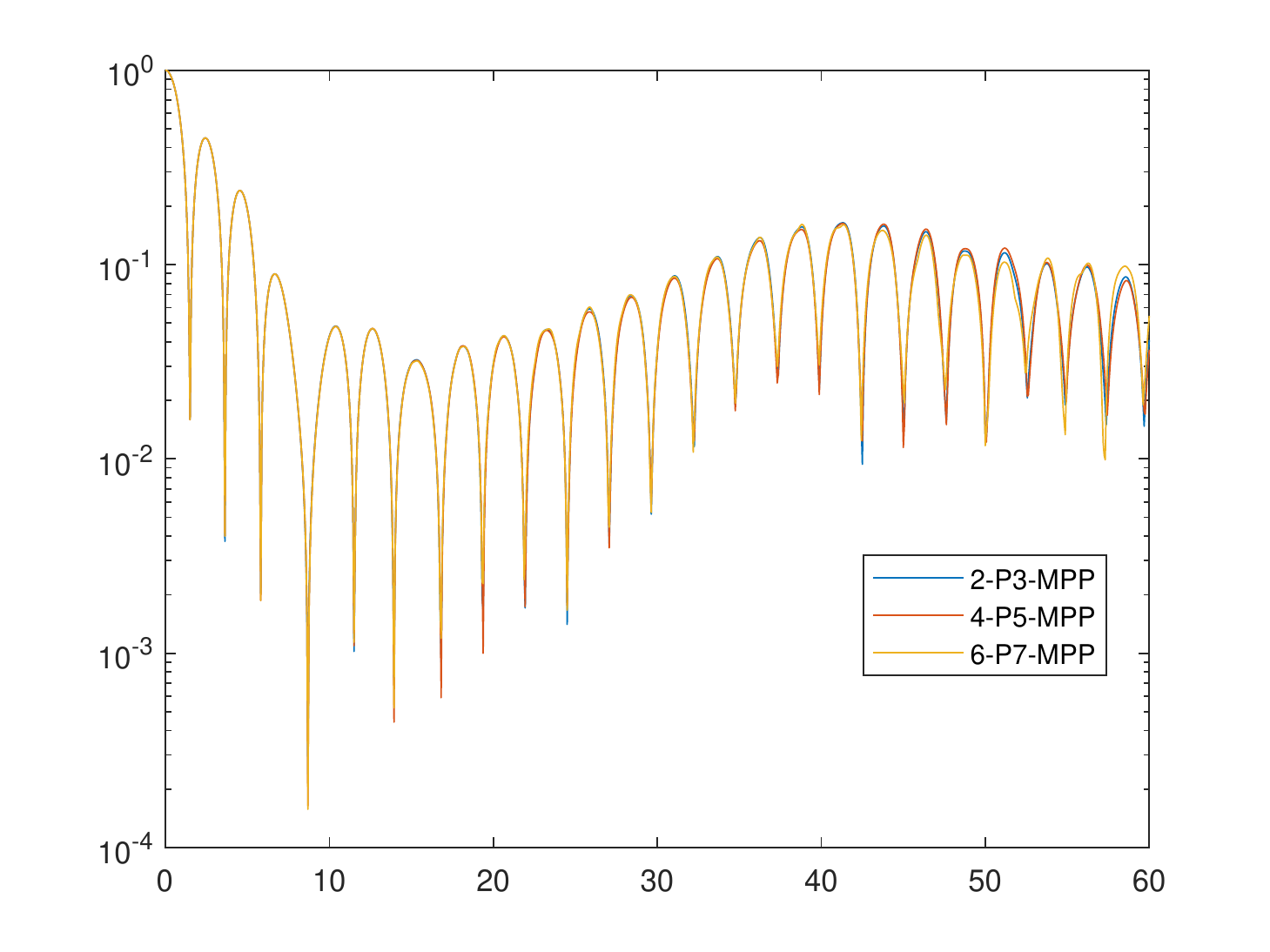}
			\subcaption{$L^{\infty}$-norm of $E$}
		\end{subfigure}	
		\vspace*{5mm}
		\begin{subfigure}[b]{0.45\linewidth}
			\includegraphics[width=1\linewidth]{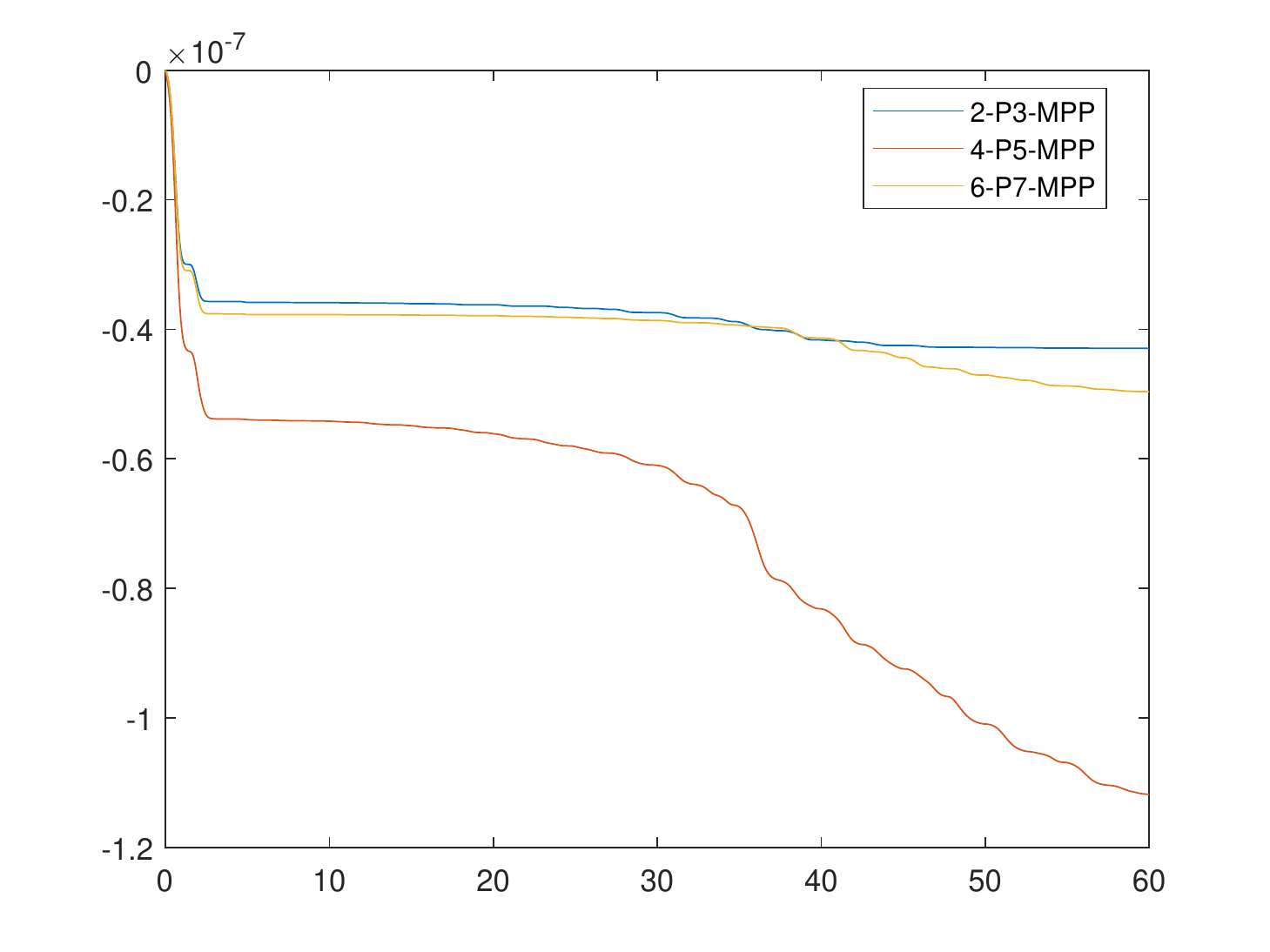}
			\subcaption{$\|f^{N_t}\|_1-\|f^0\|_1$}
			
		\end{subfigure}	
		\begin{subfigure}[b]{0.45\linewidth}
			\includegraphics[width=1\linewidth]{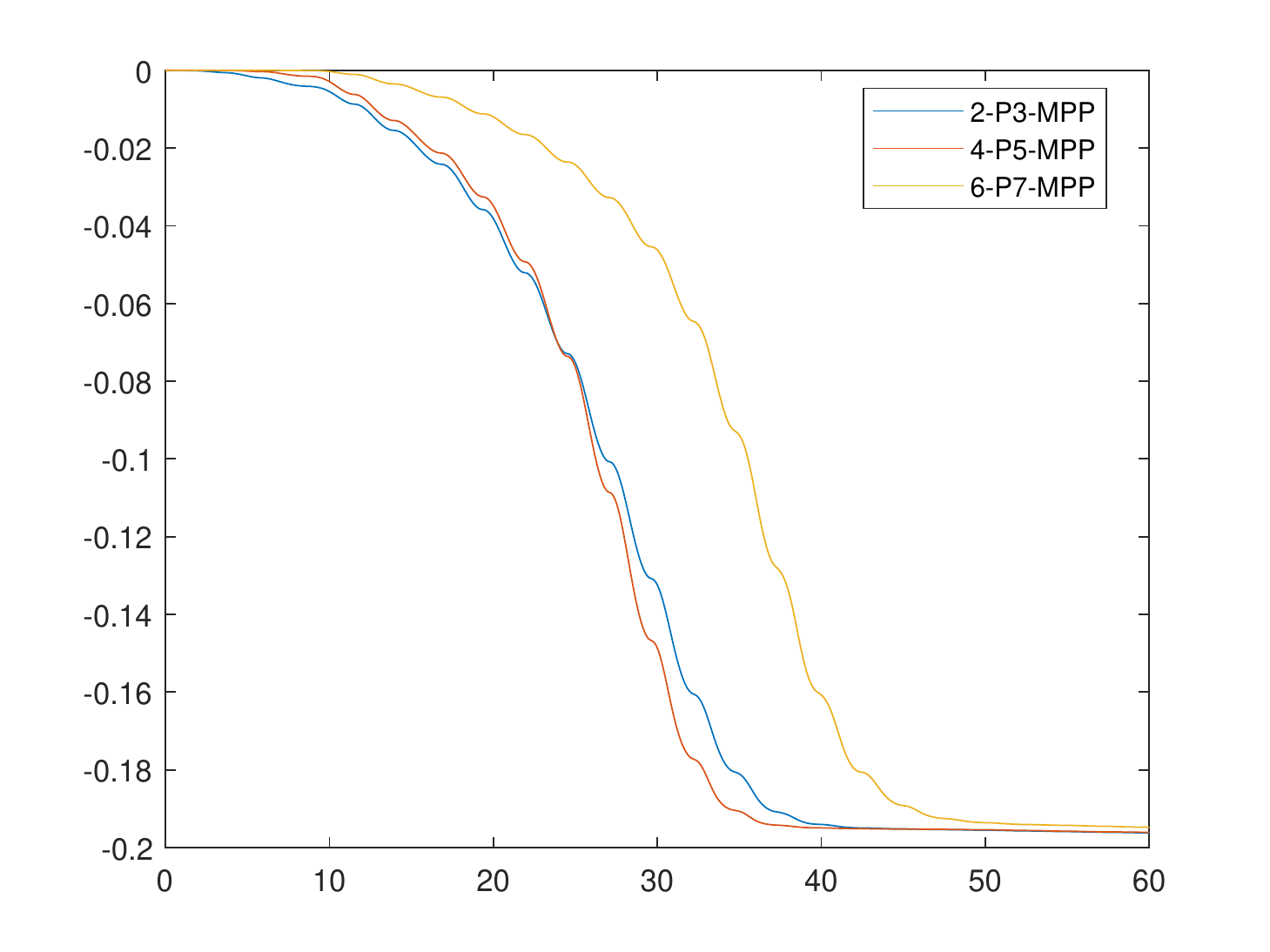}
			\subcaption{$\|f^{N_t}\|_2-\|f^0\|_2$}
			
		\end{subfigure}	
		\begin{subfigure}[b]{0.45\linewidth}
			\includegraphics[width=1\linewidth]{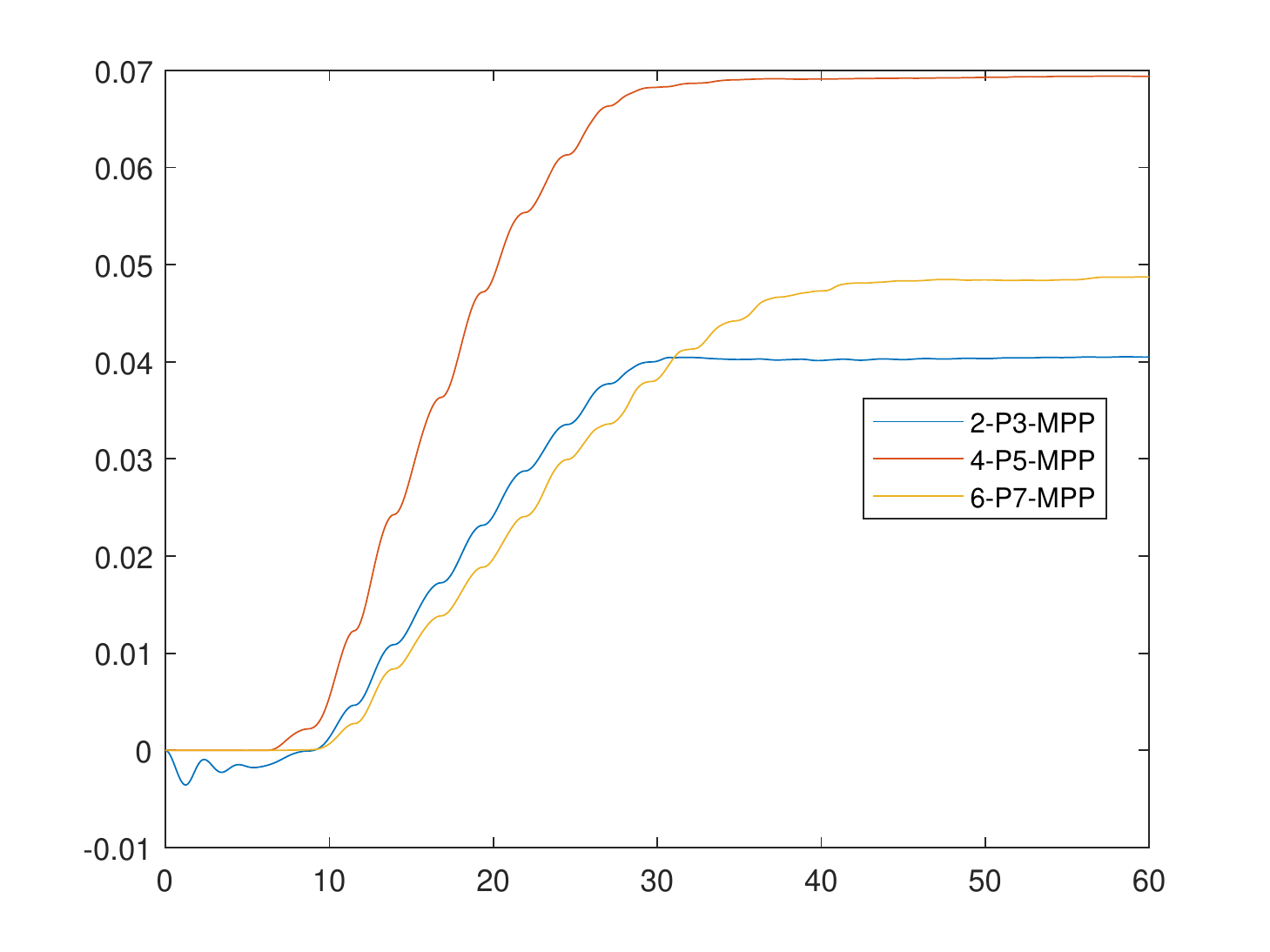}
			\subcaption{$Energy(t=T_f)-Energy(t=0)$}
		\end{subfigure}	
		\begin{subfigure}[b]{0.45\linewidth}
			\includegraphics[width=1\linewidth]{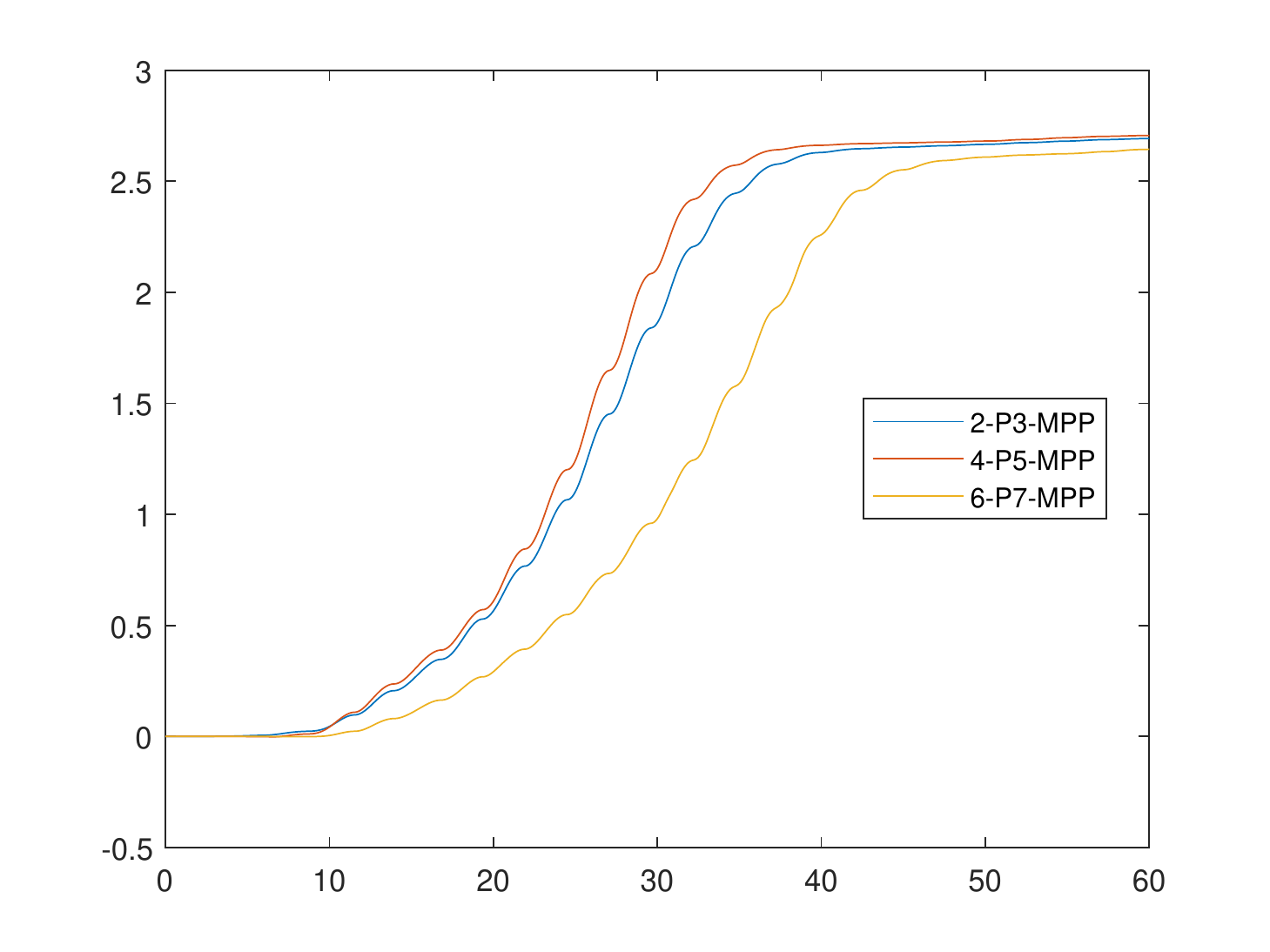}
			\subcaption{$Entropy(t=T_f)-Entropy(t=0)$}
		\end{subfigure}	
		
		\caption{1D Vlasov Poisson. Strong Landau damping with initial data \eqref{Strong Landau damping}. $N_x=64,\,N_v=128$}\label{strong landau}
	\end{figure}

	\begin{figure}[htbp]
		\centering
		\begin{subfigure}[b]{0.45\linewidth}
			\includegraphics[width=1\linewidth]{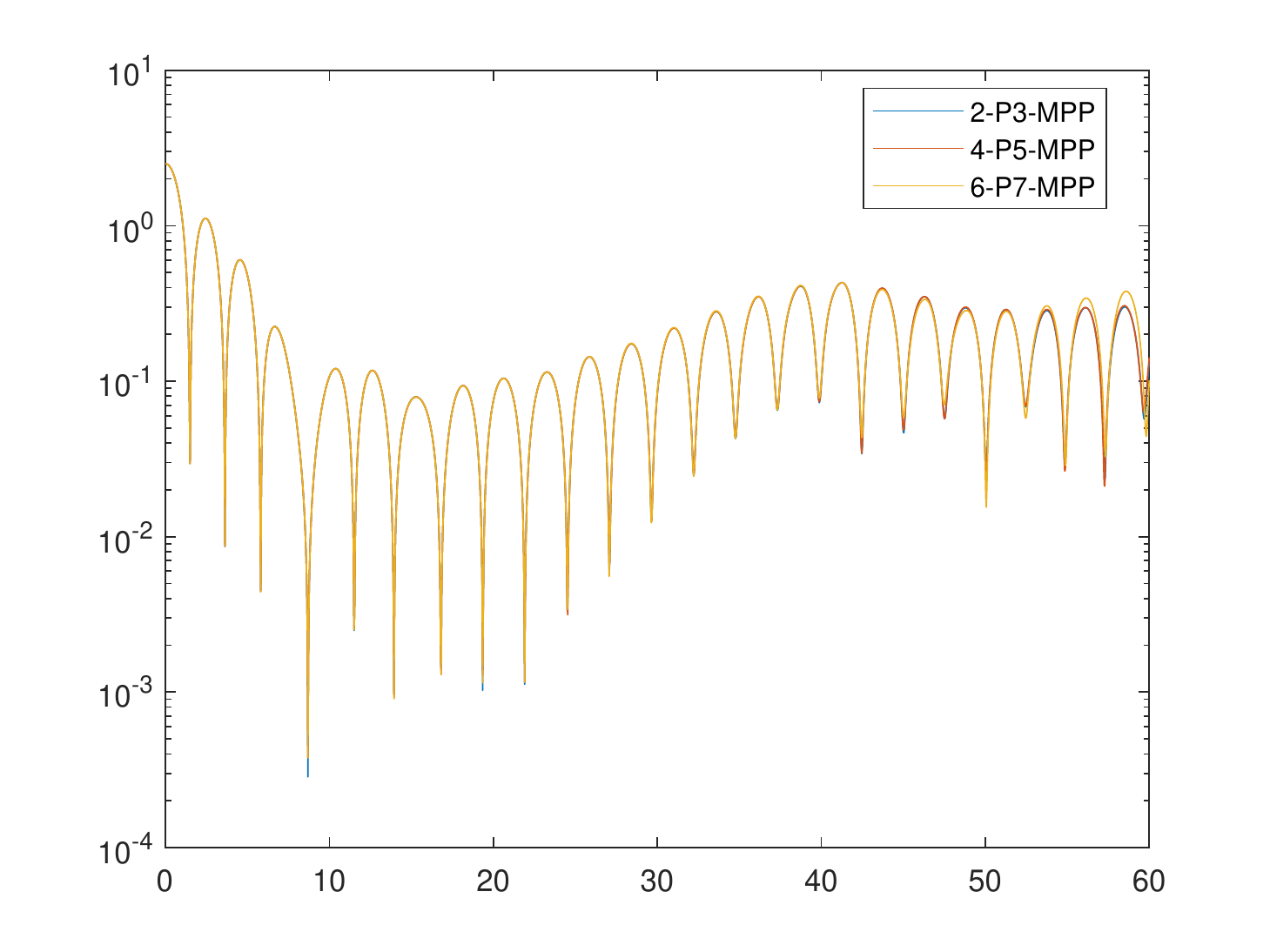}
			\subcaption{$L^2$-norm of $E$}
		\end{subfigure}	
		\vspace*{5mm}
		\begin{subfigure}[b]{0.45\linewidth}
			\includegraphics[width=1\linewidth]{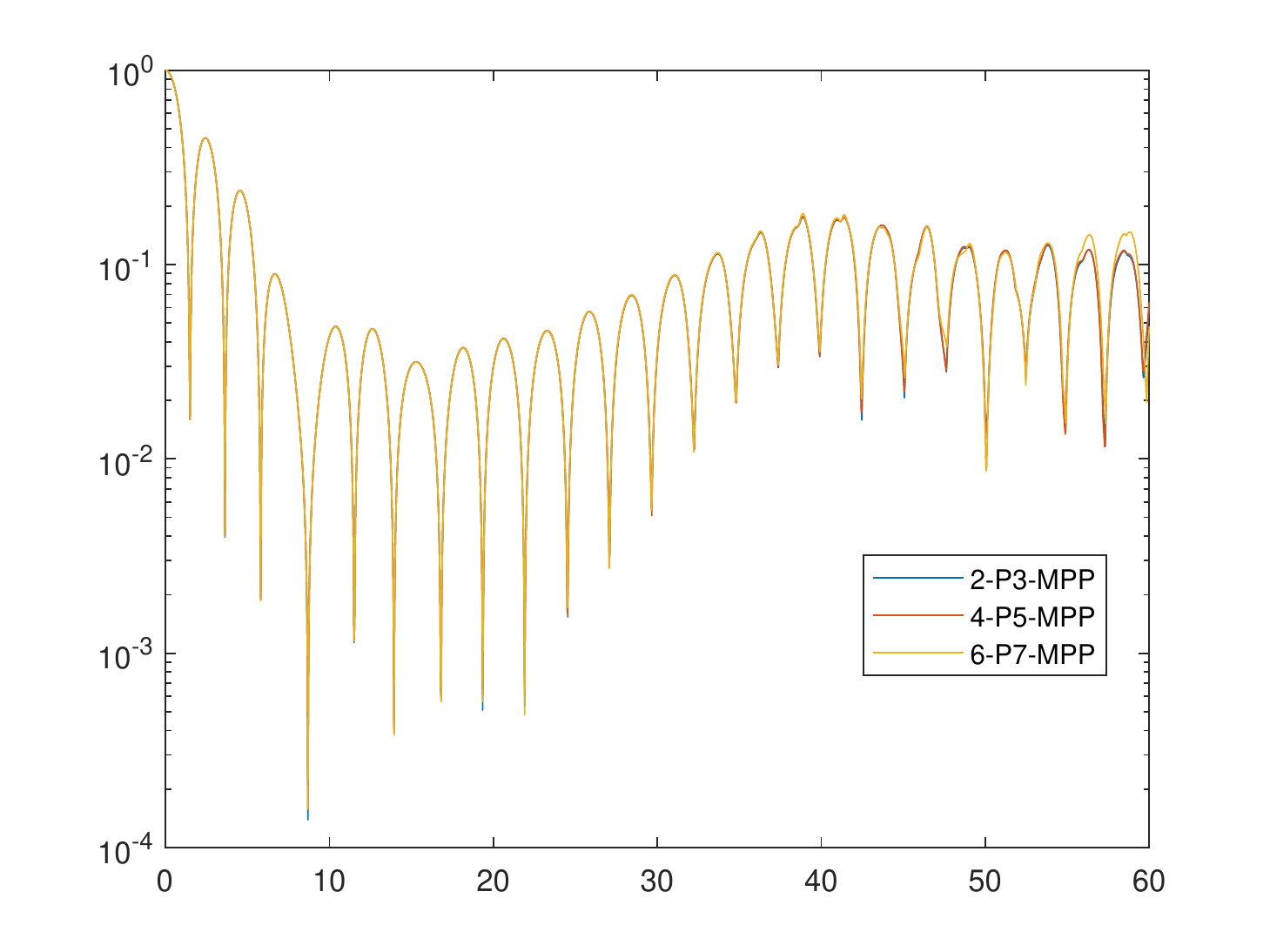}
			\subcaption{$L^{\infty}$-norm of $E$}
		\end{subfigure}	
		\vspace*{5mm}
		\begin{subfigure}[b]{0.45\linewidth}
			\includegraphics[width=1\linewidth]{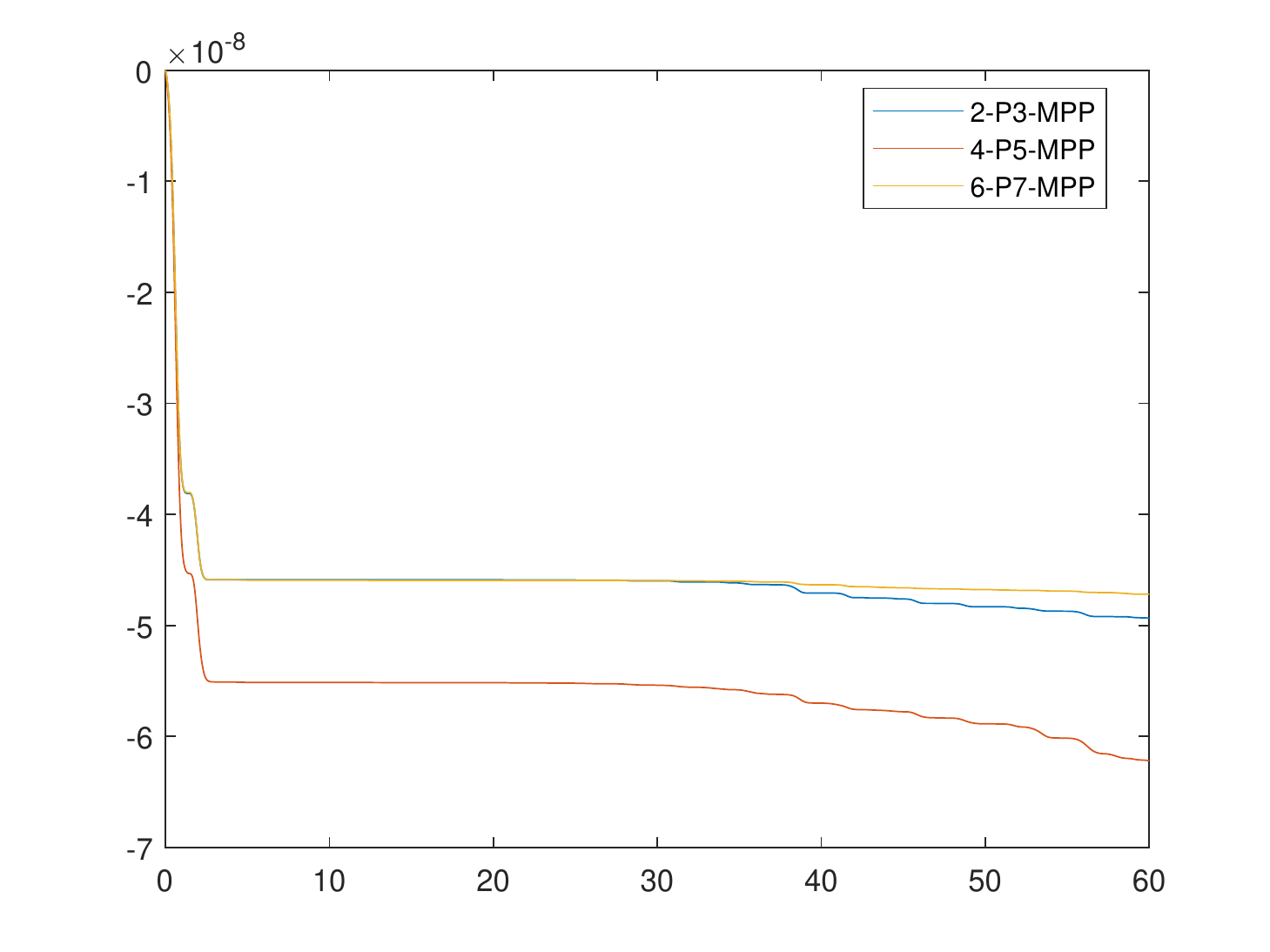}
			\subcaption{$\|f^{N_t}\|_1-\|f^0\|_1$}
			
		\end{subfigure}	
		\begin{subfigure}[b]{0.45\linewidth}
			\includegraphics[width=1\linewidth]{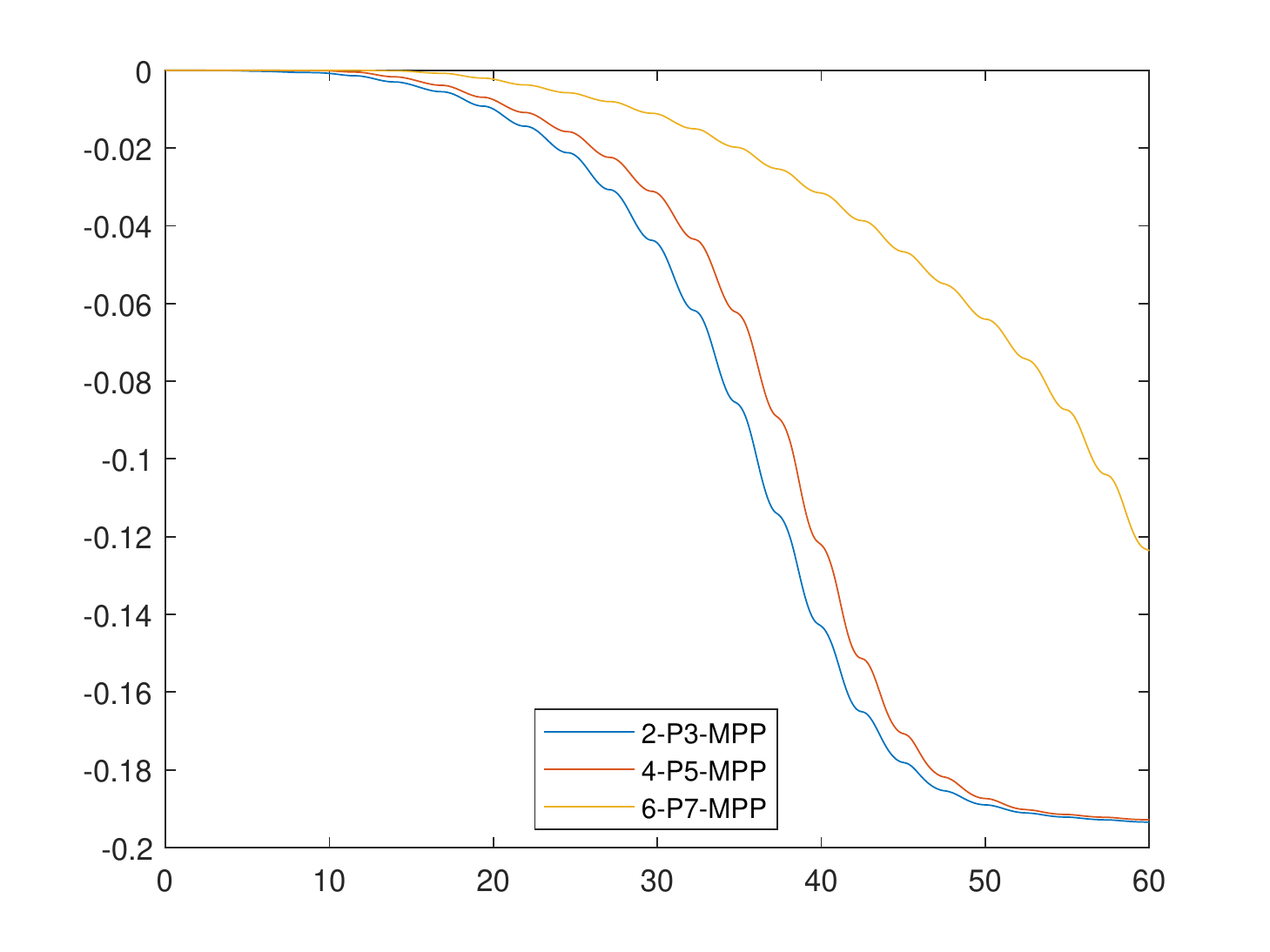}
			\subcaption{$\|f^{N_t}\|_2-\|f^0\|_2$}
			
		\end{subfigure}	
		\begin{subfigure}[b]{0.45\linewidth}
			\includegraphics[width=1\linewidth]{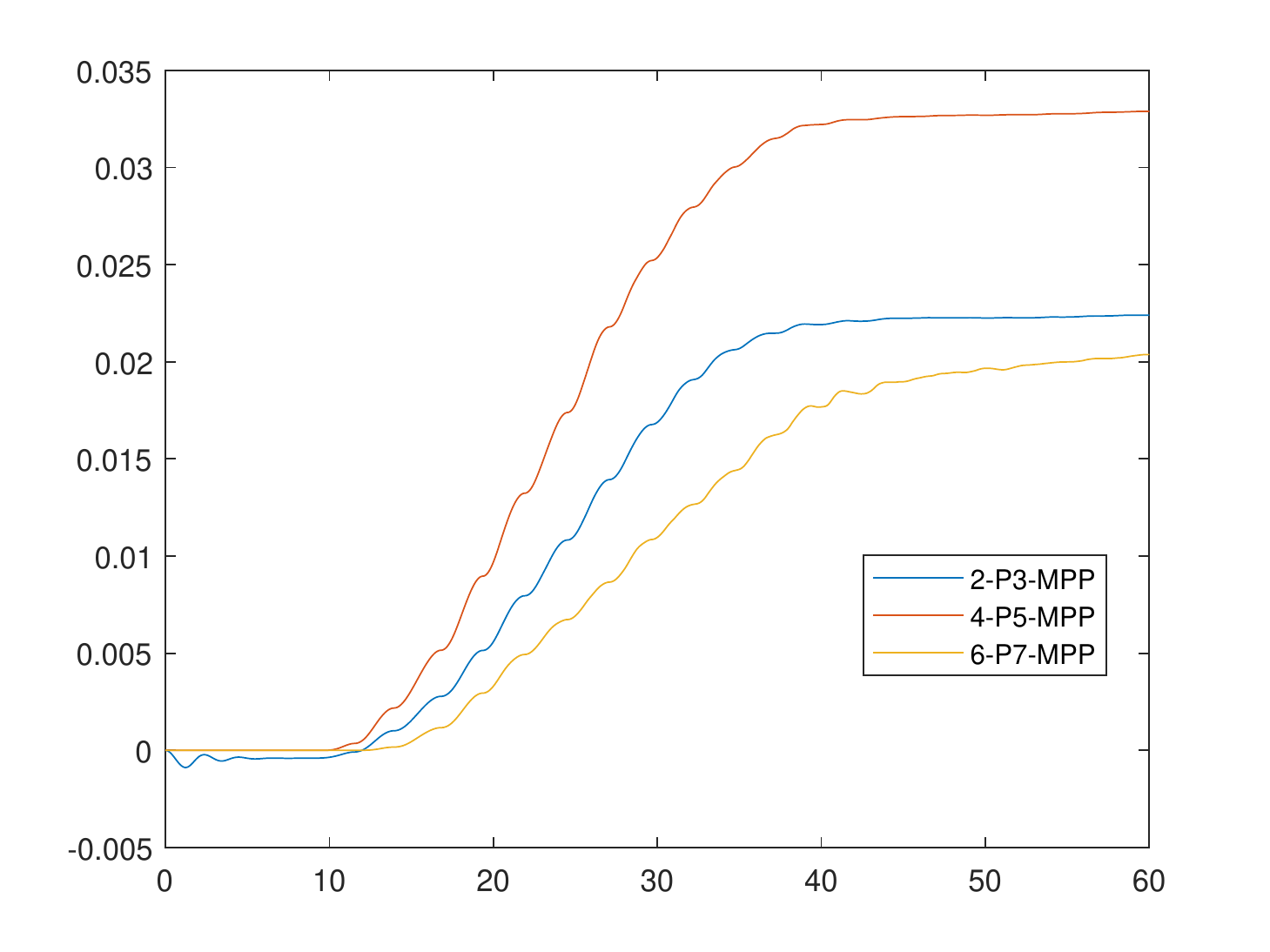}
			\subcaption{$|Energy(t=T_f)-Energy(t=0)$}
		\end{subfigure}	
		\begin{subfigure}[b]{0.45\linewidth}
			\includegraphics[width=1\linewidth]{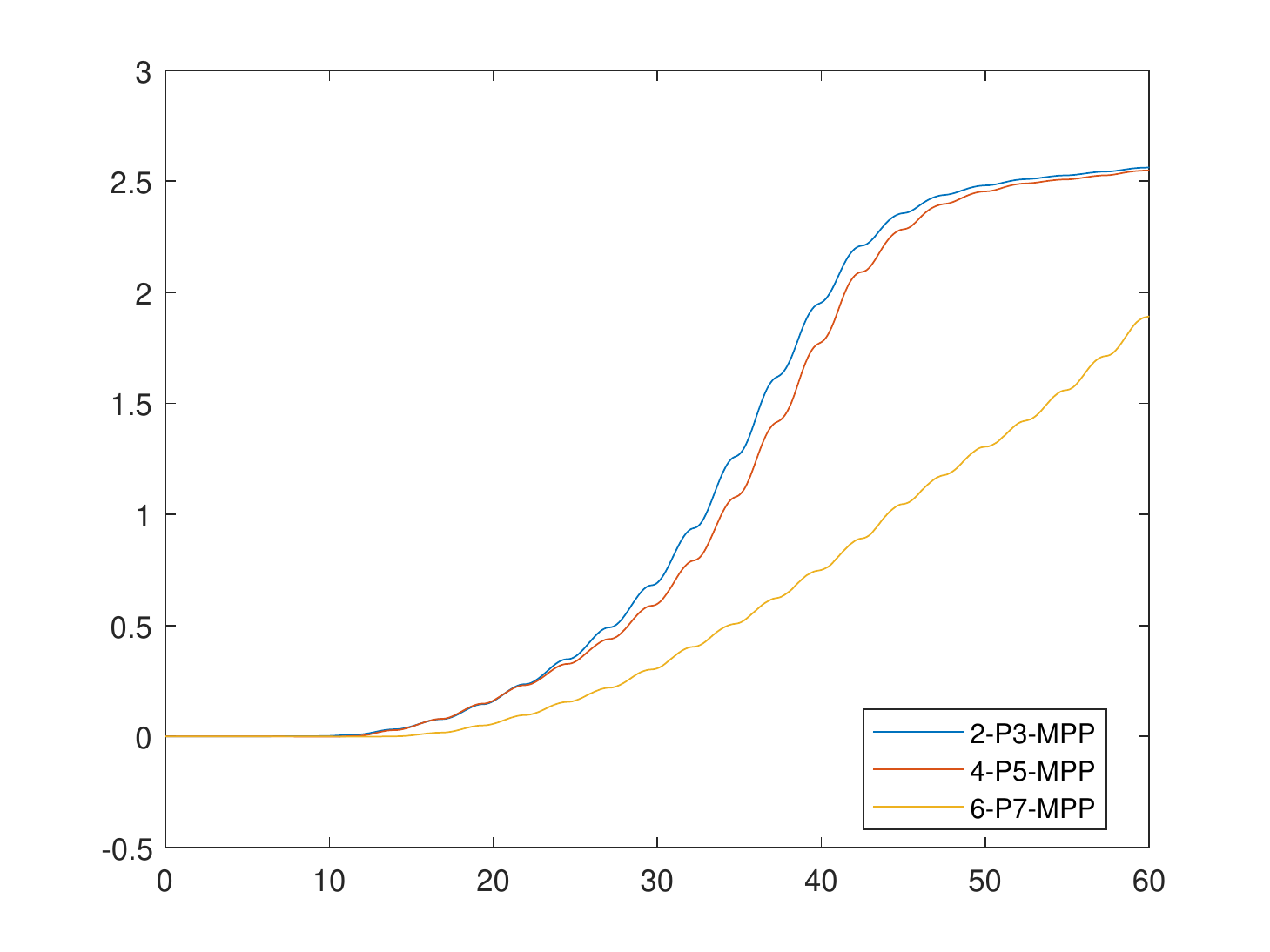}
			\subcaption{$Entropy(t=T_f)-Entropy(t=0)$}
		\end{subfigure}	
		
		\caption{1D Vlasov Poisson. Strong Landau damping with initial data \eqref{Strong Landau damping}. $N_x=128,\,N_v=256$}\label{strong landau 1}
	\end{figure}

	\begin{figure}[htbp]
		\centering
		\begin{subfigure}[b]{0.32\linewidth}
			\includegraphics[width=1.1\linewidth]{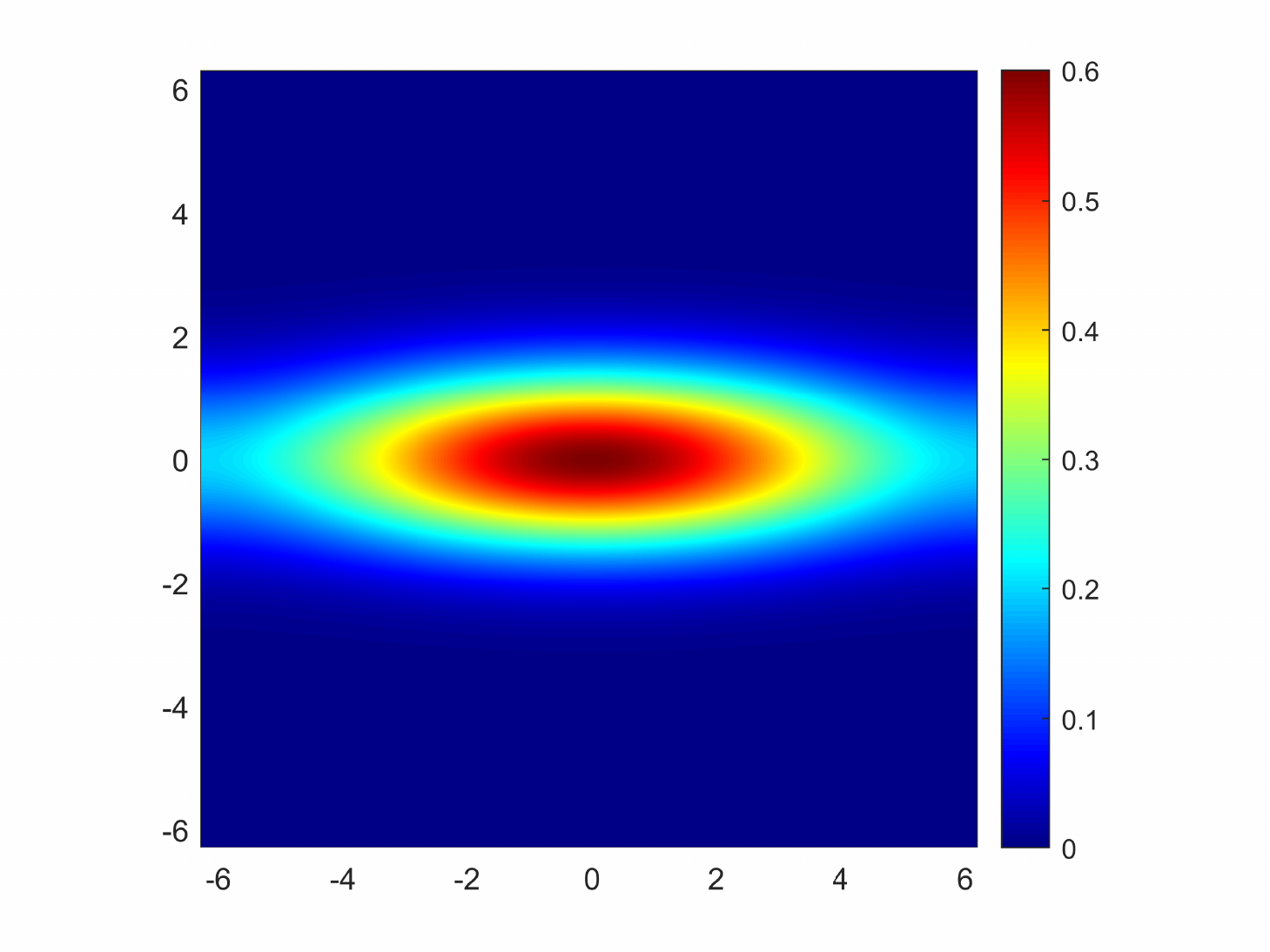}
		\end{subfigure}
		\begin{subfigure}[b]{0.32\linewidth}
			\includegraphics[width=1.1\linewidth]{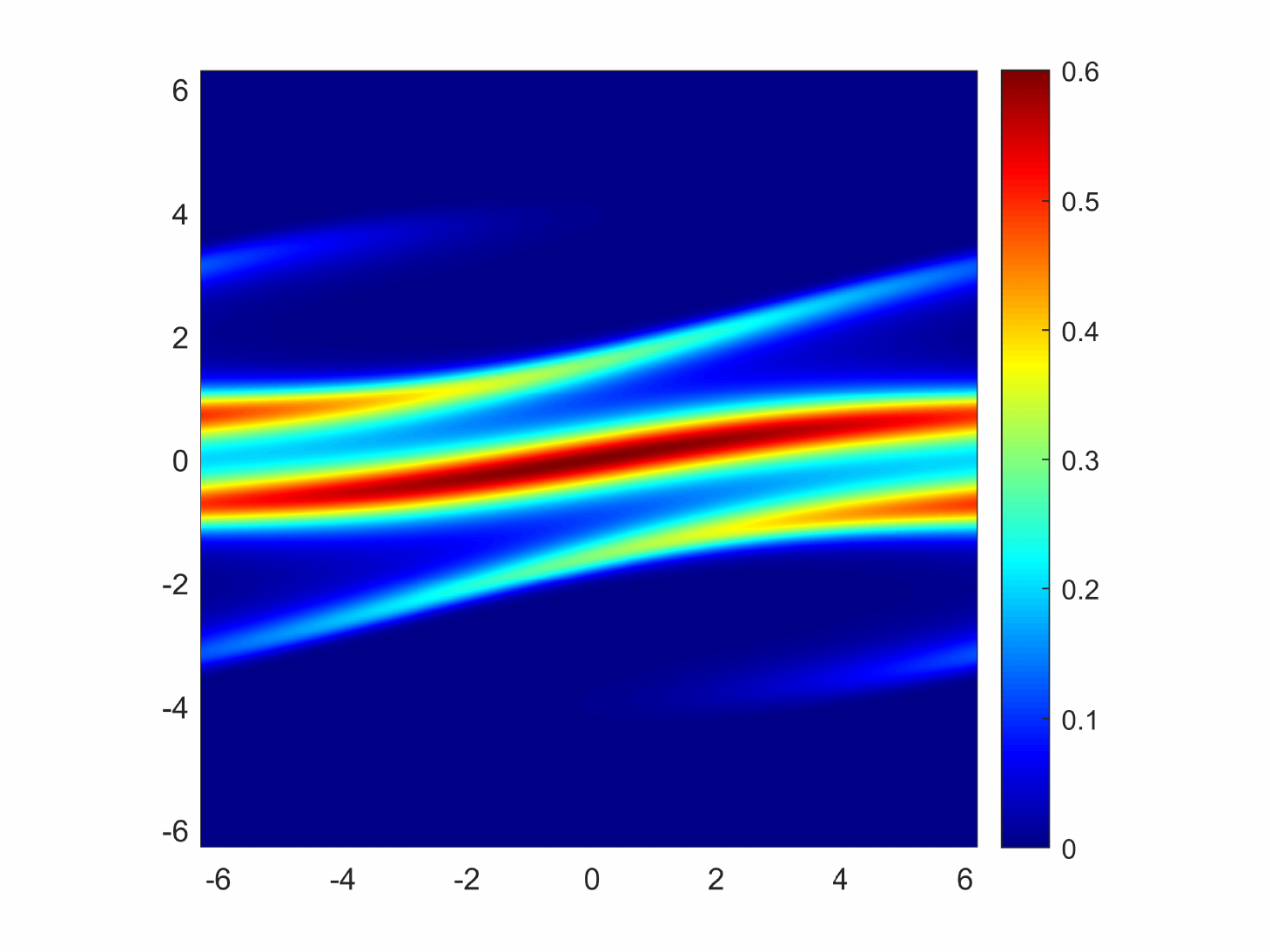}
		\end{subfigure}
		\begin{subfigure}[b]{0.32\linewidth}
			\includegraphics[width=1.1\linewidth]{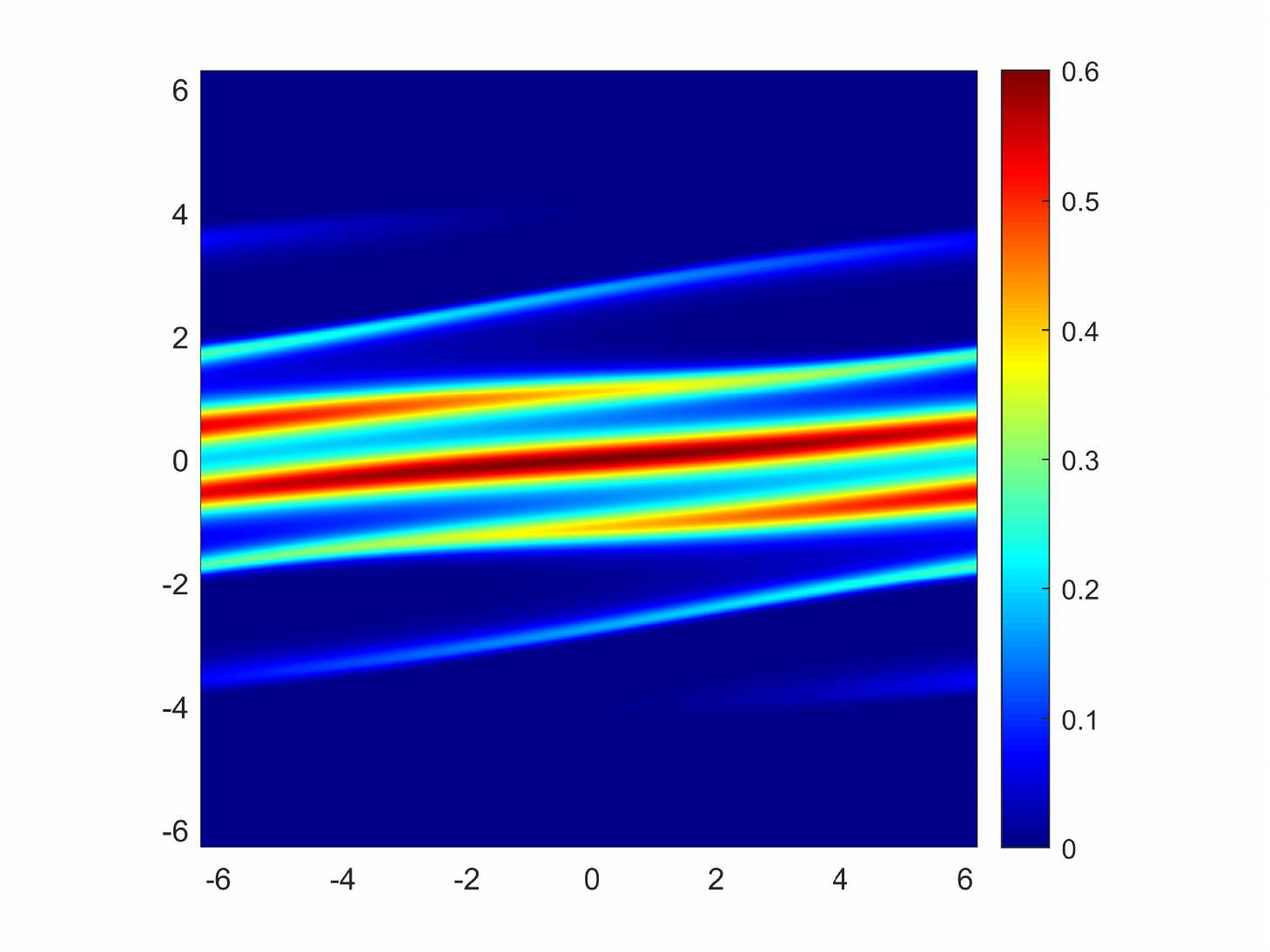}
		\end{subfigure}
		\begin{subfigure}[b]{0.32\linewidth}
			\includegraphics[width=1.1\linewidth]{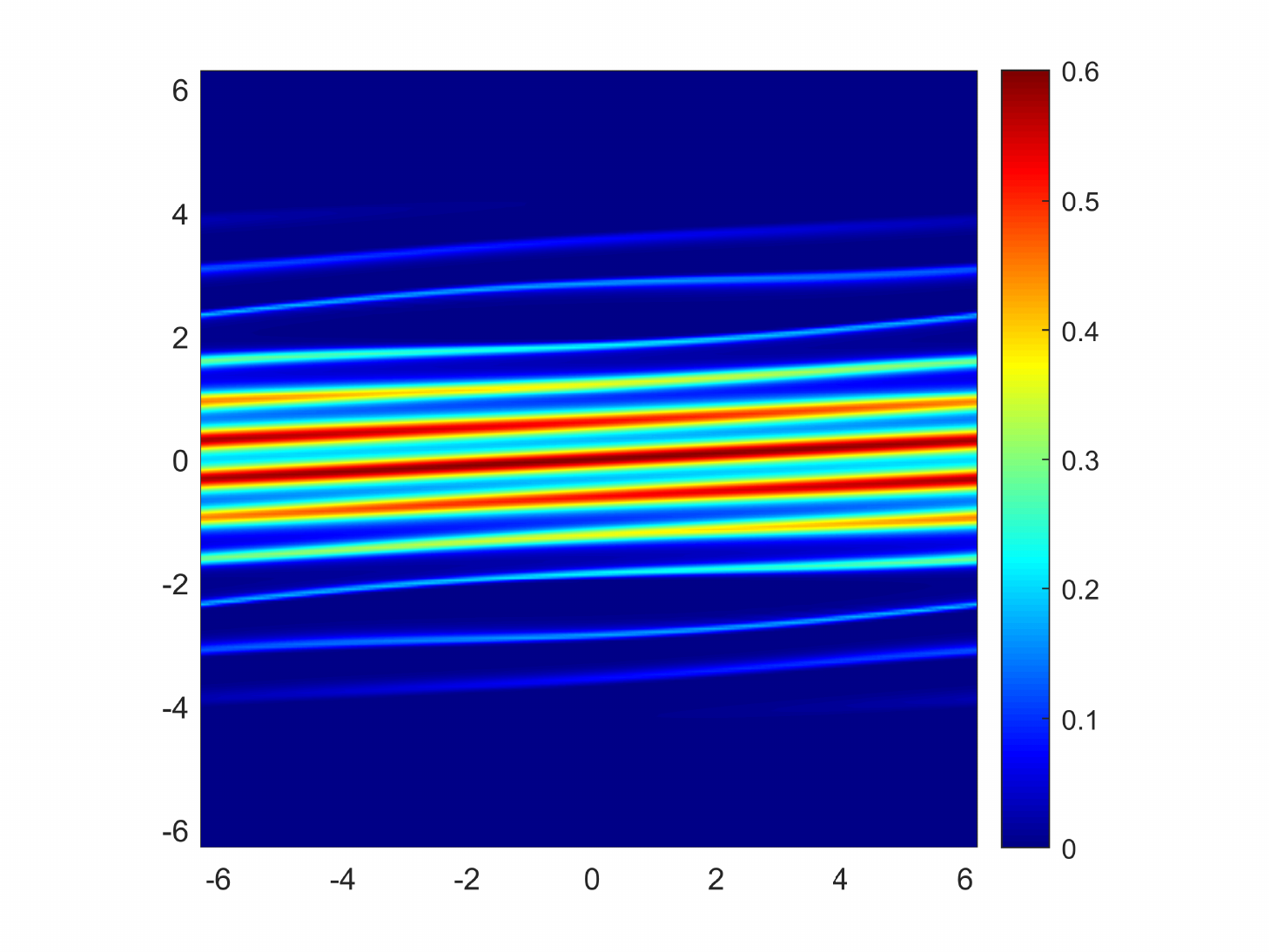}
		\end{subfigure}
		\begin{subfigure}[b]{0.32\linewidth}
			\includegraphics[width=1.1\linewidth]{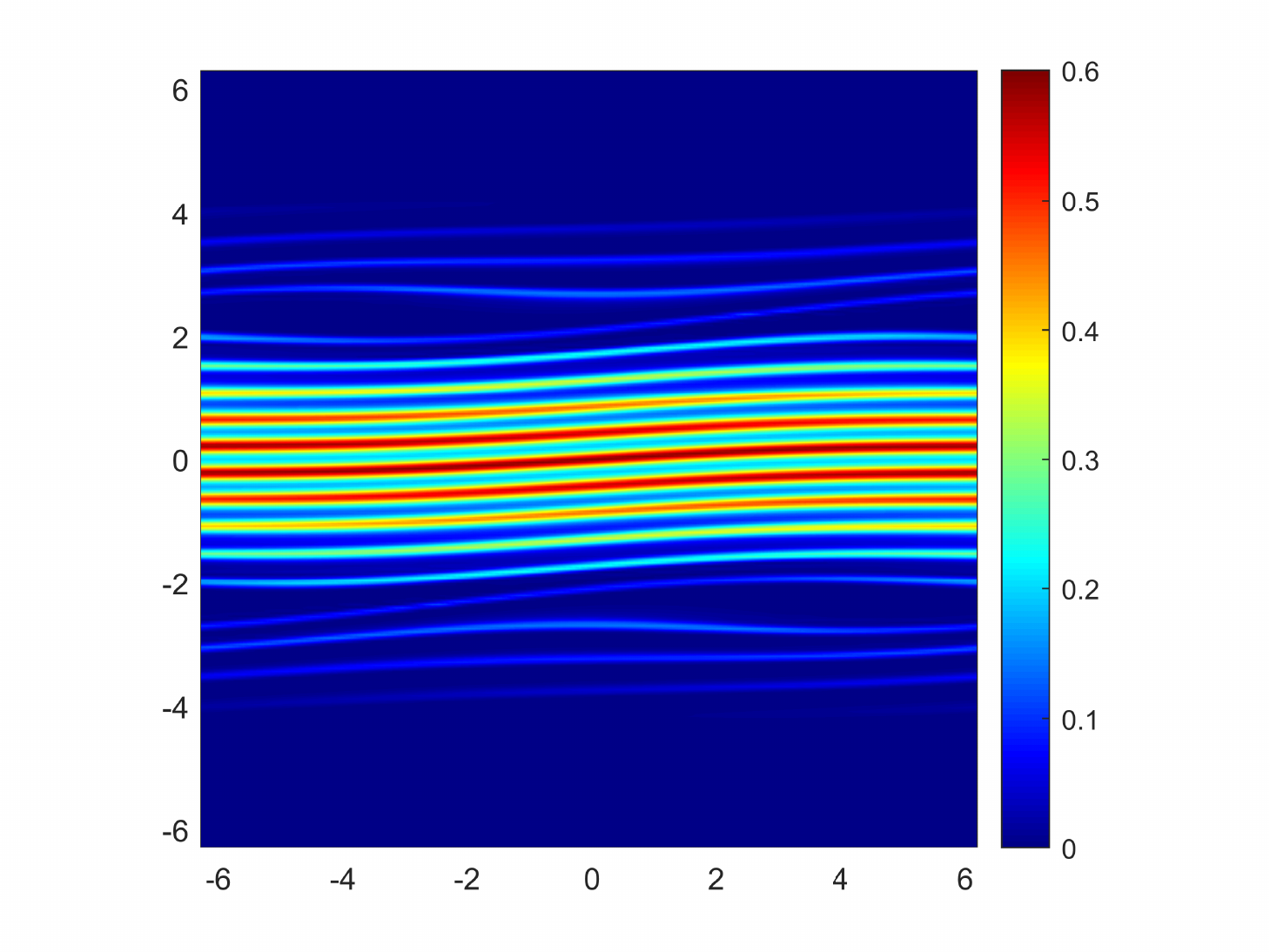}
		\end{subfigure}
		\begin{subfigure}[b]{0.32\linewidth}
			\includegraphics[width=1.1\linewidth]{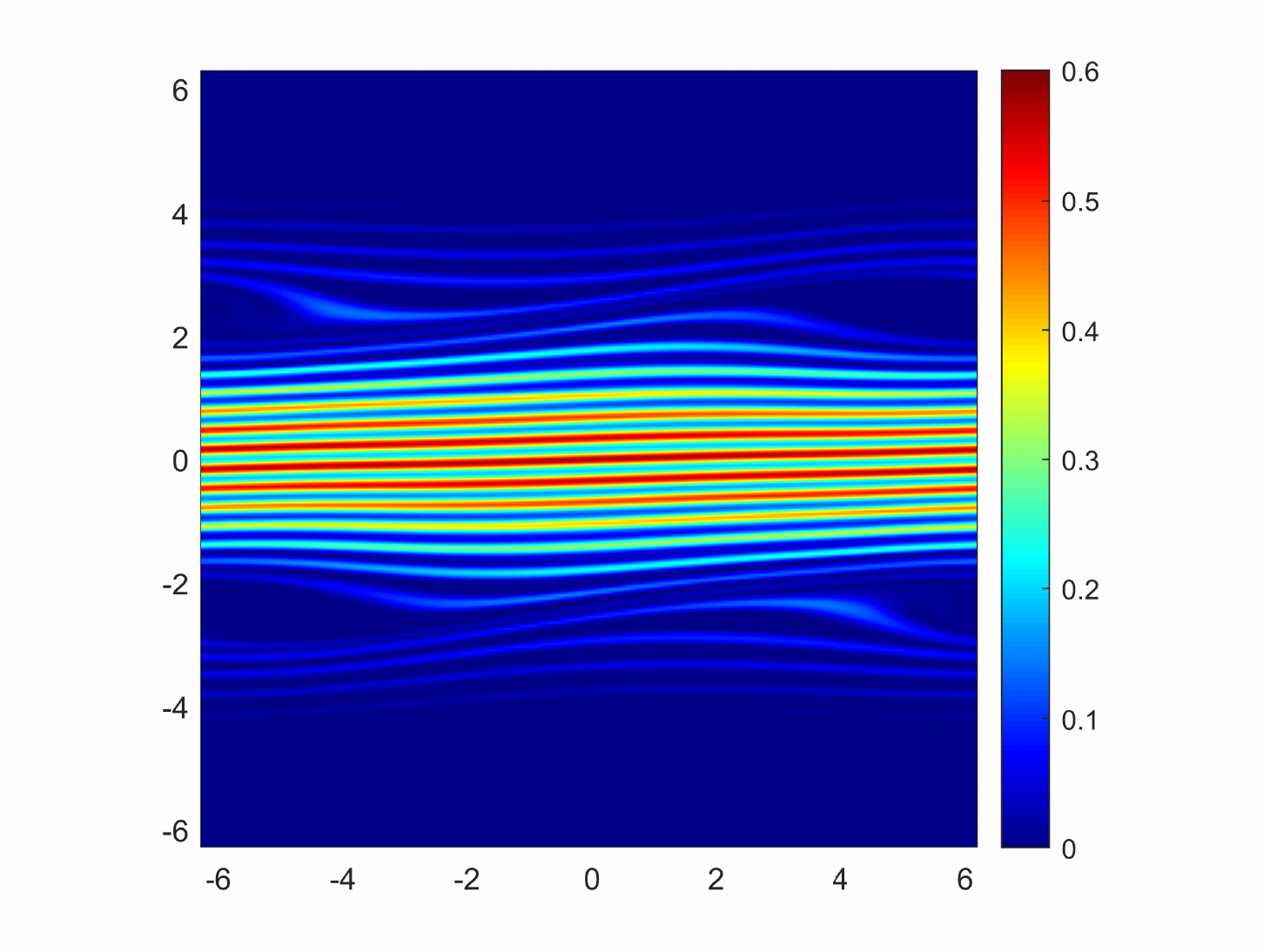}
		\end{subfigure}
		\begin{subfigure}[b]{0.32\linewidth}
			\includegraphics[width=1.1\linewidth]{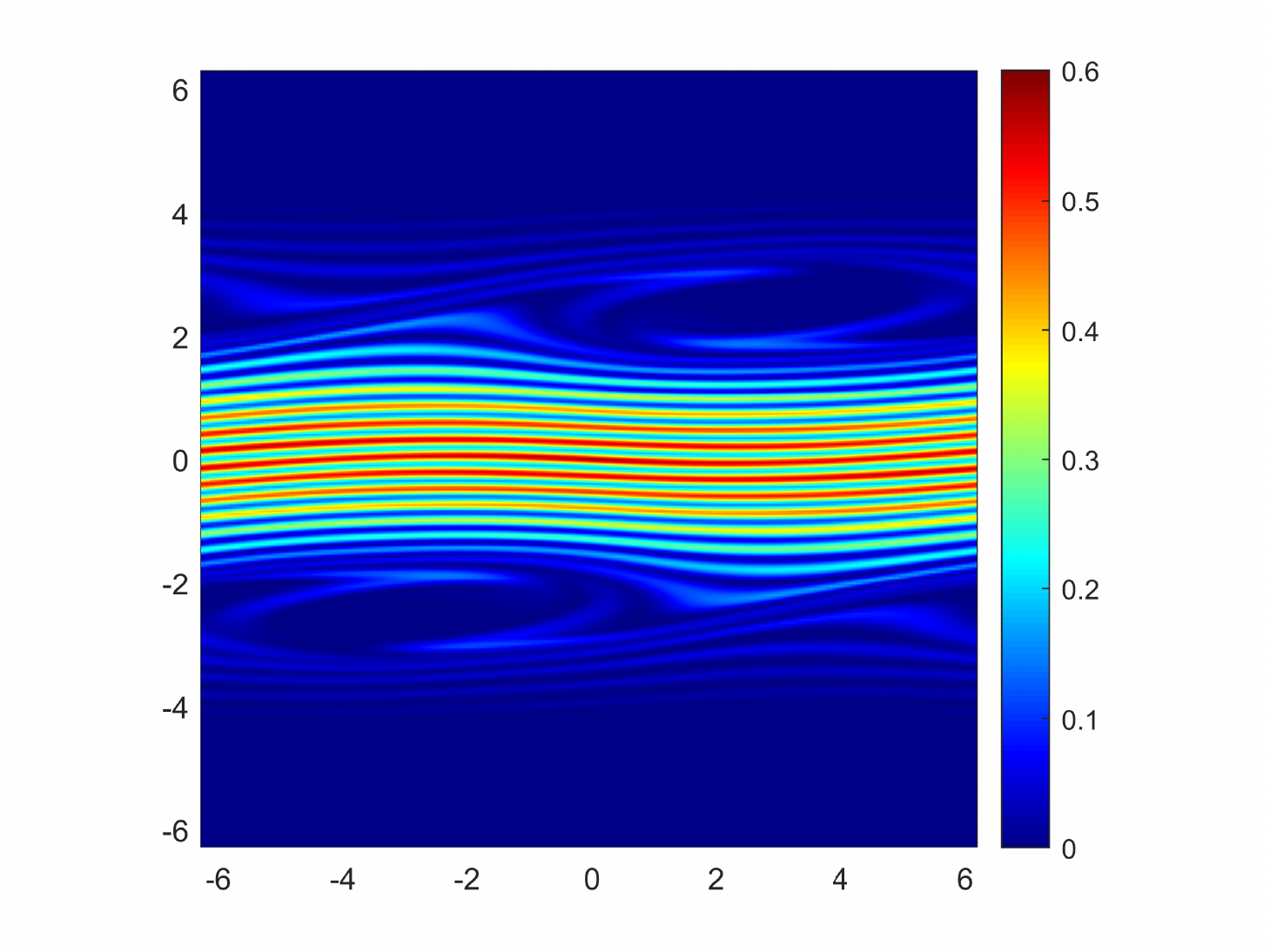}
		\end{subfigure}
		\begin{subfigure}[b]{0.32\linewidth}
			\includegraphics[width=1.1\linewidth]{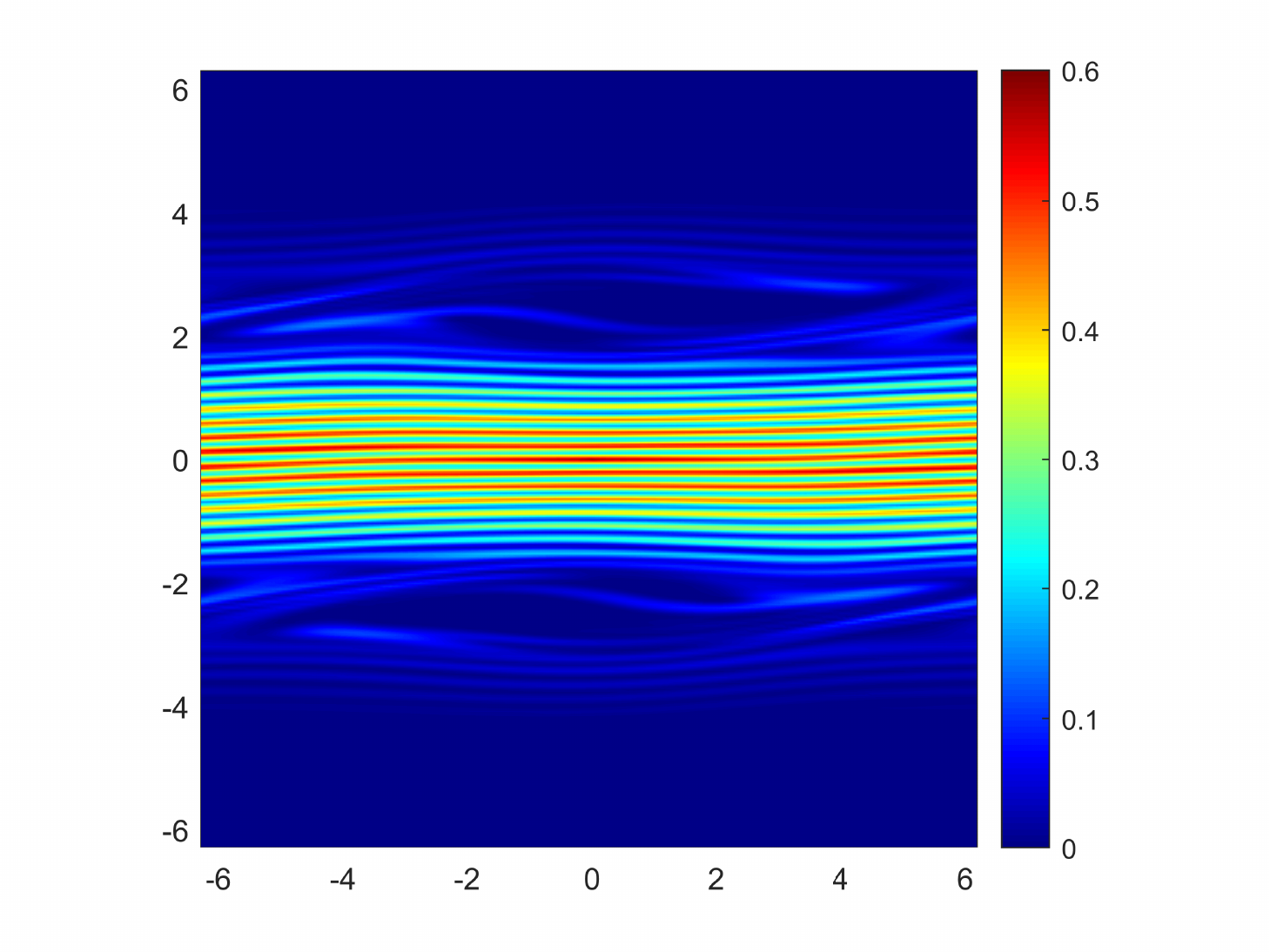}
		\end{subfigure}
		\begin{subfigure}[b]{0.32\linewidth}
			\includegraphics[width=1.1\linewidth]{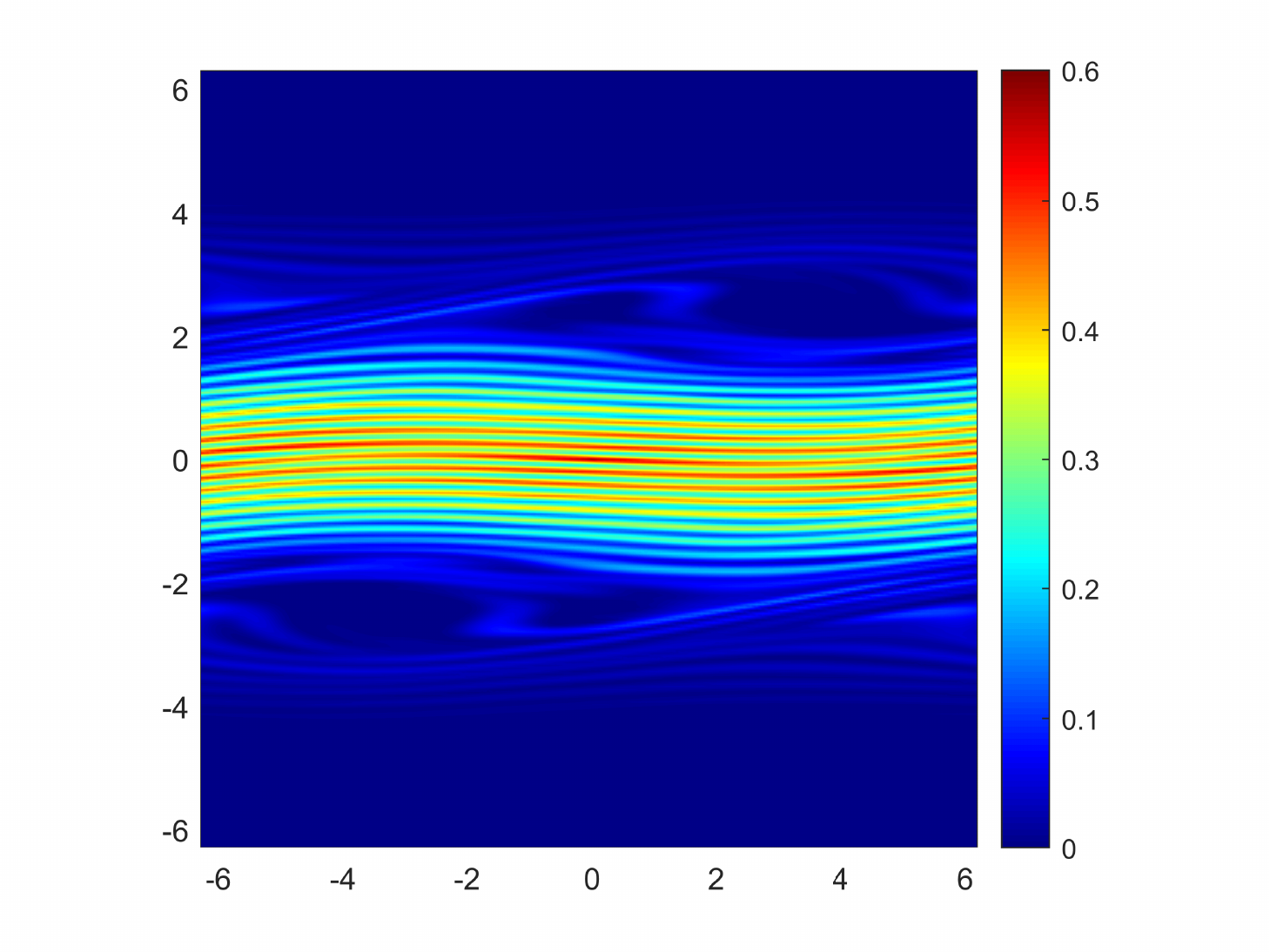}
		\end{subfigure}
		\caption{1D Vlasov Poisson. Strong Landau damping with initial data \eqref{Strong Landau damping}. Here we use a velocity domain $[-2\pi,2\pi]$ with $N_x=128$, $N_v=256$.}\label{fig strong phase space}
	\end{figure}

	\subsection{2D Vlasov-Poisson model}
	In 2D problem $(d_x, d_v)=(2,2)$, we take uniform mesh $\Delta x=\Delta y$ for space and $\Delta v$ for each velocity dimension. We use a CFL number defined by
	\begin{align}\label{CFL vp 2d}
	\text{CFL}\displaystyle :=\max\left\{ \max_j|v_j^1|\frac{\Delta t}{\Delta x} + \max_j|v_j^2|\frac{\Delta t}{\Delta y},\, \max_i|E_i^1|\frac{\Delta t}{\Delta v} + \max_i|E_i^2|\frac{\Delta t}{\Delta v} \right\},
	\end{align}
	for $E_i\equiv (E_i^1,E_i^2)$.
	
	\subsubsection{Accuracy test}\label{sec vp 2d acc}
	To check the accuracy of splitting methods for 2D Vlasov-Poisson system, we consider a setup in \cite{FSB-vlasov-2001}, where the initial data is given by 
	\begin{align}\label{weak Ladau 2d}
	f(x,y,v^1,v^2,0)= \frac{1}{2\pi}\bigg(1+ \alpha \big(\cos(k x) + \cos(k y)\big)\bigg) \exp\left(-\frac{|v^1|^2+|v^2|^2}{2}\right),\quad \alpha=0.05,\quad k=0.5.
	\end{align}	
	We impose periodic boundary condition on the space domain $[0, 4\pi]^2$ and assume {zero-boundary condition} on velocity domain $[-6, 6]$. Numerical solutions are computed up to final time $T_f=2$ taking uniform grids with $N_x=N_y=N_{v_1}=N_{v_2}=20,40,80$ and different time steps based on CFL$=1,2,3,4,5$ using \eqref{CFL vp 2d}. 
	
	In this test, we consider a semi-Lagrangian method based on 2nd order Strang splitting method, i.e. in the advection step we solve Eq.\eqref{VP2} and in the drift step we solve Eq.\eqref{VP3}. Note that we can solve each step without dimensional splitting, because during the advection step the velocities are constant in time, and during the drift step the electric field is constant in time.  For a basic reconstruction, together with MPP limiter, we again take 2D-P3, which is the two-dimensional optimal polynomial of degree 2 used for CWENO23 reconstruction \cite{LPR1}.
	\begin{center}
		\begin{table}[ht]
			\centering
			{\begin{tabular}{|ccccccccccc|}
					\hline
					\multicolumn{1}{ |c| }{}& \multicolumn{2}{ c  }{CFL$=1$} & \multicolumn{2}{ |c }{CFL$=2$}& \multicolumn{2}{ |c }{CFL$=3$} &
					\multicolumn{2}{ |c }{CFL$=4$}&
					\multicolumn{2}{ |c| }{CFL$=5$} \\ \hline
					\multicolumn{1}{ |c|  }{$(N_x^2,(2N_x)^2)$} &
					\multicolumn{1}{ c  }{error} &
					\multicolumn{1}{ c|  }{rate} &
					\multicolumn{1}{ |c  }{error} &
					\multicolumn{1}{ c  }{rate} &
					\multicolumn{1}{ |c  }{error} &
					\multicolumn{1}{ c|  }{rate} &
					\multicolumn{1}{ |c  }{error} &
					\multicolumn{1}{ c|  }{rate}&
					\multicolumn{1}{ |c  }{error} &
					\multicolumn{1}{ c|  }{rate}     \\ 
					\hline
					\hline	
					\multicolumn{1}{ |c|  }{$(20^2,40^2)$}&1.19e-03
					&2.8278          
					&1.24e-03
					&2.73&1.43e-03
					&2.46&2.22e-03
					&2.47&2.55e-03
					&2.11
					\\
					\multicolumn{1}{ |c|  }{$(40^2,80^2)$}&1.68e-04   
					&&1.87e-04&&2.61e-04&&4.01e-04&&5.92e-04&
					\\
					\hline
					\hline
			\end{tabular}}
			\caption{Accuracy test for the 2D Vlasov-Poisson system. Initial data is given in \eqref{weak Ladau 2d}.}\label{tab4}
		\end{table}
	\end{center}
	
	Table \ref{tab4} shows that, for smaller CFL numbers, numerical errors tend to be smaller, and the desired accuracy of spatial reconstruction appears. On the contrary, as CFL number gets bigger, the expected accuracy of a time splitting method is observed. As in 1D case, the use MPP limiter does not lead to order reduction.

	\subsubsection{2D Long time simulations}
	For the 2D Vlasov-Poisson model, we again consider a strong Landau damping problem as in \cite{crouseilles2008comparison}:
	\begin{align}\label{Strong Ladau 2d}
	f(x,y,v^1,v^2,0)= \frac{1}{2\pi}\bigg(1+ \alpha \big(\cos(k x) + \cos(k y)\big)\bigg) \exp\left(-\frac{|v^1|^2+|v^2|^2}{2}\right),\quad \alpha=0.5,\quad k=0.5.
	\end{align}	
	We impose periodic boundary condition in space domain $[0, 4\pi]^2$ and zero-boundary condition in velocity domain $[-6,6]^2$. Using the scheme based on Strang splitting and 2D-P3, numerical solutions are computed up to final time $T_f=40$. Compared to literature \cite{crouseilles2008comparison}, we here take a larger fixed time step $\Delta t = \frac{1}{8}$, and fewer grid points in space by setting $N_x=N_y=32$ and $N_{v^1}=N_{v^2}=128$.
	
	As in 1D case, Figure \ref{2d vp strong landau} shows that a phase transition appears for the $L^2$-norm and $L^\infty$-norm of electric field during time evolution. Even with fewer grid points and large time step, our result shows good agreement with the results in \cite{crouseilles2008comparison}. Here, due to the truncation of velocity domain the $L^1$-norm of the numerical solution is only maintained upto $\mathcal{O}(10^{-6})$, and can be reduced upto machine precision when we take a larger velocity domain. The other conservative quantities such as $L^2$-norm, entropy, energy are not fully conserved at a numerical level.

	\begin{figure}[htbp]
		\centering
		\begin{subfigure}[b]{0.45\linewidth}
			\includegraphics[width=1\linewidth]{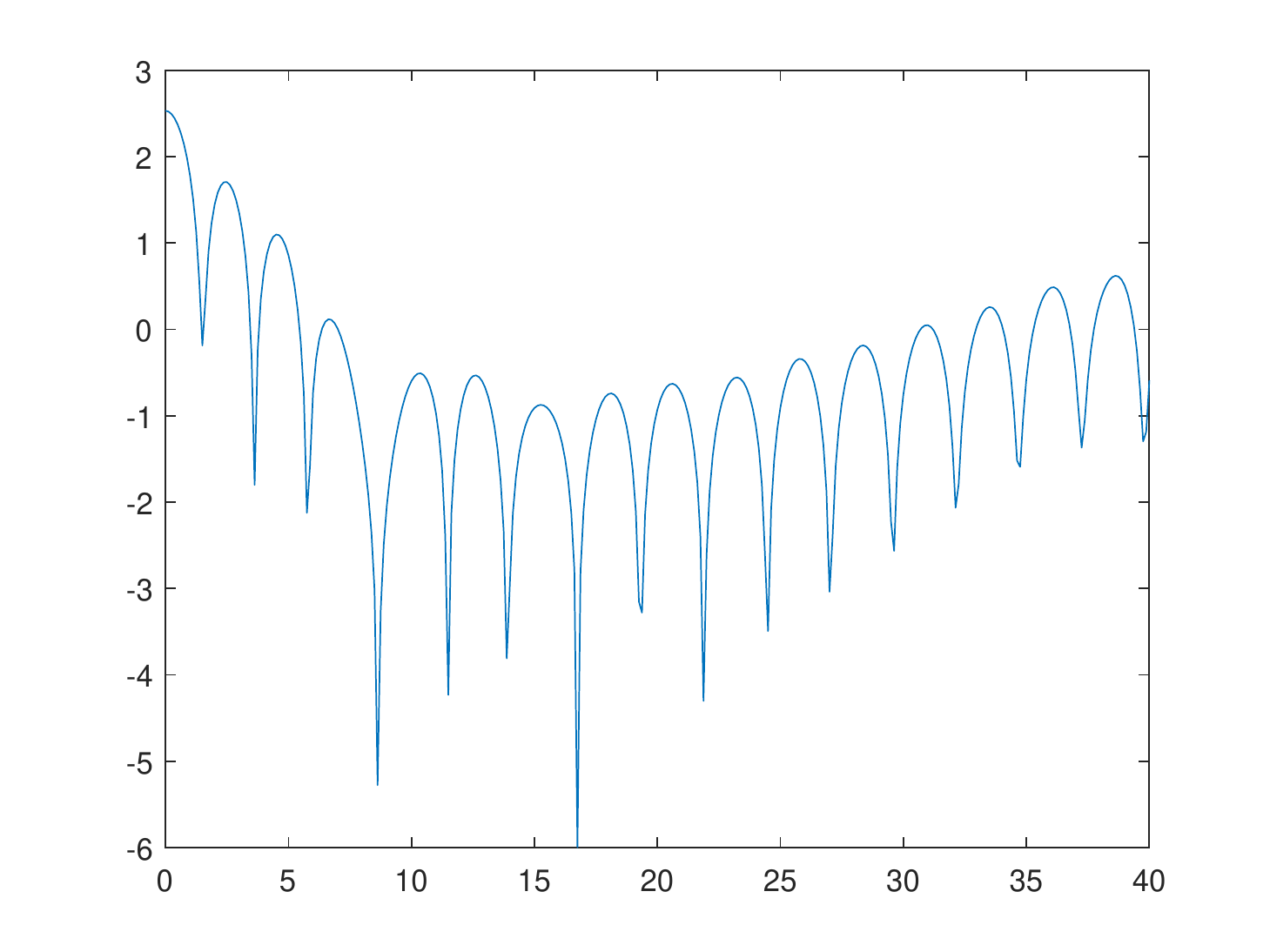}
			\subcaption{Logarithmic scale of $L^2$-norm of $E$}
		\end{subfigure}	
		\vspace*{5mm}
		\begin{subfigure}[b]{0.45\linewidth}
			\includegraphics[width=1\linewidth]{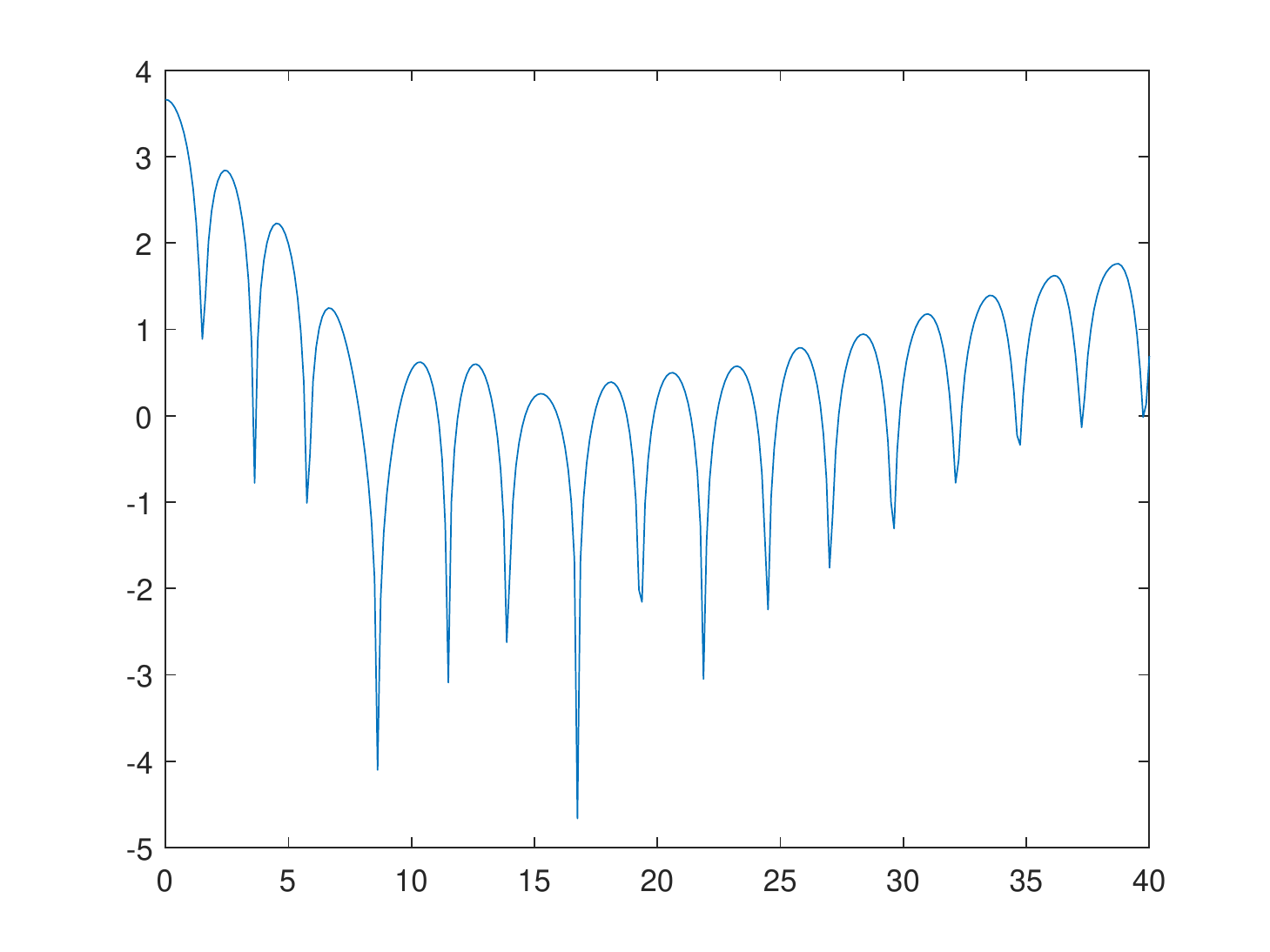}
			\subcaption{Logarithmic scale of $L^{\infty}$-norm of $E$}
		\end{subfigure}
		\vspace*{5mm}	
		\begin{subfigure}[b]{0.45\linewidth}
			\includegraphics[width=1\linewidth]{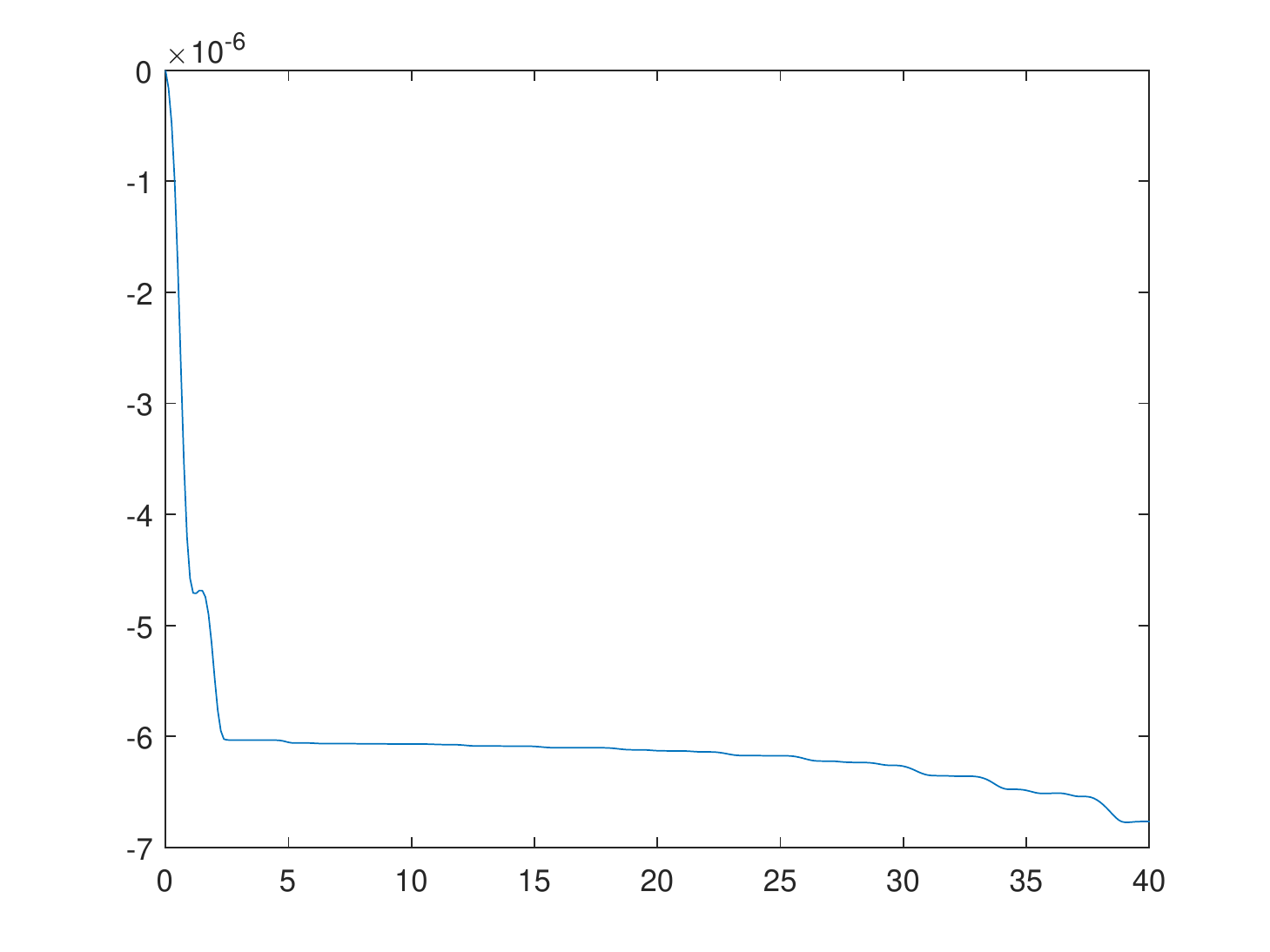}
			\subcaption{$\|f^{N_t}\|_1-\|f^0\|_1$}
			
		\end{subfigure}	
		\begin{subfigure}[b]{0.45\linewidth}
			\includegraphics[width=1\linewidth]{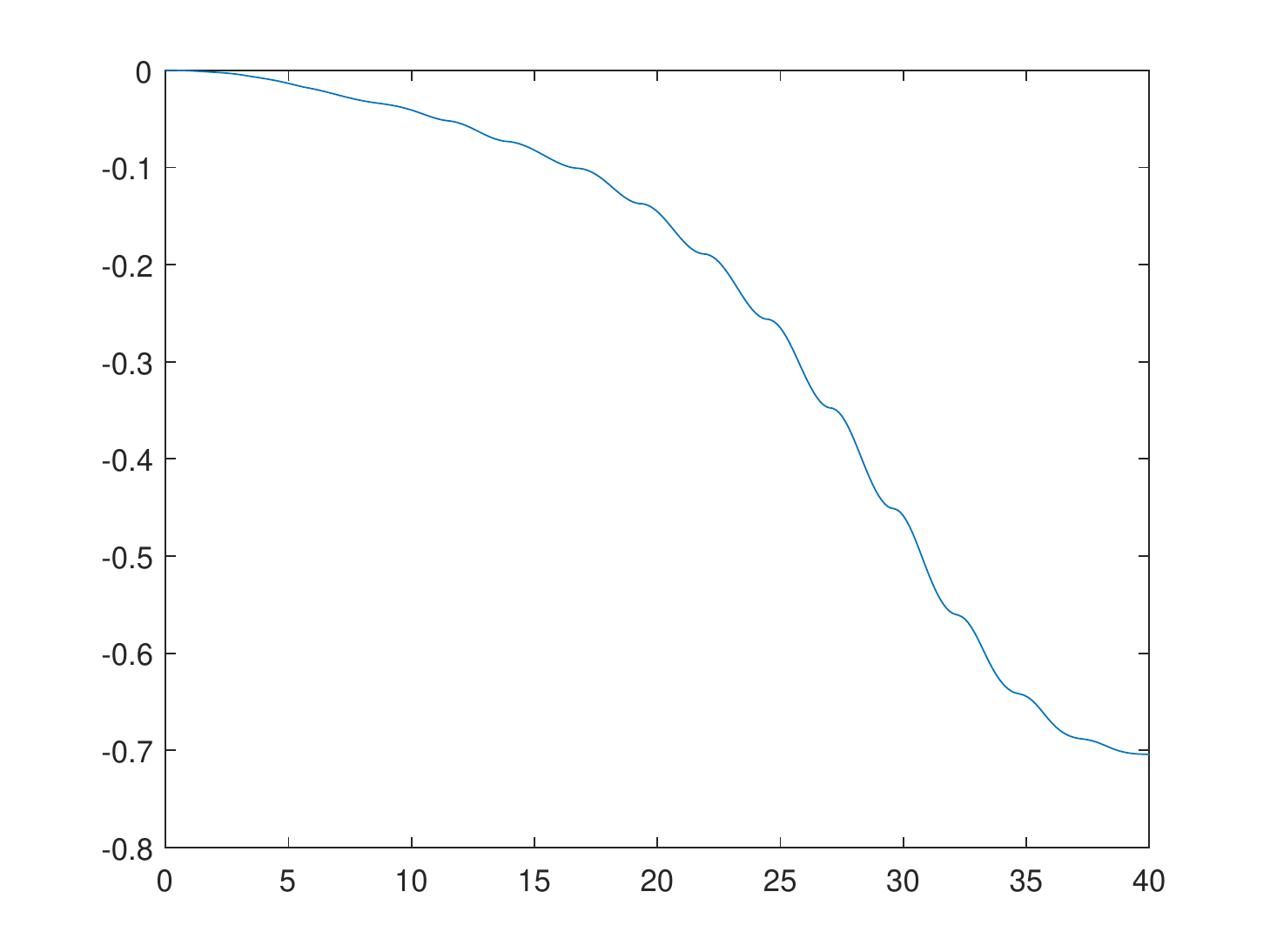}
			\subcaption{$\|f^{N_t}\|_2-\|f^0\|_2$}
		\end{subfigure}	
		\begin{subfigure}[b]{0.45\linewidth}
			\includegraphics[width=1\linewidth]{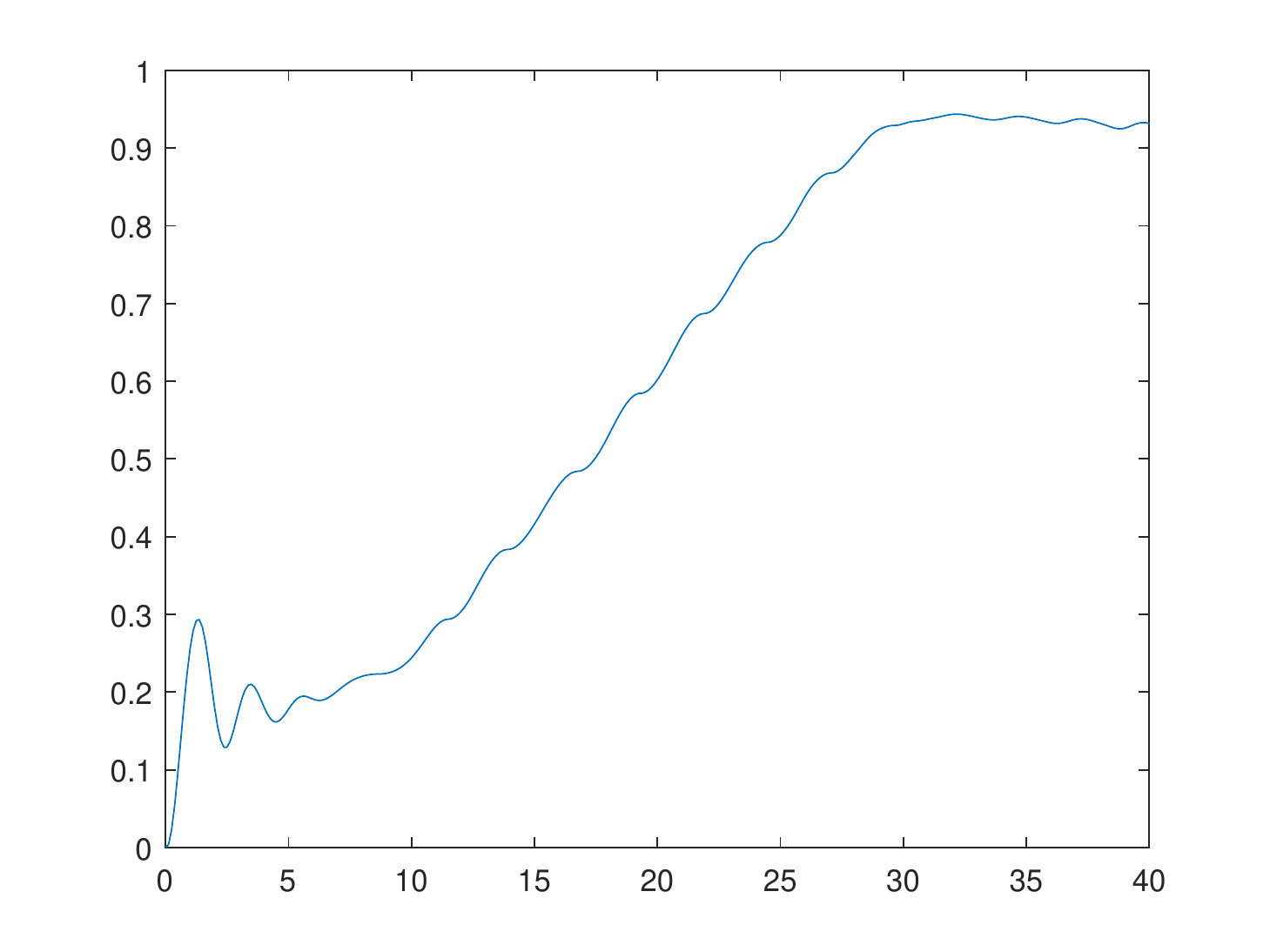}
			\subcaption{$|Energy(t=T_f)-Energy(t=0)$}
		\end{subfigure}	
		\begin{subfigure}[b]{0.45\linewidth}
			\includegraphics[width=1\linewidth]{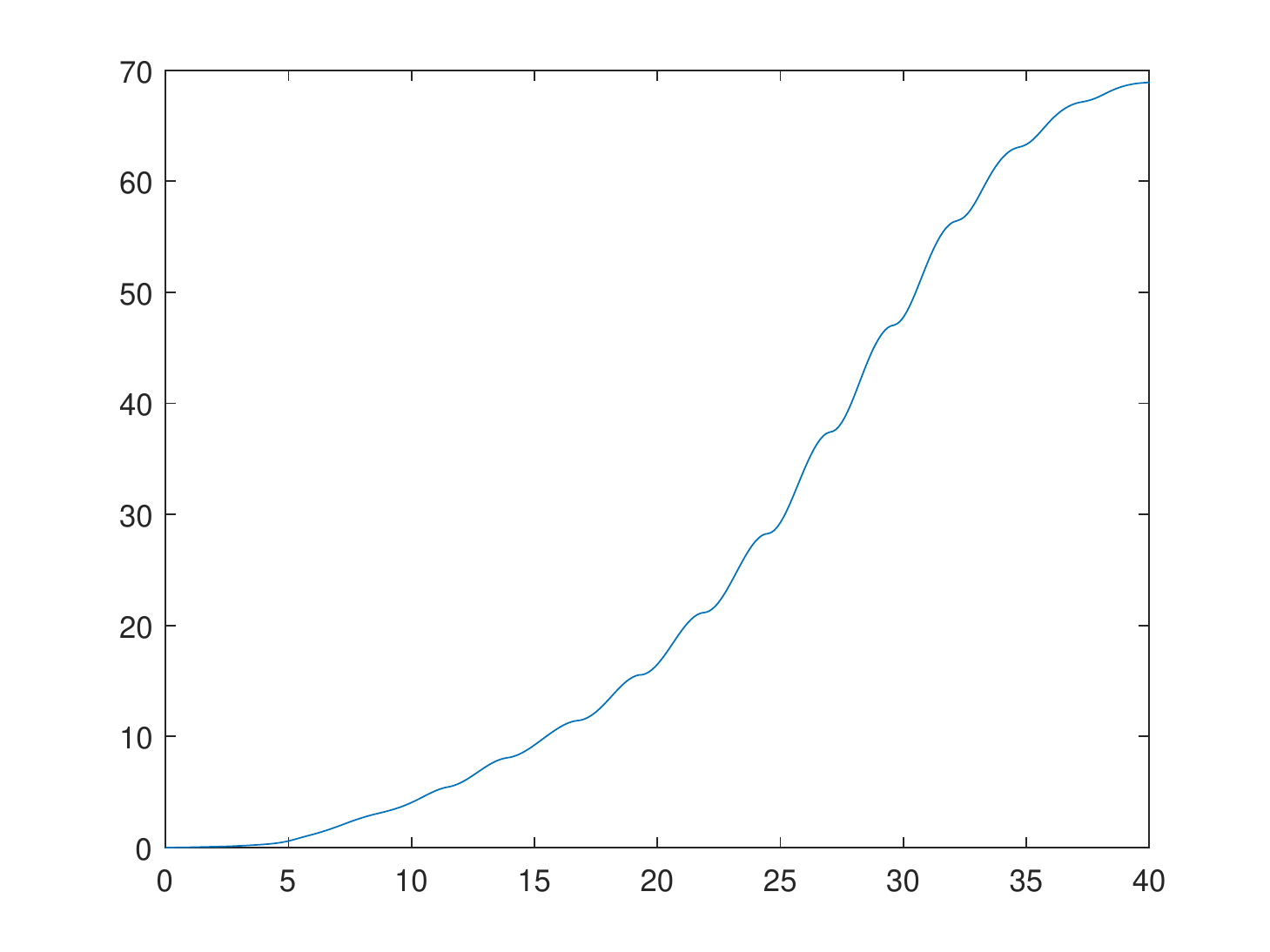}
			\subcaption{$Entropy(t=T_f)-Entropy(t=0)$}
		\end{subfigure}	
		\caption{2D Vlasov Poisson. Strong Landau damping with initial data \eqref{Strong Ladau 2d}. $N_x=N_y=32,\,N_v^1=N_v^2=128$.}\label{2d vp strong landau}
	\end{figure}

	\subsection{1D BGK model}
	Now, we move on to the numerical tests of SL schemes for the BGK model. After checking the accuracy of non-splitting semi-Lagrangian methods, we treat a related shock problem arising in the fluid limit $\kappa \to 0$.
	In 1D problem $(d_x, d_v)=(1,1)$, we again consider uniform mesh $\Delta x$ and $\Delta v$ for space and velocity domain. For 1D BGK model, we use a CFL number defined by 
	\begin{align}\label{CFL bgk 1d}
	\text{CFL}\displaystyle := \max_j|v_j|\frac{\Delta t}{\Delta x}.
	\end{align}
	\begin{center}
		\begin{table}[ht]
			\centering
			{\begin{tabular}{|ccc|}
					\hline
					\multicolumn{1}{ |c }{}&
					\multicolumn{1}{ |c| }{Time time discretization}& \multicolumn{1}{ c|  }{Basic reconstruction}   \\ \hline
					\multicolumn{1}{ |c }{RK2}&
					\multicolumn{1}{ |c|  }{2nd order DIRK method in \eqref{Butcher}}     
					&\multicolumn{1}{ |c|  }{CWENO23}
					\\ \hline
					\multicolumn{1}{ |c }{BDF2}&
					\multicolumn{1}{ |c|  }{2nd order BDF method in \eqref{BDFschemes}}     
					&\multicolumn{1}{ |c|  }{CWENO23}
					\\
					\hline
					\multicolumn{1}{ |c }{RK3}&
					\multicolumn{1}{ |c|  }{3th order DIRK method in \eqref{Butcher}}     
					&\multicolumn{1}{ |c|  }{CWENO35}
					\\
					\hline
					\multicolumn{1}{ |c }{BDF3}&
					\multicolumn{1}{ |c|  }{3rd order BDF method in \eqref{BDFschemes}}     
					&\multicolumn{1}{ |c|  }{CWENO35}
					\\
					\hline
			\end{tabular}}
			\caption{Numerical methods without MPP limiter used for 1D BGK model.}\label{tab bgk name}
		\end{table}
	\end{center}

	\subsubsection{Accuracy test}
	To check the accuracy, we consider the accuracy test in \cite{GRS}. The initial distribution is given by the Maxwellian 
	\begin{align}\label{bgk 1d initial}
	f_0(x,v)=\frac{\rho_0}{\sqrt{2 \pi T_0}}\exp\left(-\frac{|v-u_0(x)|^2}{2T_0}\right),
	\end{align}
	where initial density and temperature are assumed to be uniform, with constant value $\rho_0(x) \equiv 1 $ and $T_0(x) \equiv 1$. Initial velocity profile is given by
	\[
	u_0(x) = 0.1 \exp\left(-(10x - 1)^2\right) - 2 \exp\left(-(10x + 3)^2\right).
	\]
	For space, we assume periodic boundary conditions in the interval $[-1,1]$ with $N_x = 320,\, 640,\, 1280,\, 2560$ and $5120$. For velocity domain, we consider $[-10,10]$ with $N_v = 20$, which is enough when we use correction \eqref{first dm}, as described in Section \ref{sec:CM}. We set time steps based on CFL numbers \eqref{CFL bgk 1d}, and compute numerical solutions up to $T_f=0.32$  because for small Knudsen number shock appears at $t = 0.35$.
	
	
	For RK2 and BDF2 based methods, we use 1D CWENO23 as a basic reconstruction, and hence the expected convergence rate is between 2 and 3. For RK3 and BDF3 based methods using CWENO35 as a basic reconstruction, the convergence rate is expected to be between 3 and 5. In the test, we take CFL$=4$ for RK2, RK3 and BDF2 based methods, and CFL$=1$ for BDF3 based method. Even if larger CFL numbers are allowed by stability, the reason for this choice is to confirm the expected accuracy of time discretization and spatial reconstructions.
	
	In Table \ref{1d accuracy}, we report errors and convergence rate computed with the relative $L^1$-norm for density $\rho$. (See \eqref{rel norm} and \eqref{rate}.) The result shows that RK2 and BDF2 based methods attain the expected order between 2 and 3. Although order reductions appear for RK3 based method as $\kappa \to 0$, it attains the expected order 5 of reconstruction for $\kappa = 1$. On the contrary, BDF3 does not suffer from order reduction and gives the desired convergence rate between 3 and 5 for all values of $\kappa$.
	We remark that similar results are obtained for other macroscopic quantities such as momentum and energy, but we omit to show them here.
	\begin{center}
		\begin{table}[htbp]
			\centering
			{\begin{tabular}{|cccccccccc|}
					\hline
					\multicolumn{10}{ |c| }{Relative $L^1$ error and order of density} \\ \hline
					\multicolumn{1}{ |c }{}&
					\multicolumn{1}{ |c| }{}& \multicolumn{2}{ c  }{$\kappa=10^{-6}$} & \multicolumn{2}{ |c }{$\kappa=10^{-4}$}& \multicolumn{2}{ |c| }{$\kappa=10^{-2}$} &
					\multicolumn{2}{ |c| }{$\kappa=10^{-0}$} \\ \hline
					\multicolumn{1}{ |c }{}&
					\multicolumn{1}{ |c|  }{$(N_x,2N_x)$} &
					\multicolumn{1}{ c  }{error} &
					\multicolumn{1}{ c|  }{rate} &
					\multicolumn{1}{ |c  }{error} &
					\multicolumn{1}{ c  }{rate} &
					\multicolumn{1}{ |c  }{error} &
					\multicolumn{1}{ c|  }{rate} &
					\multicolumn{1}{ |c  }{error} &
					\multicolumn{1}{ c|  }{rate}     \\ 
					\hline
					\hline
					\multicolumn{1}{ |c }{}&
					\multicolumn{1}{ |c|  }{$(320,640)$}&2.29e-05
					&2.20	    
					&1.91e-05	
					&2.35
					&2.36e-06
					&2.27
					&5.31e-06	
					&3.00
					\\
					\multicolumn{1}{ |c }{RK2}&
					\multicolumn{1}{ |c|  }{$(640,1280)$}&4.97e-06
					&2.07          
					&3.75e-06
					&2.26&4.90e-07
					&2.14&6.66e-07
					&3.00
					\\
					\multicolumn{1}{ |c }{}&
					\multicolumn{1}{ |c|  }{$(1280,2560)$}&1.18e-06   
					&&7.84e-07&&1.11e-07&&8.34e-08&
					\\
					\hline
					\multicolumn{1}{ |c }{}&
					\multicolumn{1}{ |c|  }{$(320,640)$}&7.72e-05
					&2.02    
					&7.39e-05
					&2.03
					&1.48e-05
					&2.07
					&8.33e-06	
					&2.99
					
					\\
					\multicolumn{1}{ |c }{BDF2}&
					\multicolumn{1}{ |c|  }{$(640,1280)$}&1.90e-05		
					&2.01        
					&1.81e-05
					&2.01&3.53e-06
					&2.04&1.05e-06
					&3.00
					\\
					\multicolumn{1}{ |c }{}&
					\multicolumn{1}{ |c|  }{$(1280,2560)$}&4.72e-06  
					&&4.50e-06&&8.60e-07&&1.31e-07&
					\\
					\hline
					\multicolumn{1}{ |c }{}&
					\multicolumn{1}{ |c|  }{$(640,1280)$}&4.88e-06
					&2.00      
					&3.73e-06
					&2.17&1.87e-07
					&2.73&1.53e-09
					&4.96
					\\
					\multicolumn{1}{ |c }{RK3}&
					\multicolumn{1}{ |c|  }{$(1280,2560)$}&1.22e-06
					&2.01        
					&8.28e-07
					&2.19&2.83e-08
					&2.84&4.92e-11
					&4.85
					\\
					\multicolumn{1}{ |c }{}&
					\multicolumn{1}{ |c|  }{$(2560,5120)$}&3.04e-07 
					&&1.82e-07
					&&3.96e-09&&1.72e-12&
					\\
					\hline
					\multicolumn{1}{ |c }{}&
					\multicolumn{1}{ |c|  }{$(640,1280)$}&1.68e-07
					&4.18  
					&1.15e-07
					&3.81
					&4.27e-09
					&2.92
					&3.69e-09	
					&5.00
					
					\\
					\multicolumn{1}{ |c }{BDF3}&
					\multicolumn{1}{ |c|  }{$(1280,2560)$}&9.26e-09
					&2.82     
					&8.16e-09
					&2.81&5.62e-10
					&2.98&1.15e-10
					&5.00
					\\
					\multicolumn{1}{ |c }{}&
					\multicolumn{1}{ |c|  }{$(2560,5120)$}&1.31e-09  
					&&1.16e-09&&7.14e-11&&3.61e-12&
					\\
					\hline
					\hline
			\end{tabular}}
			\caption{Accuracy test for the 1D BGK equation. Initial data is given in \eqref{bgk 1d initial}.}\label{1d accuracy}
		\end{table}
	\end{center}
	
	{Since semi-Lagrangian schemes allow large CFL numbers, the choice of CFL number is important to secure both accuracy and efficiency. For such a reason we compute the the errors obtained with different values of CFL numbers for various Knudsen numbers. In Figure \ref{fig optimal CFL}, we report the relative $L^1$ errors obtained by $N_x=160,~320$ for each scheme. 
		\begin{figure}[htbp]
			\centering
			\begin{subfigure}[b]{0.45\linewidth}
				\includegraphics[width=1\linewidth]{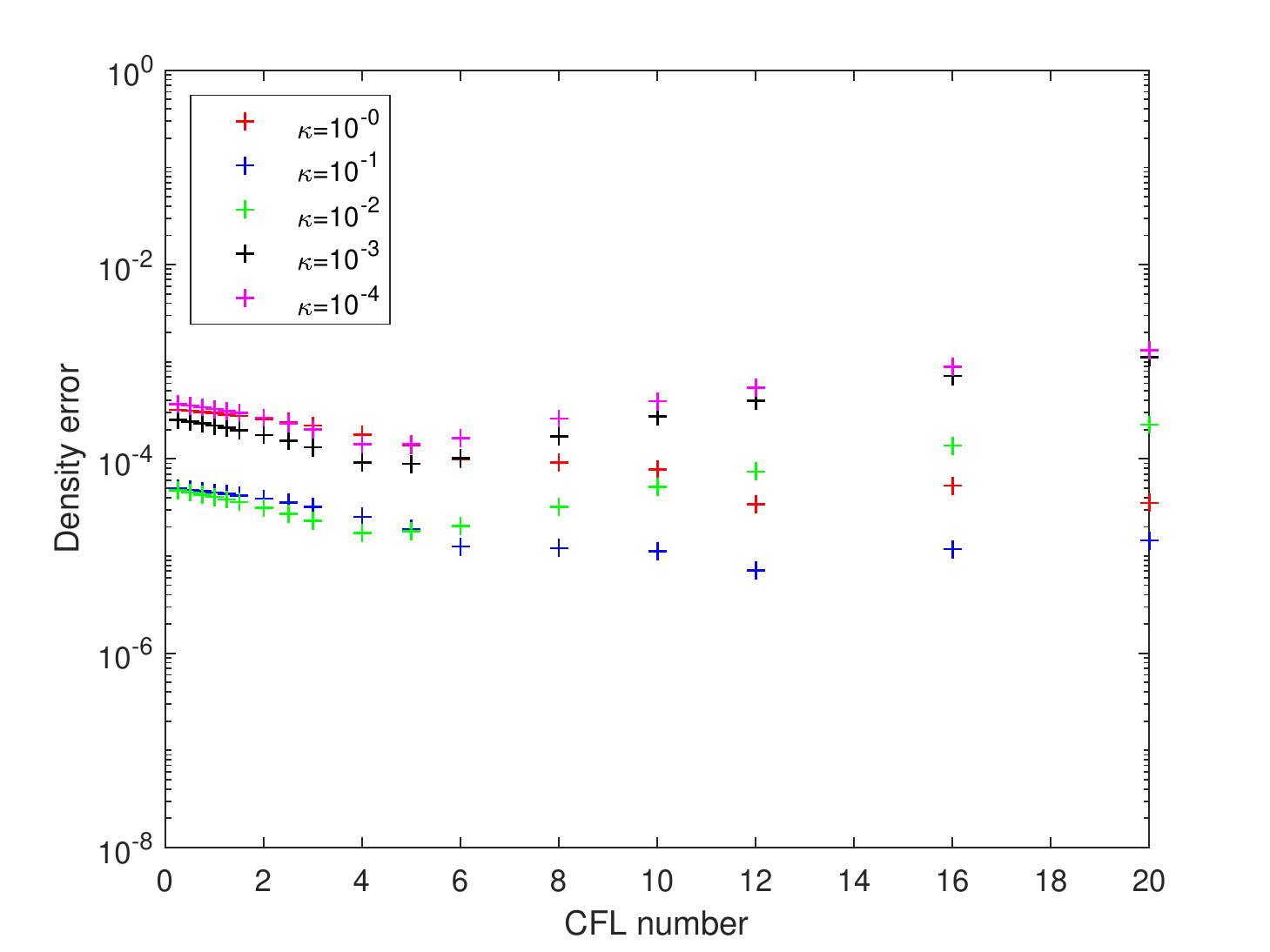}
				\subcaption{RK2}
			\end{subfigure}	
			\vspace*{5mm}
			\begin{subfigure}[b]{0.45\linewidth}
				\includegraphics[width=1\linewidth]{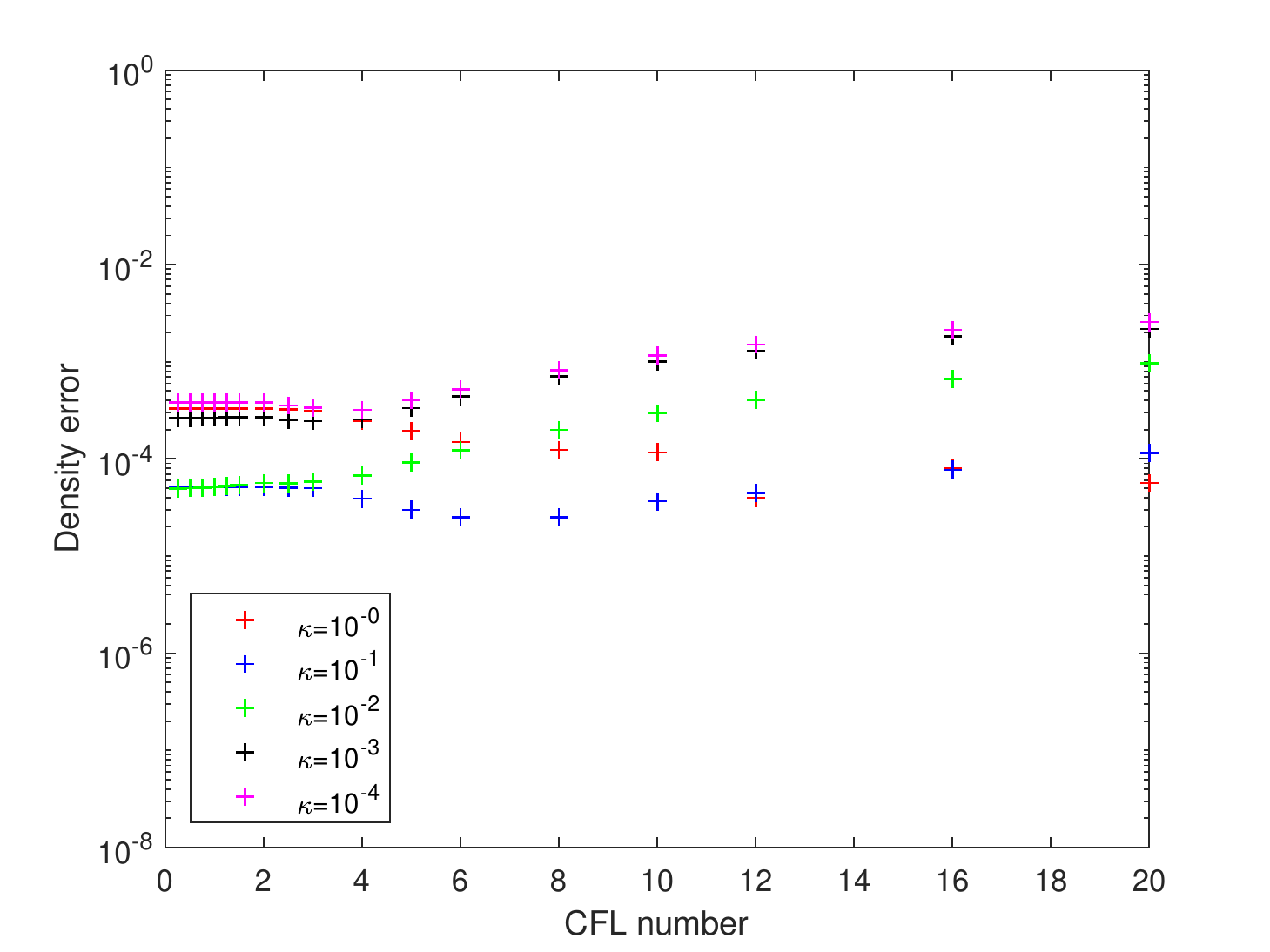}	
				\subcaption{BDF2}
			\end{subfigure}		
			\begin{subfigure}[b]{0.45\linewidth}
				\includegraphics[width=1\linewidth]{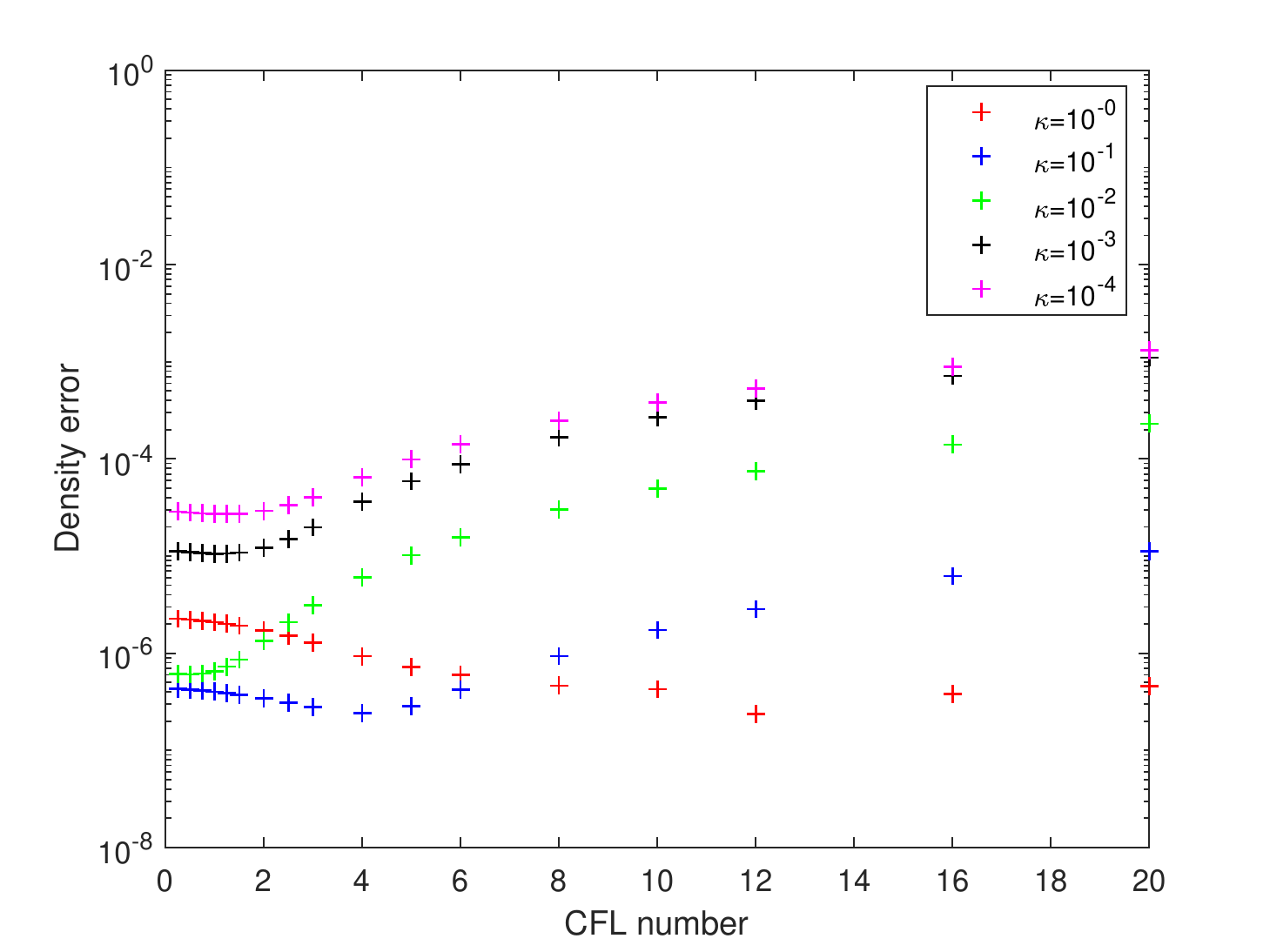}
				\subcaption{RK3}
			\end{subfigure}	
			\vspace*{5mm}
			\begin{subfigure}[b]{0.45\linewidth}
				\includegraphics[width=1\linewidth]{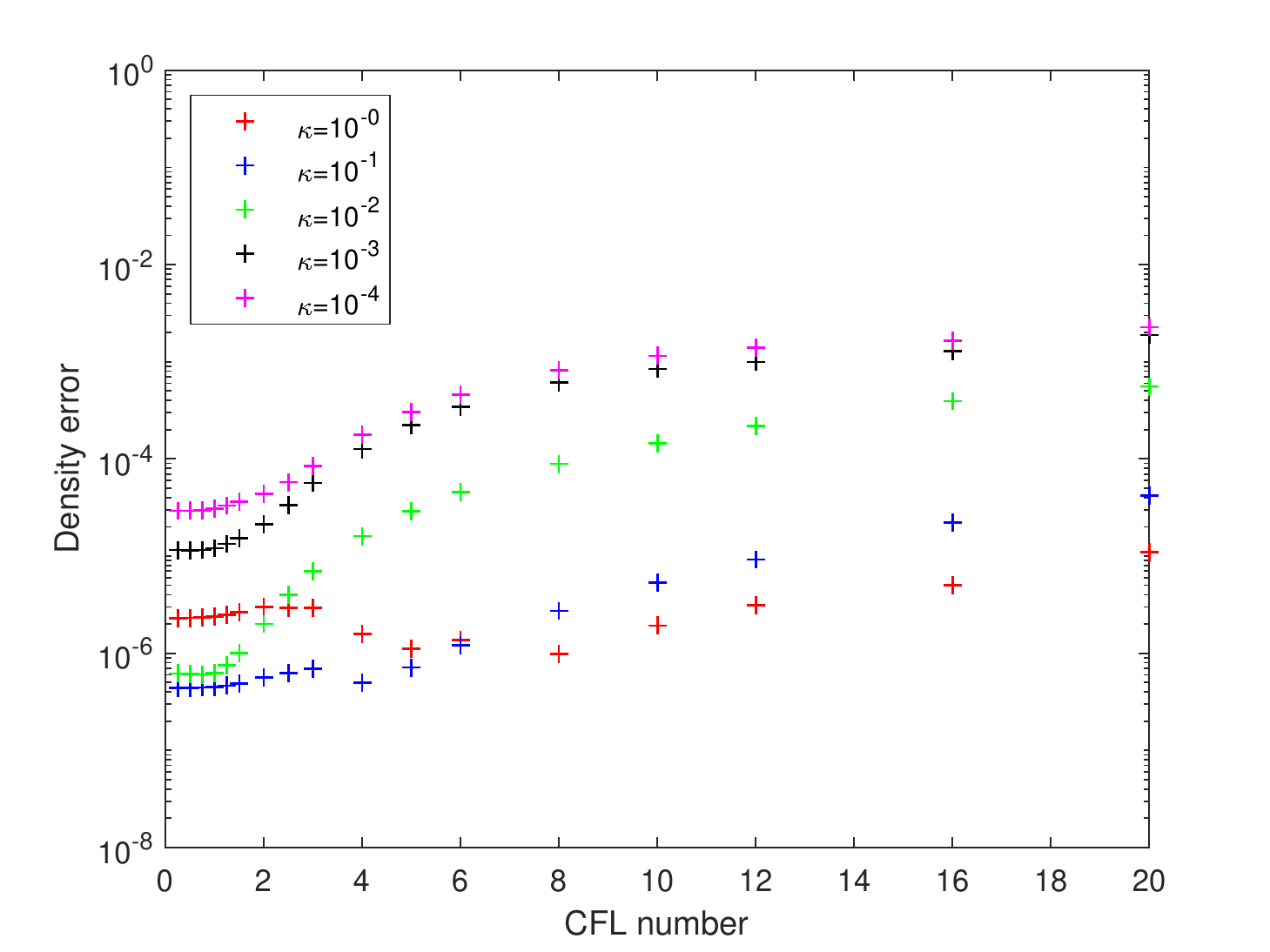}	
				\subcaption{BDF3}
			\end{subfigure}		
			\caption{Optimal CFL number. RK based schemes (left) and BDF based schemes (right). Initial data is given in \eqref{bgk 1d initial}.}\label{fig optimal CFL}  	
		\end{figure}
		From Figure \ref{fig optimal CFL}, we observe that the error has a non monotone behaviour with respect to the CFL number, so there is an optimal CFL number that depends on the Knudsen number, and is in general larger than one. For large enough Courant number the error increases, as expected. Also, we can see that higher order schemes give smaller errors. 
		We note that RK based schemes have smaller errors compared to the BDF based schemes of same order. Furthermore, compared to the case of BDF based schemes, in RK based schemes the optimal CFL number for which an error becomes the smallest one appears for larger CFL numbers. This can be explained by the smaller error constants of RK based schemes.}

	\subsubsection{1D Riemann problem}\label{sec 1d riemann}
	Here we consider the classical Riemann problem to check the shock capturing capability of of high order conservative semi-Lagrangian schemes in the fluid limit $\kappa \rightarrow 0$. We compare our scheme with non-conservative semi-Lagrangian schemes \cite{GRS}  and with conservative correction based semi-Lagrangian methods \cite{BCRY}. 
	As initial data we consider local Maxwellians with piecewise constant macroscopic quantities used in \cite{GRS}:
	\begin{align}\label{initial 1d riemann}
	(\rho_0,u_0,p_0) = \left\{\begin{array}{lr}
	(2.25,0,1.125), & \text{for } x\le 0.5\\
	(3/7,0,1/6), & \text{for } x>0.5\\
	\end{array}\right\}.
	\end{align}
	We impose freeflow boundary condition on the interval $ x \in [0, 1]$, and velocity domain $ v \in [-10, 10]$ upto final time $T_f=0.16$. We take uniform grids with $N_x=200$, $N_v=30$ and CFL$=2$, which means a time step by $\Delta t = \Delta x/5.$
	For this test, we use Eq.~\eqref{second dm}, described in Section \ref{sec:CM}, to construct the local discrete Maxwellian, which we prefer to the approach based on $L^2$ minimization,  because it guarantees the positivity of local Maxwellian with relatively small number of velocity grid points $N_v=30$, for which the approach in \eqref{first dm} may produce negative values, thus leading to instability. In Figure \ref{Sod fig}, the results obtained with BDF3 based schemes are plotted for comparison. The results show that our schemes are comparable to the high order conservative method introduced in \cite{BCRY} in capturing the exact shock position, while non conservative semi-lagrangian schemes produced shocks with incorrect speed.

	\begin{figure}[htbp]
		\centering
		\begin{subfigure}[b]{0.45\linewidth}
			\includegraphics[width=1\linewidth]{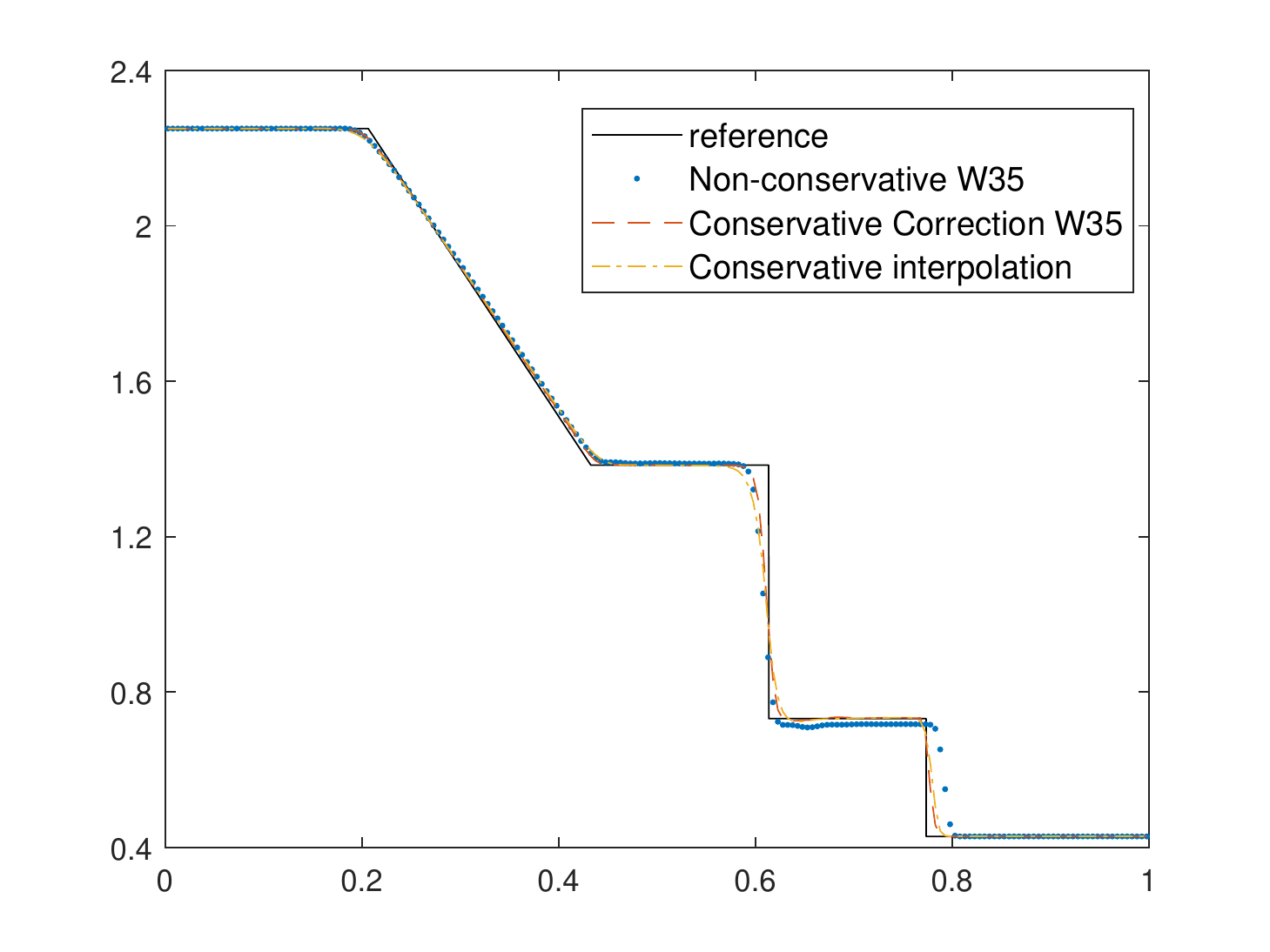}
			\subcaption{Density}
		\end{subfigure}	
		\vspace*{5mm}
		\begin{subfigure}[b]{0.45\linewidth}
			\includegraphics[width=1\linewidth]{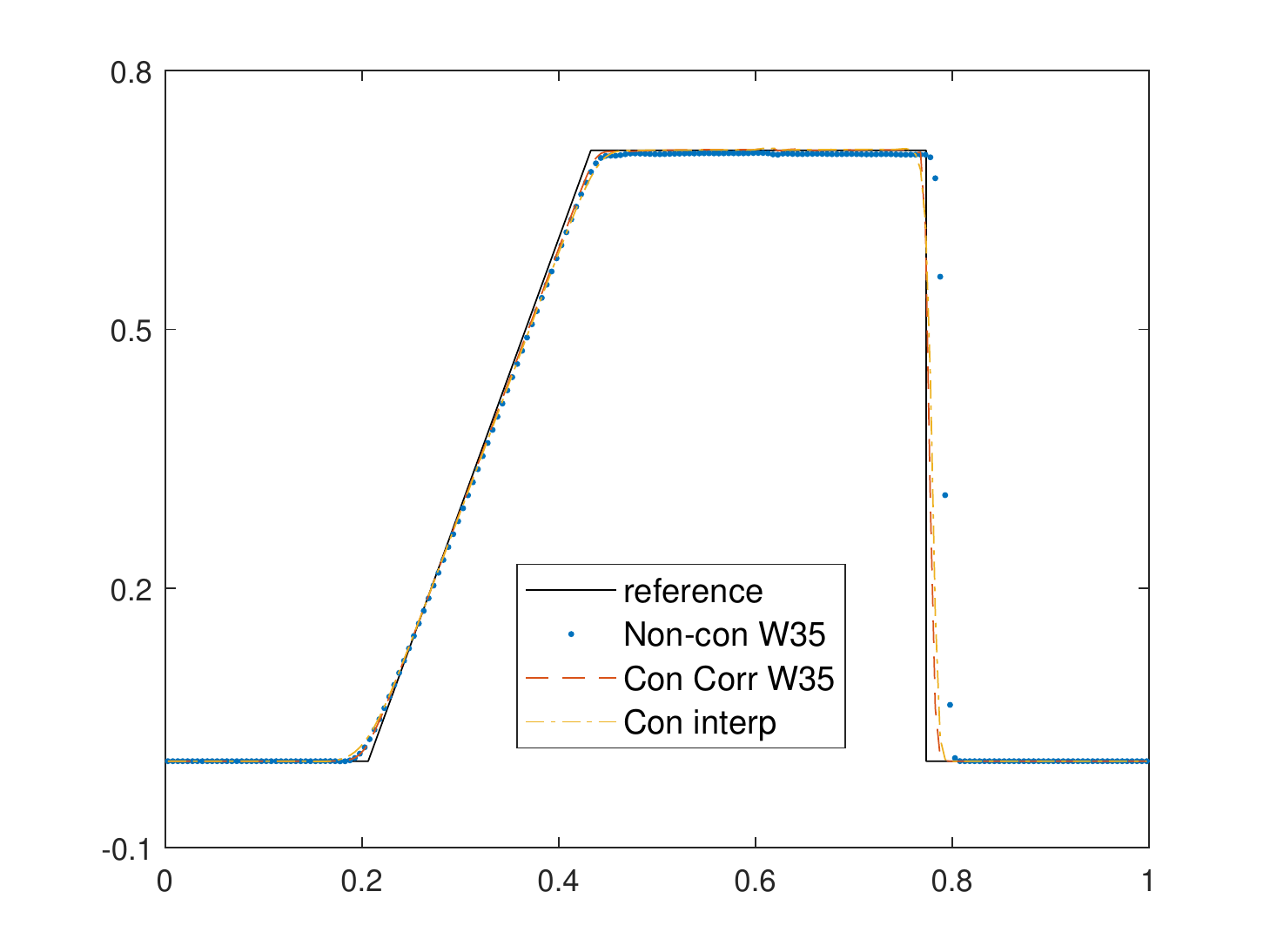}	
			\subcaption{Velocity}
		\end{subfigure}	
		
		\begin{subfigure}[b]{0.45\linewidth}
			\includegraphics[width=1\linewidth]{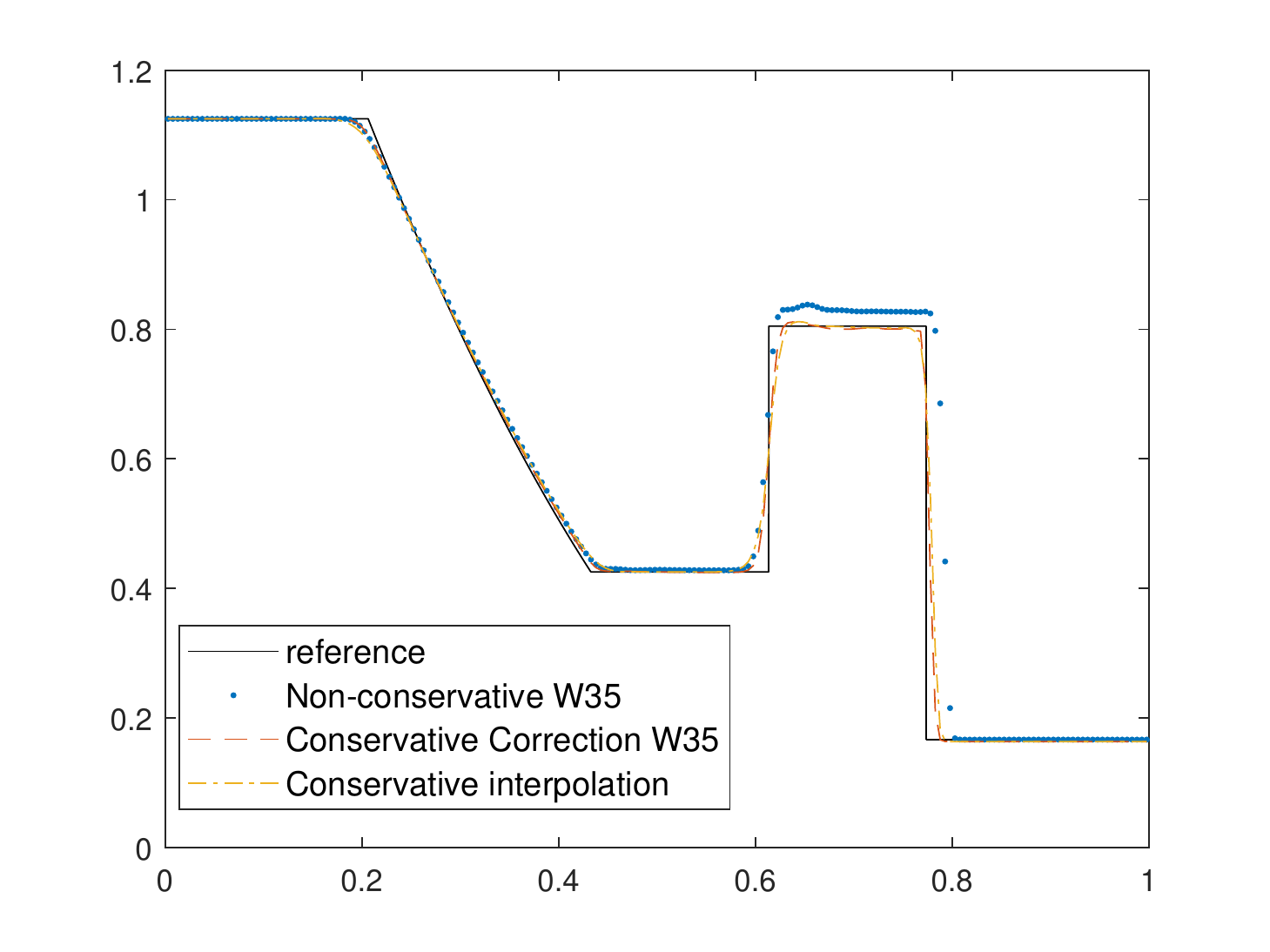}	
			\subcaption{Temperature}
		\end{subfigure}		
		\begin{subfigure}[b]{0.45\linewidth}
			\includegraphics[width=1\linewidth]{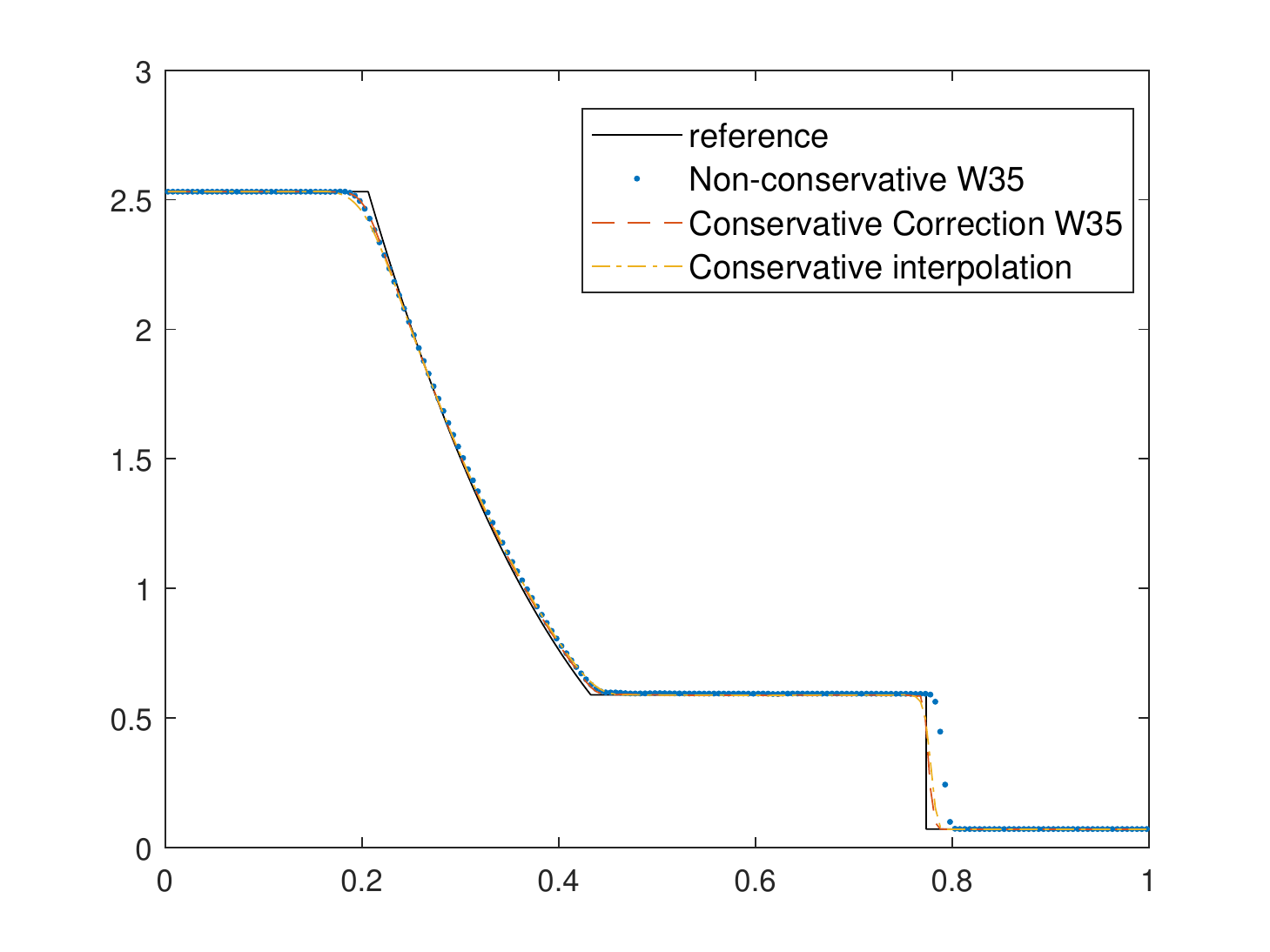}	
			\subcaption{Pressure}
		\end{subfigure}	
		
		\caption{Riemann problem for 1D BGK model. BDF3 based method with QCWENO35. Initial data is given in \eqref{initial 1d riemann}.}\label{Sod fig}  	
	\end{figure}
	
	%

	\subsection{2D BGK model}
	In 2D problem $(d_x, d_v)=(2,2)$, we take uniform meshes both in space ($\Delta x=\Delta y=h$) and in velocity ($\Delta v_x = \Delta v_y = \Delta v$). For the BGK model, we use the CFL number defined by 
	\begin{align}\label{CFL bgk 2d}
	\text{CFL}\displaystyle :=\max_{1\leq \ell \leq d_v} \max_j|v_j^{(\ell)}|\frac{\Delta t}{\Delta x}.
	\end{align}
	
	\subsubsection{Accuracy test}
	Here, we check the accuracy of semi-Lagrangian methods for 2D BGK model. Initial conditions are taken by the local Maxwellian with the macroscopic quantities:  
	\begin{align}\label{bgk accuracy 2d}
	&\rho_0(x,y)= 1 + (\sin{\pi x})^2 (\sin{\pi y})^2, \quad u_0(x,y) = (0,0), \quad T_0(x,y) = 1, 
	\end{align}
	where periodic boundary conditions are imposed on $ [0,1]^2$. We take $N_{v}=32$ for the velocity domain $[-10,10]^2$, which means that we use $(32+1)^2$ grid points. The final time is  $T^f=0.6$ up to which the solution remains smooth. We use a fixed time step $\Delta t$ with CFL$=4$ using \eqref{CFL bgk 2d}. In this test, the positivity of Maxwellian is not strictly demanded, hence, we compute the Maxwellian using \eqref{first dm}. Here we consider RK2 and BDF2 based semi-Lagrangian schemes combined with the conservative reconstruction based on the 2D CWENO23 method, so the expected convergence rates are between 2 and 3. Note that 2D CWENO23 is a genuinely 2D reconstruction, so we do not consider dimension by dimension interpolation here.
	\begin{center}
		\begin{table}[htbp]
			\centering	
			{\begin{tabular}{|cccccccccc|}
					\hline
					\multicolumn{10}{ |c| }{Relative $L^1$ error and order of density} \\ \hline
					\multicolumn{1}{ |c }{}&
					\multicolumn{1}{ |c| }{}& \multicolumn{2}{ c  }{$\kappa=10^{-6}$} & \multicolumn{2}{ |c }{$\kappa=10^{-4}$}& \multicolumn{2}{ |c| }{$\kappa=10^{-2}$} &
					\multicolumn{2}{ |c| }{$\kappa=10^{-0}$} \\ \hline
					\multicolumn{1}{ |c }{}&
					\multicolumn{1}{ |c|  }{$(N_x^2,(2N_x)^2)$} &
					\multicolumn{1}{ c  }{error} &
					\multicolumn{1}{ c|  }{rate} &
					\multicolumn{1}{ |c  }{error} &
					\multicolumn{1}{ c  }{rate} &
					\multicolumn{1}{ |c  }{error} &
					\multicolumn{1}{ c|  }{rate} &
					\multicolumn{1}{ |c  }{error} &
					\multicolumn{1}{ c|  }{rate}     \\ 
					\hline
					\hline	
					\multicolumn{1}{ |c }{}&
					\multicolumn{1}{ |c|  }{$(40^2,80^2)$}&5.16e-04
					
					&2.54     
					&4.89e-04
					&2.59&8.56e-05
					&2.54&1.74e-05
					&3.03
					\\
					\multicolumn{1}{ |c }{RK2}&
					\multicolumn{1}{ |c|  }{$(80^2,160^2)$}&8.90e-05
					&2.64         
					&8.14e-05
					&2.71&1.47e-05
					&2.29&2.12e-06
					&3.03
					\\
					\multicolumn{1}{ |c }{}&
					\multicolumn{1}{ |c|  }{$(160^2,320^2)$}&1.43e-05   
					&&1.25e-05&&3.01e-06&&2.60e-07&
					\\
					\hline
					\hline	
					\multicolumn{1}{ |c }{}&
					\multicolumn{1}{ |c|  }{$(40^2,80^2)$}&1.31e-03
					&2.35      
					&1.28e-03
					&2.36&5.00e-04
					&2.22&2.42e-05	
					&3.07
					\\
					\multicolumn{1}{ |c }{BDF2}&
					\multicolumn{1}{ |c|  }{$(80^2,160^2)$}&2.57e-04
					&2.42        
					&2.50e-04
					&2.42&1.08e-04
					&2.15&2.89e-06
					&3.08
					\\
					\multicolumn{1}{ |c }{}&
					\multicolumn{1}{ |c|  }{$(160^2,320^2)$}&4.83e-05  
					&&4.68e-05&&2.42e-05&&3.41e-07&
					\\
					\hline
					\hline	
					\multicolumn{1}{ |c }{}&
					\multicolumn{1}{ |c|  }{$(40^2,80^2)$}&5.39e-04
					&2.53      
					&5.13e-04
					&2.58&8.71e-05
					&2.70&1.74e-05
					&3.03
					\\
					\multicolumn{1}{ |c }{RK3}&
					\multicolumn{1}{ |c|  }{$(80^2,160^2)$}&9.31e-05
					&2.63         
					&8.57e-05
					&2.70&1.34e-05
					&2.66&2.14e-06
					&3.02
					\\
					\multicolumn{1}{ |c }{}&
					\multicolumn{1}{ |c|  }{$(160^2,320^2)$}&1.50e-05   
					&&1.32e-05
					&&2.12e-06&&2.64e-07&
					\\
					\hline
					\hline
					\multicolumn{1}{ |c }{}&
					\multicolumn{1}{ |c|  }{$(40^2,80^2)$}&7.18e-04
					&2.34     
					&6.94e-04
					&2.37&1.35e-04
					&2.68&2.74e-05	
					&3.03
					\\
					\multicolumn{1}{ |c }{BDF3}&
					\multicolumn{1}{ |c|  }{$(80^2,160^2)$}&1.42e-04
					&2.92     
					&1.35e-04
					&2.94&2.11e-05
					&2.89&3.36e-06
					&3.02
					\\
					\multicolumn{1}{ |c }{}&
					\multicolumn{1}{ |c|  }{$(160^2,320^2)$}&1.87e-05   
					&&1.75e-05&&2.84e-06&&4.15e-07&
					\\
					\hline
			\end{tabular}}
			\caption{Accuracy test for the 2D BGK equation. Initial data is given in \eqref{bgk accuracy 2d}.}\label{2d accuracy}
		\end{table}
	\end{center}
	In Table \ref{2d accuracy}, we report the numerical errors and the corresponding accuracy orders for density. The result shows the desired convergence rate between 2 and 3 for all values of Knudsen number. Notice that for large values of the Knudsen number most error is due to interpolation, so there is no improvement in using high order time discretization. Note that in case of the RK3 based method, we observe order reduction in the limit $\kappa=10^{-6}$ as in the 1D accuracy test (see Table \eqref{1d accuracy}). 
	A similar result can be found, for example, in \cite{GRS}.

	\subsubsection{2D Riemann problem}\label{sec 2d riemann}
	Now, we perform a test on the 2D Riemann problem presented in \cite{dimarco2013towards}. We use $\kappa = 10^{-4}$ so that the solution is rather close to a local Maxwellian.
	The initial distribution function is a local Maxwellian with macroscopic moments:
	\[
	(\rho_0,{u_1}_0,{u_2}_0,T_0) = \left\{\begin{array}{lr}
	(1,0,0,5), & \text{for } (x-1)^2+(y-1)^2 \le 0.2\\
	(0.125,0,0,4), & \text{for } (x-1)^2+(y-1)^2 > 0.2\\
	\end{array}\right\}.
	\]
	The computation is performed on the space domain $[0, 2]^2$ with free flow boundary conditions; the velocity domain is $[-15, 15]^2$ , and the final time $T_f=0.07$. We take uniform grids with $N_x=200$, $N_y=200$, and $N_v=20$ and fix a time step $\displaystyle \Delta t$ with CFL$=2.1$ using \eqref{CFL bgk 2d}. Here we construct Maxwellian with the approach in \eqref{first dm} ($L^2$-projection method).

	\begin{figure}[htbp]
		\centering
		\begin{subfigure}[b]{0.40\linewidth}
			\includegraphics[width=1\linewidth]{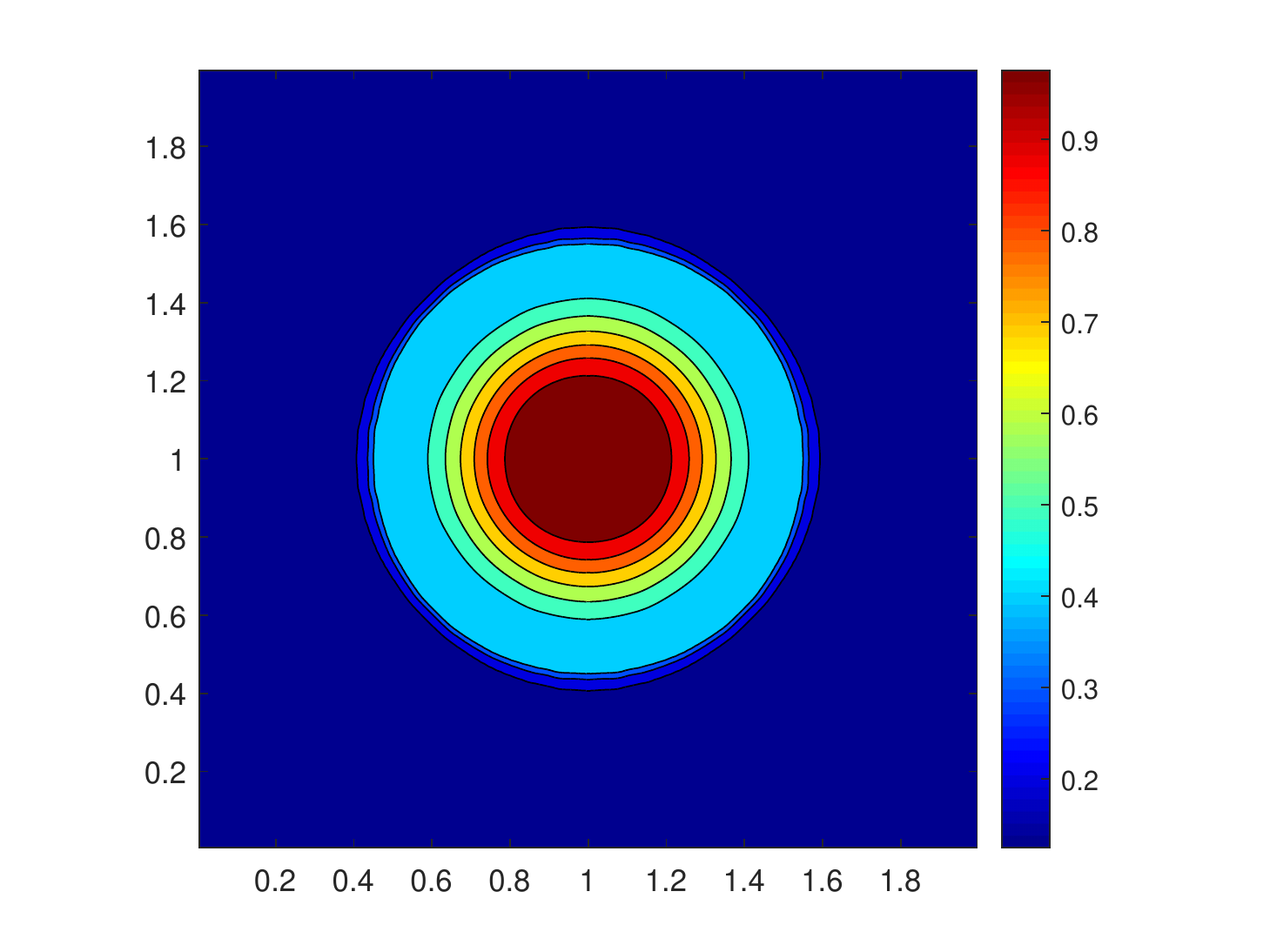}
			\subcaption{Density}
		\end{subfigure}	
		\begin{subfigure}[b]{0.40\linewidth}
			\includegraphics[width=1\linewidth]{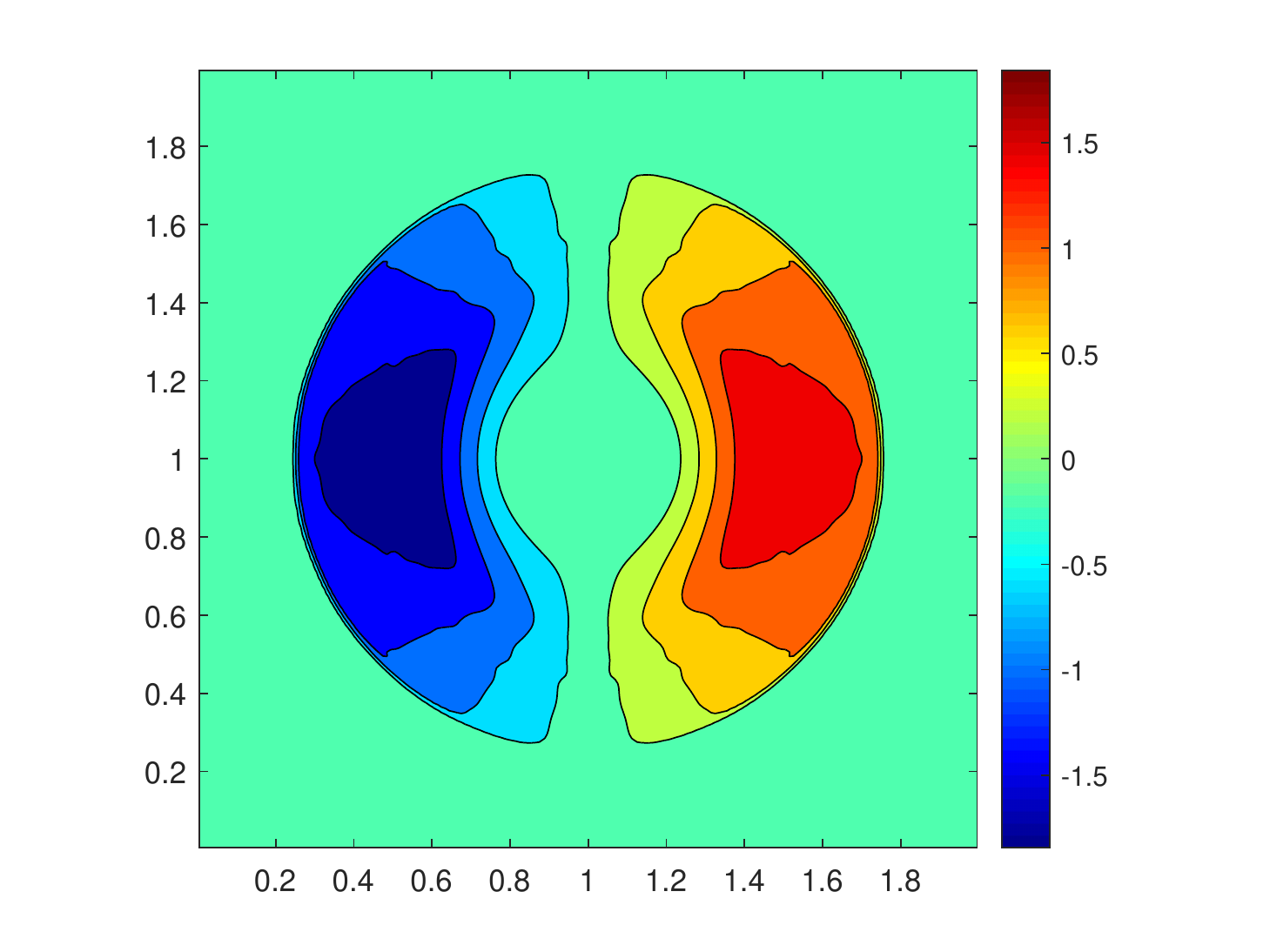}
			\subcaption{Velocity in x-direction}
		\end{subfigure}	
		
		\begin{subfigure}[b]{0.40\linewidth}
			\includegraphics[width=1\linewidth]{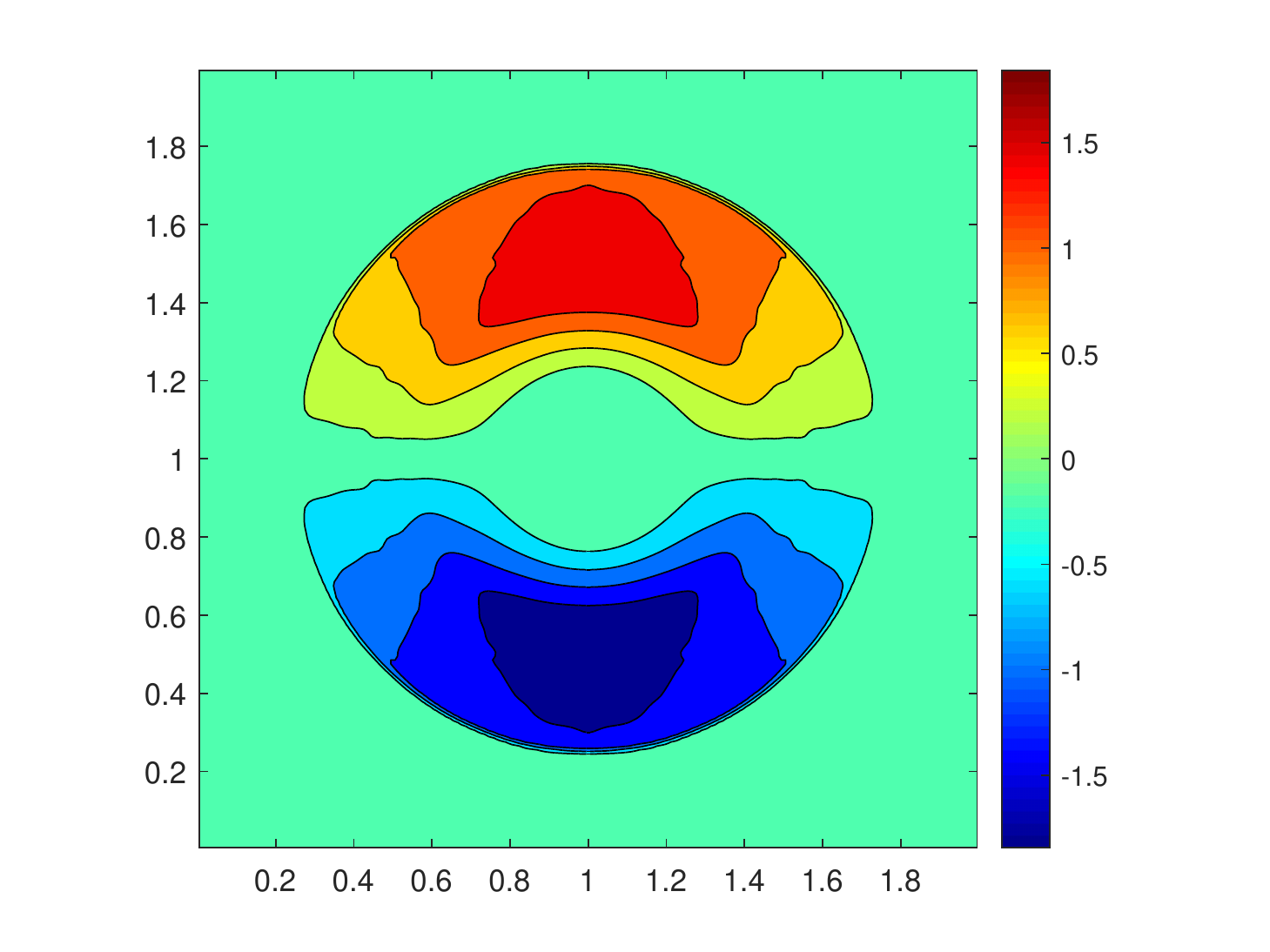}	
			\subcaption{Velocity in y-direction}
		\end{subfigure}		
		\begin{subfigure}[b]{0.40\linewidth}
			\includegraphics[width=1\linewidth]{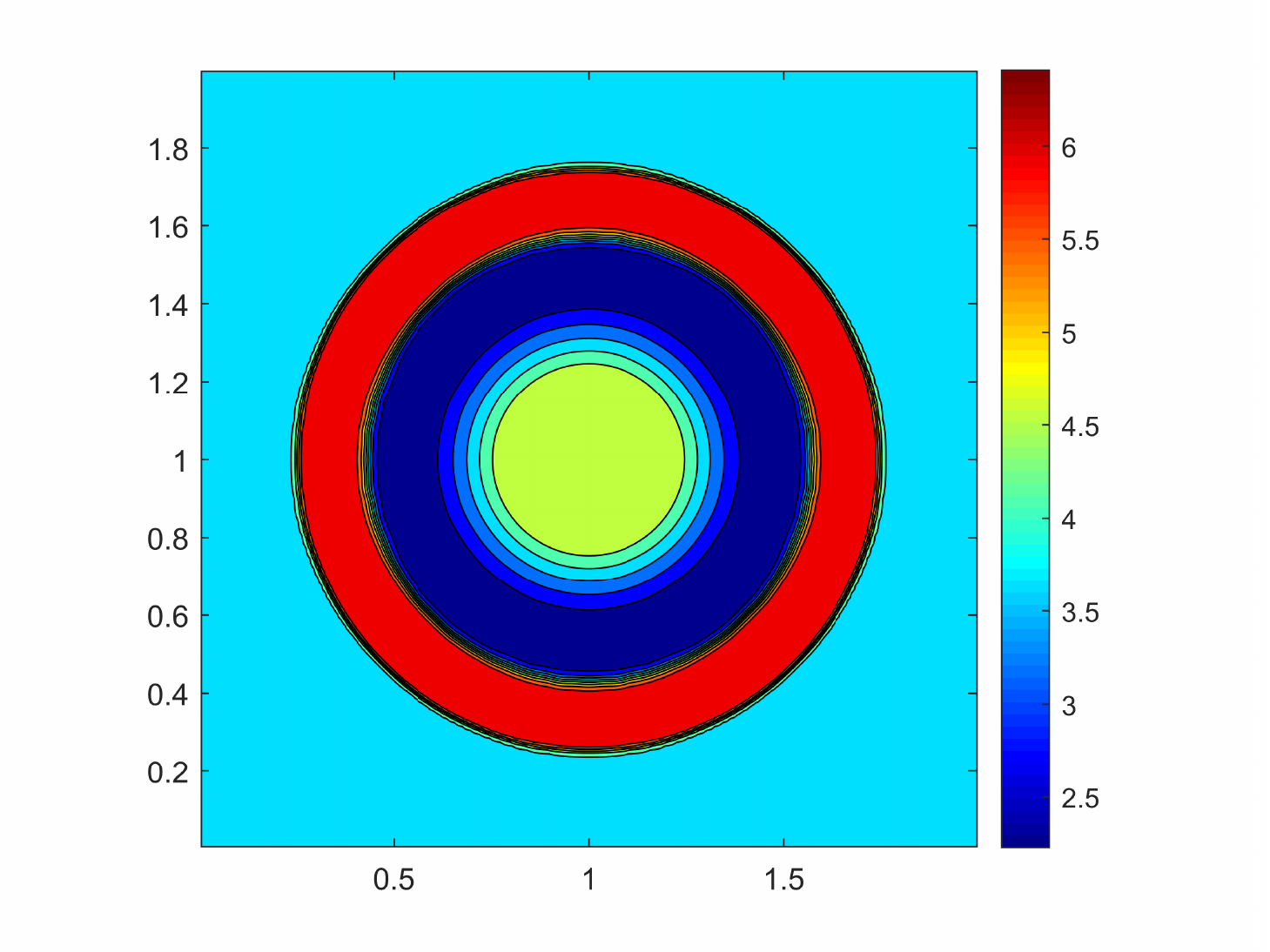}	
			\subcaption{Temperature}
		\end{subfigure}	
		\caption{2D Riemann problem with $\kappa=10^{-4}$, Macroscopic moments. } \label{Riemann 2d} 	
		
		{\small \begin{subfigure}[b]{0.40\linewidth}
				\includegraphics[width=1\linewidth]{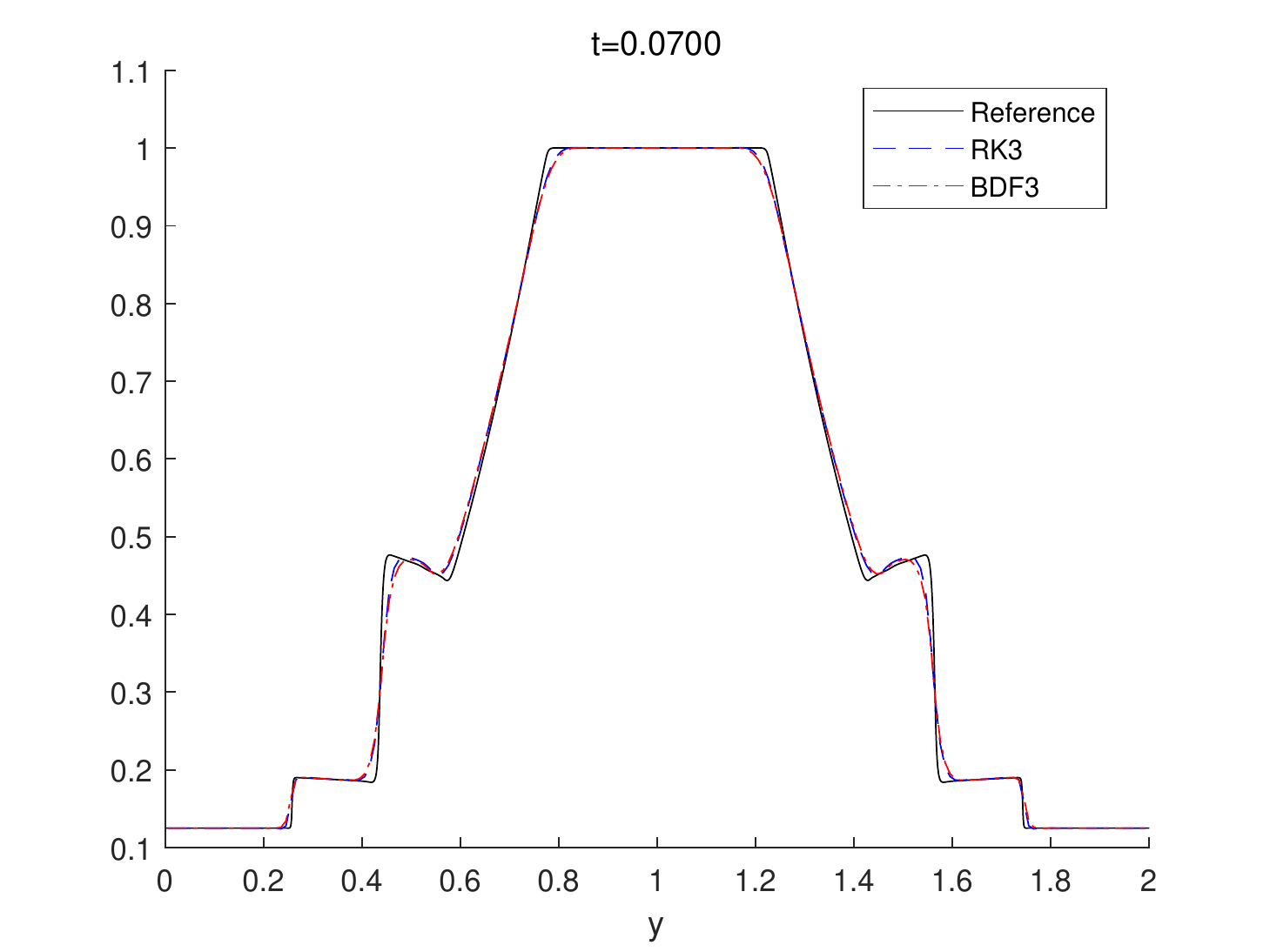}
				\subcaption{Density}
			\end{subfigure}	
			
			\begin{subfigure}[b]{0.40\linewidth}
				\includegraphics[width=1\linewidth]{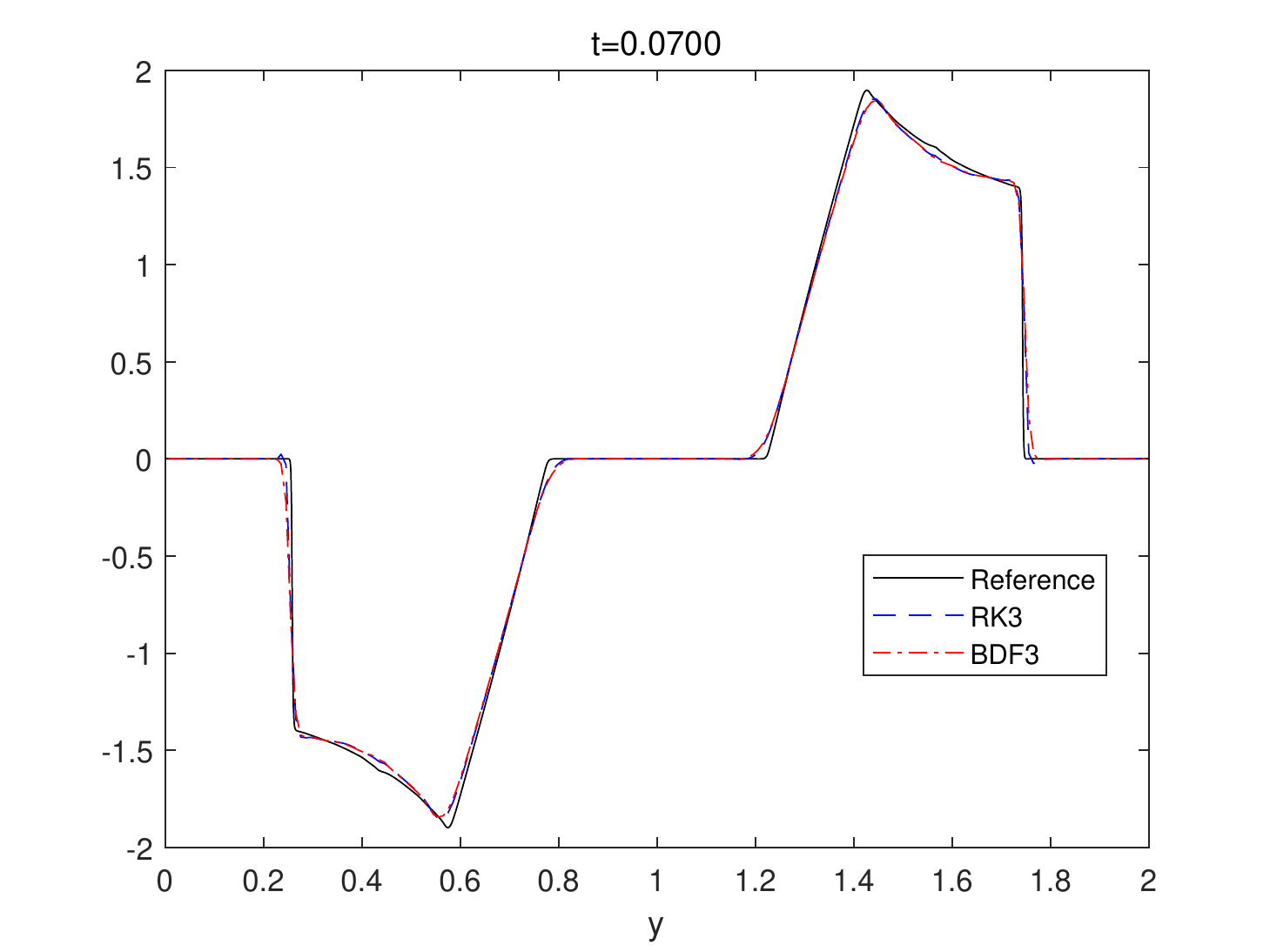}	
				\subcaption{Velocity in y-direction}
			\end{subfigure}		
			\begin{subfigure}[b]{0.40\linewidth}
				\includegraphics[width=1\linewidth]{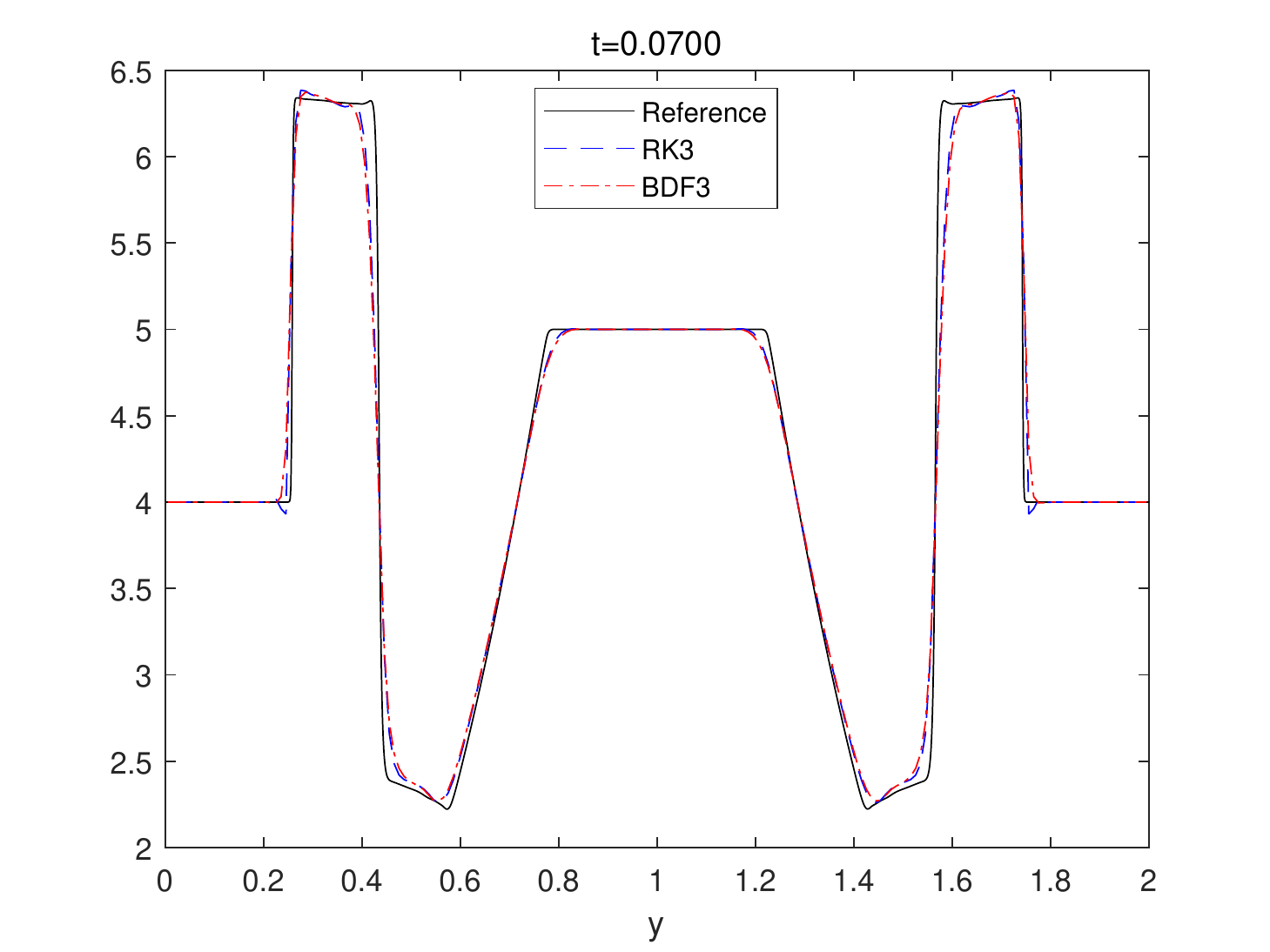}	
				\subcaption{Temperature}
		\end{subfigure}}	
		\caption{Riemann problem for 2D BGK model with $\displaystyle \kappa=10^{-4}$. Macroscopic moments at $x=x_{100}$}\label{Riemann 2d slice}
	\end{figure}
	
	In Figure \ref{Riemann 2d}, we present a contour plot of the macroscopic variables obtained by the BDF3-QCWENO23 method. In the following Figure \ref{Riemann 2d slice}, we compare the results obtained by BDF3-QCWENO23 and 
	RK3-QCWENO23 schemes with a reference solution obtained by solving compressible 2D Euler equations by explicit conservative finite difference method with WENO23 and RK3 (with a grid of $1600\times 1600$ points). For comparison, all solutions are plotted at the location $x=1$. In particular, we plot solution at $x=x_{100}$ which is the closest grid point to $x=1$.
	The result shows that high order conservative schemes are able to capture correct shock positions which appear near $y \approx 0.26, 1.74$, and reproduce the profile of the reference solution. We omit the figure for the slice of x-directional velocity because the x-directional velocity is zero at $x=1$, by symmetry.

	%
	%

\end{section}


\section{Conclusions}
The conservative reconstruction, $Q[u^n]$ derived and analysed in the first part of the this paper \cite{BCRY1} has been adopted to propose new conservative semi-Lagrangian schemes, which are then applied to a variety of problems, such as rigid rotation, Vlasov-Poisson system, and BGK model.

$Q$ is obtained by sliding average of a basic reconstruction, $R$, therefore it inherits the properties of $R$, such as, for example, non-oscillatory behaviour or positivity. 
The choice of $R$ depends on the properties one would like to maintain for $Q$. 

High order operator splitting is adopted in the case of rigid rotation and Vlasov-Poisson system, while no splitting is necessary to solve the BGK model. 

The resulting schemes are stable, i.e.\ they allow large time steps with no CFL restriction,
and are high order accurate both in space and time. 

In the case of the BGK model, strict conservation is obtained by combining the conservative reconstruction with a conservative treatment of the collision term (which is obtained either using a discrete Maxwellian or a least-square conservative projection). 
The implicit treatment of the collision term, together with the strict conservation properties of the method, guarantee asymptotic preserving property: as the Knudsen number vanishes, the method becomes a high resolution shock capturing scheme for the underlying Euler equations. 

We plan to extend the method to a wider range of problems, including BGK models for polyatomic gases and mixtures, and to analyse in greater detail the resulting limit schemes for compressible Navier-Stokes and Euler equations. Another important issue concerns adaptivity in velocity: regions with small values of the Knudsen number require fewer grid points in velocity. All these issues will be subject of future investigation.

\label{sec:conclusions}



%

\section*{Acknowledgments}
S. Y. Cho has been supported by ITN-ETN Horizon 2020 Project ModCompShock, Modeling and Computation on Shocks and Interfaces, Project Reference 642768. 
S.-B. Yun has been supported by Samsung Science and Technology Foundation under Project Number SSTF-BA1801-02. 
All the authors would like to thank the Italian Ministry of Instruction, University and Research (MIUR) to support this research with funds coming from PRIN Project 2017 (No. 2017KKJP4X entitled \emph{Innovative numerical methods for evolutionary partial differential equations and applications}). 
S. Boscarino has been supported by the University of Catania (\emph{Piano della Ricerca 2016/2018, Linea di intervento 2}). S. Boscarino and G. Russo are members of the INdAM Research group GNCS. 



	\bibliographystyle{amsplain}
	\bibliography{references}

\end{document}